\documentclass{article}
\usepackage{array}
\usepackage{amsfonts}
\usepackage{amsmath}
\usepackage{amssymb}
\usepackage{authblk}
\usepackage[colorlinks=true, linkcolor=blue, citecolor=blue, urlcolor=blue]{hyperref}
\usepackage{subcaption}
\usepackage{tikz}
\usepackage[normalem]{ulem}
\usepackage{color}
\usepackage[scr=boondox]{mathalpha}

\usepackage{algorithm}      
\usepackage{algpseudocode}  

\usepackage[square,numbers]{natbib}
\usepackage[left=2.5cm,top=2.5cm,right=2.5cm, bottom=2.5cm]{geometry}

\usepackage{graphicx}
\usepackage{enumerate}
\usepackage{soul} 
\usepackage[final, commentmarkup=none, authormarkup=none]{changes}
\definechangesauthor[name={Review 2}, color=blue]{R2}
\definechangesauthor[name={Review 1}, color=red ]{R1}

\newcommand{\bm}[1]{\boldsymbol{#1}}
\newcommand{\dd}{\mathrm{d}}
\newcommand{\pd}{\partial}
\newcommand{\Tr}{\intercal}

\newcommand{\fpd}[2]{\frac{\pd #1}{\pd #2}}

\usepackage{ifthen} 
\usepackage{cleveref} 

\newcommand{\figref}[2][]{\ifthenelse{\equal{#1}{}}{\ref{#2}}{\ref{#2}#1}}

\title{A ROM-based BDDC solver for unfitted p-FEM level-set-based \added[id=R2]{two-dimensional} lattice structures}
\author[1]{Gonzalo Bonilla Moreno \thanks{gonzalo.bonilla.moreno@gmail.com}}
\author[1]{Giuliano Guarino \thanks{giuliano.guarino@epfl.ch}}
\author[1]{Pablo Antolin \thanks{pablo.antolin@epfl.ch}}
\affil[1]{\small Institute of Mathematics, École polytechnique fédérale de Lausanne, Station 8, CH-1015 Lausanne, Switzerland}

\begin{document}

\maketitle

\begin{abstract}
We present a domain decomposition method for the fast simulation of large \added[id=R2]{two-dimensional} lattice structures described by level set functions. The method does not rely on homogenization or multiscale techniques, and therefore avoids their underlying assumptions such as scale separation and periodicity. Individual cells are defined through level set functions and mapped into physical space using arbitrary order mappings, which allows the creation of complex graded designs with varying geometries and topologies. The discretization is based on unfitted p-FEM, where each cell is approximated by a single high order element. This choice naturally handles the implicit geometric description and provides high accuracy with a moderate number of degrees of freedom. The solver is built on the Balanced Domain Decomposition by Constraints (BDDC) method, where each cell corresponds to one subdomain. To accelerate the assembly of the cell stiffness matrices, we combine a fast assembly technique that separates the contributions of the geometric mapping from the trimmed domain with a reduced order model (ROM) based on the matrix discrete empirical interpolation method (MDEIM). The ROM surrogate is trained offline and can be reused for any geometric mapping, restricting the expensive quadrature on cut elements to the training stage. A stabilization term is introduced to ensure the scalability of the solver when using the ROM approximation, at the cost of a small and controllable error. We validate the method through a series of numerical experiments and demonstrate its performance on a\deleted[id=R2]{ complex} 2D problem with more than 17,000 cells of varying geometry, which is solved in approximately 30 seconds on a standard laptop. The number of solver iterations grows only mildly as the number of subdomains increases, provided the ratio between subdomain and mesh sizes is kept constant, consistent with the scalability properties of BDDC methods.
\end{abstract}

{\small
\noindent\textbf{Note added in version 3.}
The version of record is the article published in \emph{Computer Methods in Applied Mechanics
and Engineering}, article 119304 (2026), \texttt{doi:10.1016/j.cma.2026.119304}. This version of
the preprint corrects three defects of the printed text. None of them changes a reported result.
\begin{enumerate}[(i)]
    \item In the first of Equations \eqref{eq: ELA constitutive}, the two inverse-Jacobian factors
    were transposed with respect to the form given by the chain rule stated above them.
    \item In Section \ref{sssec: RES wing}, the threshold parameters of the sandwich wing were
    attributed to $\mathbb{P}_1^\mathrm{SD}$, while the interval printed there is that of
    $\mathbb{P}_2^\mathrm{SD}$. The label has been corrected to match the interval.
    \item In Algorithm \ref{alg: Local solves a}, the matrix $\hat{T}^{(i)}$ was used without being
    defined, and the assembly step referred to step 5 instead of step 3. The reduced right-hand
    side, which shared its symbol with the pulled-back traction $\hat{t}$, has been renamed
    $\hat{r}^{(i)}$.
\end{enumerate}
\noindent A few typographical corrections are also included.
}

\section{Introduction}\label{sec: INTRO}

Over the last few decades, the cost of additive manufacturing has decreased while the quality of their products has improved~\cite{ngo2018additive}.
This, combined with the emergence of techniques for the fabrication of lattices at the nanoscale~\cite{vyatskikh2018additive}, has paved the way to the creation of architected materials with unprecedented strength to weight ratio, heat exchange capacity, among others.
While porous geometries are abundant in nature (e.g., bones and bird beaks), artificial porous (heterogeneous) artifacts were very difficult to create prior to the additive manufacturing era.
Now, their cell topology, geometry, and material can be tailored to achieve specific performance requirements~\cite{pan2020design}: significant weight savings while maintaining the stiffness and strength of homogeneous structures~\cite{zheng2014ultralight,bauer2014high,han2015new,meza2015resilient,shaikeea2022Toughness}.
In addition, lattices can behave mechanically in very atypical ways (e.g., highly stretchable, auxetic~\cite{ren2018auxetic}) and can be multi-functional, which makes them attractive also for applications such as energy absorption~\cite{tancogne2016additively} and storage~\cite{shan2015multistable}, vibration reduction, thermal management, etc.~\cite{surjadi2019mechanical}.

From a simulation perspective, lattice structures presenting a large number of cells are quite challenging~\cite{kochmann2019multiscale}.
A straightforward application of finite element methods is greedy in terms of computing resources (both CPU time and memory).
To circumvent this problem, different numerical approaches are applied in practice, namely:
multiscale FEM~\cite{buck2013multiscale,hou1997multiscale,castelletto2017multiscale}, Generalized and Extended FEM~\cite{dovskavr2021microstructure,savvas2014homogenization,strouboulis2000design}, multilevel FEM (FE$^2$)~\cite{feyel2003multilevel,schroder2014numerical}, or numerical homogenization~\cite{charalambakis2010homogenization,moulinec1995fft,zeman2010accelerating,muller2015homogenization,glaesener2020continuum},
where the macroscopic behavior of the heterogeneous materials is characterized through numerical simulations on representative volume elements.
However, to alleviate undesirable size effects~\cite{onck2001size,yoder2018size}, all these techniques rely on i) the separation of scales, which is usually not the case in practice due to the achievable length scale of current 3D printers; and ii) periodic cells, an assumption which does not hold for graded lattices.
In addition, macro-geometries are often slender geometries (plate, shell) that present just a few cells along the thickness direction.
We refer the interested reader to~\cite{kochmann2019multiscale} for an in-depth discussion of such approaches in the context of hierarchical metamaterials.
In the case of truss-based lattices, an appealing alternative is the use of 3D beam models~\cite{portela2018impact,jamshidian2020multiscale,weeger2022isogeometric}. However, those methods are no longer suitable in the case of thick trusses, and cannot be applied to cellular lattices.

On the other hand, full scale 3D finite element simulations remain rare due to their high computational cost.
They are typically limited to a few cells, or even a single one, and are used to estimate the macro-behavior properties of unit micro-structures or small samples~\cite{tancogne2016additively,bonatti2017large,kochmann2019multiscale,nguyen2025efficient}.
For instance, in~\cite{tancogne2016additively} the authors performed a nonlinear simulation using Abaqus, including plasticity effects, with $7\times7\times7$ octet cells and, even though no computing times were reported, they required a cluster computer with 120 cores.
Or~\cite{korshunova2021image,korshunova2021bending}, where the authors performed a high-fidelity simulation of 24 CT-scanned octets (98M unknowns) using the finite cell method~\cite{duster2008finite}, which required 52 minutes using 1120 processors of a supercomputer.

Recent works aim to overcome these limitations using full scale finite element methods that leverage state of the art domain decomposition techniques combined with reduced order modeling (ROM) that exploits cell similarities to accelerate the solution.
In \cite{hirschler2022fast,hirschler2024reduced} the authors proposed an accelerated inexact FETI-DP preconditioner that allows to analyze linear problems with thousands of cells (millions of degrees of freedom) in just a few minutes, and with very low memory requirements, using an off-the-shelf laptop. 
This work has been recently extended in~\cite{guillet2026efficient} to the case of nonlinear hyperelastic materials undergoing large deformations.
In \cite{guillet2025multilevel} similar ideas are exploited in a matrix-free multigrid solver able to simulate designs with hundreds of thousands of cells (billions of degrees of freedom) in just a few minutes using a few thousand processors from a supercomputer. 

However, the application of this family of fast methods is limited to cell geometries simple enough to be described with conforming discretizations and parametric morphings, which is essential for exploiting cell similarities through ROM techniques.
This prevents its applicability to the case of cellular like structures, with arbitrarily varying geometries and topologies that cannot be parameterized using such mappings.
This is the case of, for instance, cellular structures~\cite{han2015new,kumar2020inverse} or Triply Periodic Minimal Surfaces (TPMS)~\cite{feng2022triply}.

A natural way of handling such complex geometries is through level-set functions, which turn the geometric parameterization from explicit (parametric) to implicit.
Such representations are particularly well suited for unfitted discretizations, in which the geometry is embedded in a background grid that serves as base for discretizing the PDE at hand~\cite{duster2008finite,badia2018aggregated,main2018shifted,burman2025cut}.
In this way, all cell geometries share the same discretization grid, regardless of their shape or topology.

This common grid structure enables the use of projection based ROM techniques.
In particular, it allows the construction of surrogate models for the fast assembly of cell stiffness matrices using the matrix version of the discrete empirical interpolation method (MDEIM)~\cite{chaturantabut2010nonlinear,negri2015}. 
This bypasses the expensive tailored quadrature rules for cut elements required by unfitted methods, restricting their use to an offline training stage.
The combination of MDEIM with unfitted methods has been previously addressed in~\cite{karatzas2020projection,chasapi2023localized,mueller2026reduced}.

In this work, we propose a novel fast domain decomposition method for the linear elasticity problem on \added[id=R2]{two-dimensional} lattice structures described by level-set functions.
The method relies on unfitted p-FEM~\cite{parvizian2007finite,lehrenfeld2018analysis,martorell2024high} (see also~\cite{idesman2025oltem} for very high order unfitted discretizations of 3D elasticity; the motivation for high order discretizations will be addressed in Section~\ref{sec:modeling}) and exploits cell similarities through the ROM techniques described above.
The domain decomposition solver is based on the Balanced Domain Decomposition by Constraints method (BDDC), \added[id=R2]{introduced in \cite{dohrmann2003preconditioner,mandel2003convergence,li2006feti} for FEM and subsequently adapted to various numerical methods, including isogeometrical analysis \cite{beirao2012bddc,beirao2014isogeometric} and unfitted methods \cite{badia2018robust}.}
\deleted[id=R2]{previously studied in~\cite{badia2018robust} in the context of unfitted methods}. \added[id=R2]{Here, the BDDC preconditioner is} combined with the acceleration ideas first introduced in~\cite{hirschler2022fast,hirschler2024reduced} for \added[id=R2]{two-dimensional} lattice structures.
\deleted[id=R2]{As will be shown in Section~\ref{sec: RES}, the resulting method is able to solve 2D problems with more than 17,000 cells in around 30 seconds on a laptop. Moreover, the number of solver iterations remains asymptotically bounded as the number of subdomains grows, provided the subdomain-to-mesh-size ratio $H/h$ is kept constant, in agreement with the scalability properties of BDDC methods.}

The remainder of the paper is organized as follows.
Section~\ref{sec:modeling} introduces the geometric modeling of single unit cells as well as their assembly into full lattice structures. The elasticity problem and its discretization by means of unfitted $p$-FEM are also discussed.
In Section~\ref{sec: DD} the BDDC domain decomposition method considered is introduced, alongside the required hypotheses and specific characteristics for this problem.
The fast assembly of the cell stiffness matrices is discussed in Section~\ref{sec: SPEED}.
Finally, a series of numerical experiments that validate the considered hypotheses and assess the numerical performance of the proposed solver are presented in Section~\ref{sec: RES}, while main conclusions are drawn in Section~\ref{sec: CONC}.

\section{Modeling and Discretization} \label{sec:modeling}

\subsection{Geometric modeling}\label{ssec:geometry}
This section describes the geometric construction of the lattice cells investigated in this work. The geometry is obtained by juxtaposing cells, each constructed using a two-step procedure: first, mapping a parametric domain, and second, trimming it with the desired porous shape. To ensure the compatibility of the final assembled structure, the mapping and trimming processes must satisfy certain conditions, which are detailed in the remainder of this section.

\begin{figure}[ht]
\centering
\includegraphics[width=0.5\textwidth]{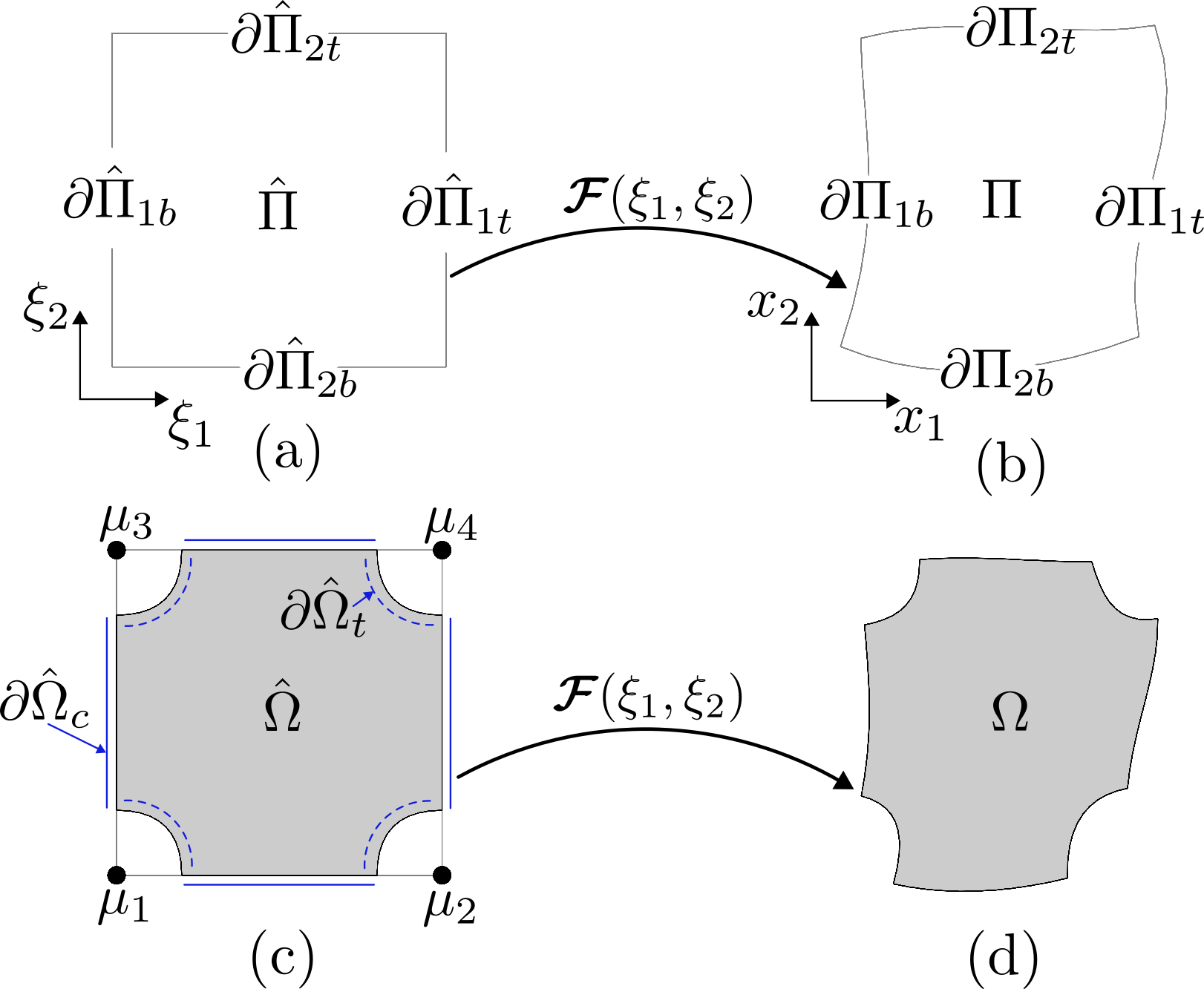}
\caption{Construction of a single cell of the lattice structure. The untrimmed parametric domain $\hat{\Pi}$ in (a) is mapped through $\bm{\mathcal{F}}$ into the physical domain $\Pi$ (b). The boundary of the untrimmed parametric domain is composed by four edges $\partial\hat{\Pi}_{1b}$, $\partial\hat{\Pi}_{1t}$, $\partial\hat{\Pi}_{2b}$, and $\partial\hat{\Pi}_{2t}$, which get mapped into the four edges of the untrimmed physical domain. The trimmed parametric domain $\hat{\Omega}$ is shown in (a). The threshold parameters $\mu_1$, $\mu_2$, $\mu_3$, and $\mu_4$ correspond to the nodes of first-order Lagrangian functions therefore the corners of $\hat{\Pi}$. The boundary of the trimmed domain is partitioned in its conformal part $\partial\hat{\Omega}_c$ and its trimmed part $\partial\hat{\Omega}_t$. Finally, in (d) it is shown the geometry of the trimmed physical domain $\Omega$.}
\label{fig: GEO single cell}
\end{figure}

\subsubsection{Mapping of a single cell} \label{sssec:single cell map}
The starting point for the geometric definition of an individual cell is the untrimmed parametric domain, $\hat{\Pi}$, which is defined as the unit square:
\begin{equation}
\hat{\Pi} = [0,1] \times [0,1].
\end{equation}
This domain represents the set where the curvilinear coordinates $(\xi_1, \xi_2)$ take values. These coordinates are used to explicitly define the map $\mathcal{F}: \hat{\Pi} \to \mathbb{R}^2$ as:
\begin{equation}\label{eq: GEO map}
\bm{\mathcal{F}} \left( \xi_1, \xi_2 \right) = \mathcal{F}_i \left( \xi_1, \xi_2 \right) \bm{e}_i \quad \text{with} \quad i = 1, 2,
\end{equation}
where $\bm{e}_i$ represents the vectors of the standard Euclidean basis.
The four edges of the parametric domain are denoted as: the left vertical edge $\partial\hat{\Pi}_{1b} = \{(0,\xi)\}$, the right vertical edge $\partial\hat{\Pi}_{1t} = \{(1,\xi)\}$, the bottom horizontal edge $\partial\hat{\Pi}_{2b} = \{(\xi,0)\}$, and the top horizontal edge $\partial\hat{\Pi}_{2t} = \{(\xi,1)\}$, where $\xi \in [0,1]$.  The untrimmed physical domain in $\mathbb{R}^2$ is the image of the unit square $\hat{\Pi}$ under this map:
\begin{equation}
{\Pi} = \bm{\mathcal{F}}(\hat{\Pi}).
\end{equation}
Figure \ref{fig: GEO single cell} shows an example of such construction. At this stage, no particular assumptions are made concerning the function $\bm{\mathcal{F}}$. However, a discussion regarding its required regularity is provided in the remainder of this section.

\subsubsection{Trimming of a single cell} \label{sssec:single cell trimming}

\begin{figure}[ht]
     \begin{subfigure}[b]{0.21\textwidth}
         \centering
         \includegraphics[width=\textwidth]{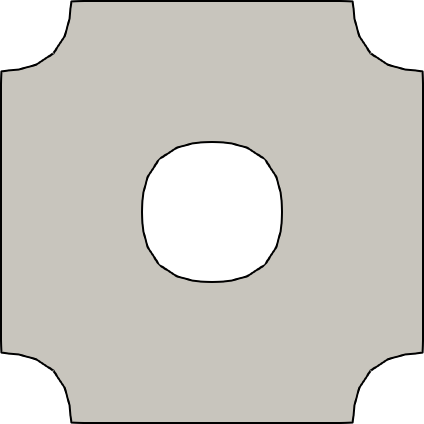}
         \caption{Schwarz Diamond}
     \end{subfigure}\hfill
     \centering
     \begin{subfigure}[b]{0.21\textwidth}
         \centering
         \includegraphics[width=\textwidth]{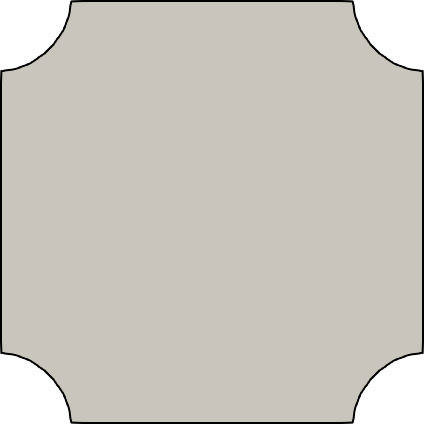}
         \caption{Schwarz Primitive}
     \end{subfigure}\hfill
     \begin{subfigure}[b]{0.21\textwidth}
         \centering
         \includegraphics[width=\textwidth]{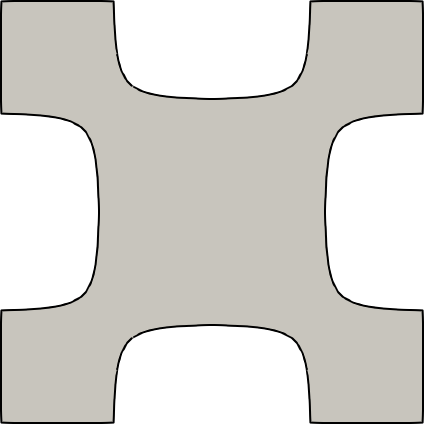}
         \caption{Schoen FRD}
     \end{subfigure}\hfill
     \begin{subfigure}[b]{0.21\textwidth}
         \centering
         \includegraphics[width=\textwidth]{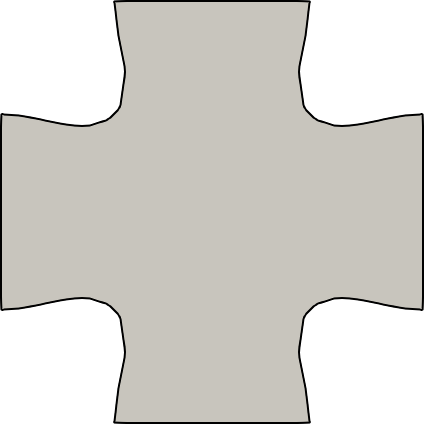}
         \caption{Schoen IWP}
     \end{subfigure}
     \caption{2D TPMS corresponding to the level-set functions in Table \ref{tab: GEO tpms}.}\label{fig: GEO TPMS geometries}
\end{figure}

\begin{table}[ht]\centering
\caption{Common TPMS in 2D. The level-sets have been adapted from their 3D counterparts by evaluation at $\xi_3=0$. The sole exception is the Schwarz Primitive that is evaluated at $\xi_3=1/2$ to ensure the geometry remains connected. For the same reason, the sign of the Schoen FRD is the opposite of its standard definition \added{and the standard threshold parameter for the Schwarz Diamond is taken as 0.1}.}

\renewcommand{\arraystretch}{1.2} 
\begin{tabular}{ll}\hline
\textbf{Name} & \textbf{level-set} $\phi_0(\xi_1,\xi_2)$ \\ \hline
Schwarz Diamond & 
$\cos(2\pi\xi_1)\cos(2\pi\xi_2)$ \\
Schwarz Primitive & 
$\cos(2\pi\xi_1) + \cos(2\pi\xi_2) - 1$  \\
Schoen FRD & 
$\cos(4\pi\xi_1)\cos(4\pi\xi_2) + \cos(4\pi\xi_2) + \cos(4\pi\xi_1) -4 \cos(2\pi\xi_1)\cos(2\pi\xi_2)$ \\
Schoen IWP &
$2 \big( \cos(2\pi\xi_1)\cos(2\pi\xi_2) + \cos(2\pi\xi_2) + \cos(2\pi\xi_1) \big) - \cos(4\pi\xi_1) - \cos(4\pi\xi_2) - 1$ \\ \hline
\end{tabular} 
\label{tab: GEO tpms} 
\renewcommand{\arraystretch}{1.0} 
\end{table}

The mapping introduced above is combined with an immersed boundary approach to define the active portion of the cell. Specifically, the trimmed parametric domain is defined as:
\begin{equation}
\hat{\Omega} = \left\{\left(\xi_1,\xi_2\right) \in \hat{\Pi} : \phi_0\left(\xi_1,\xi_2\right) < \mu \left(\xi_1,\xi_2\right)\right\} \;,
\end{equation}
where $\phi_0$ and $\mu$ are the \emph{level-set} and \emph{threshold} functions, respectively. Or equivalently, by introducing the function $\phi\left(\xi_1,\xi_2\right)=\phi_0\left(\xi_1,\xi_2\right)-\mu\left(\xi_1,\xi_2\right)$ as
\begin{equation}
\hat{\Omega} = \left\{\left(\xi_1,\xi_2\right) \in \hat{\Pi} : \phi\left(\xi_1,\xi_2\right) <0\right\} \;.
\end{equation}

While $\phi_0$ can be any arbitrary implicit function, for illustrative purposes, in this work it is selected from the class of triply periodic minimal surfaces (TPMS) adapted for 2D applications, which does not constitute a limitation for the proposed method. Common examples are defined in Table \ref{tab: GEO tpms} and depicted in Figure \figref{fig: GEO TPMS geometries}.

 Typically, the threshold is kept uniform, therefore determining which isoline of $\phi_0$ is adopted as domain boundary. However, to introduce more flexibility in the design of the cells, this work allows $\mu$ to vary linearly as follows:
\begin{equation}
\mu\left(\xi_1,\xi_2\right) = \sum_{i=1}^4 \mathscr{l}_i\left(\xi_1,\xi_2\right)\mu_i \;,
\end{equation}
where $\mathscr{l}_i$ are the classical first-degree Lagrangian shape functions, and the threshold parameters $\mu_i$ correspond to the nodal values at the bottom-left, bottom-right, top-left, and top-right corners of $\hat{\Pi}$. As for its untrimmed counterpart, the trimmed physical domain is obtained through the map:
\begin{equation}
\Omega = \bm{\mathcal{F}}(\hat{\Omega})\;.
\end{equation}
The boundary of the parametric domain is identified as two separate portions: the conformal portion lying on the untrimmed parametric domain boundary, and the trimmed portion that is internal to the domain:
\begin{subequations}
\begin{align}
    \partial\hat{\Omega}_c &= \left\{(\xi_1,\xi_2)\in\partial\hat\Pi:\phi_0\left(\xi_1,\xi_2\right) < \mu \left(\xi_1,\xi_2\right)\right\}\;, \\
    \partial\hat{\Omega}_t &= \left\{(\xi_1,\xi_2)\in\mathrm{int}(\hat{\Pi}):\phi_0\left(\xi_1,\xi_2\right) = \mu \left(\xi_1,\xi_2\right)\right\}\;,
\end{align}
\end{subequations}
where $\mathrm{int}(\bullet)$ denotes the internal part of the domain $\bullet$, such that $\partial\hat\Omega = \partial\hat\Omega_c \cup \partial\hat\Omega_t$ and $\partial\hat\Omega_c \cap \partial\hat\Omega_t = \emptyset$. In Figure \ref{fig: GEO single cell}, an example of a trimmed domain is shown in parametric and physical coordinates and its correspondent threshold parameters.

\subsubsection{Assembled lattice structure}
\begin{figure}[ht]
\centering
\includegraphics[width=0.6\textwidth]{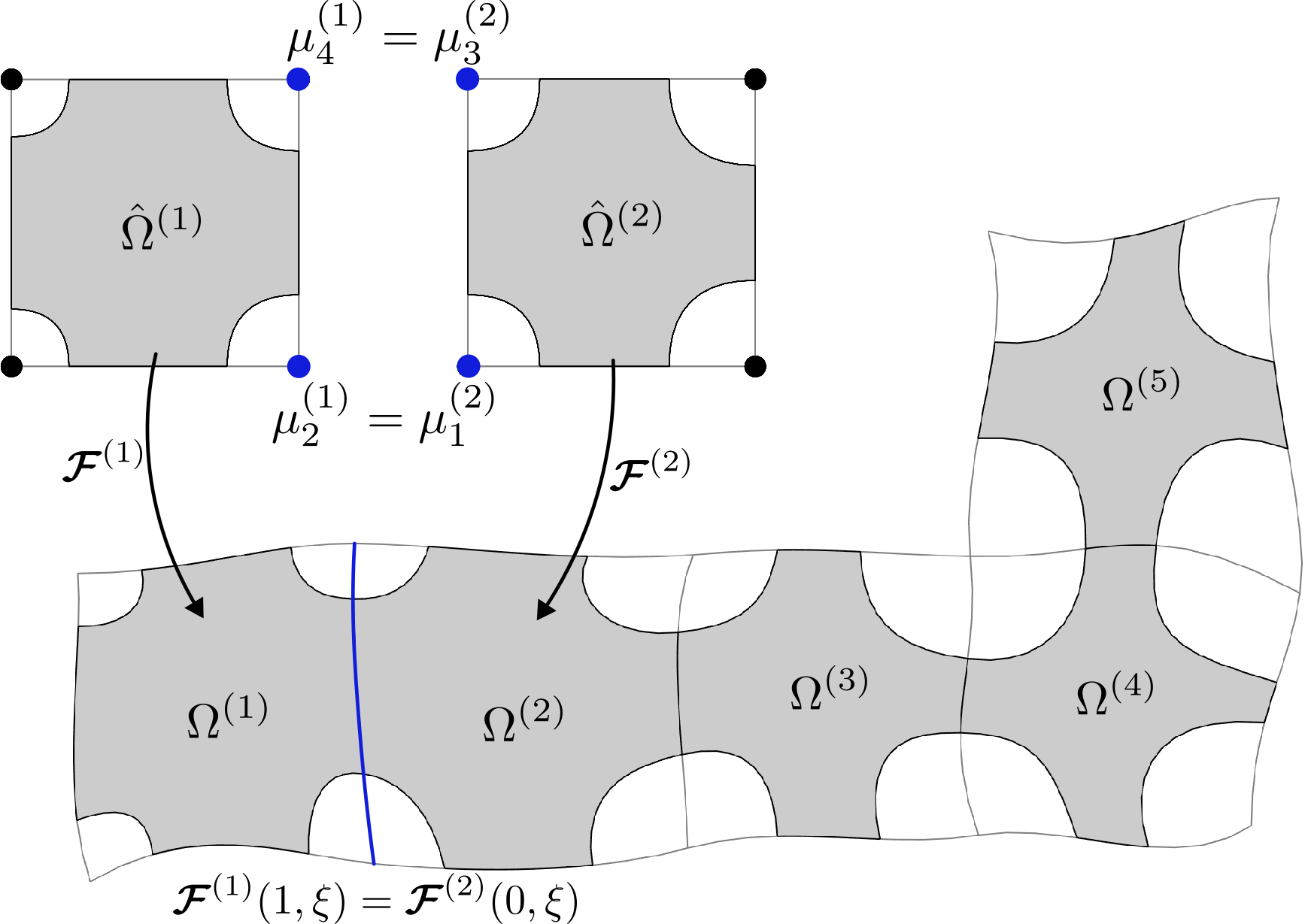}
\caption{A lattice structure is constructed by juxtaposing six cells in a compatible way. The compatibility conditions are illustrated for cells $(1)$ and $(2)$. The map of the common edge has to be the same for both cells, meaning $\bm{\mathcal{F}}^{(1)}(1,\xi)=\bm{\mathcal{F}}^{(2)}(0,\xi)$ with $\xi\in[0,1]$. To guarantee that the threshold function varies in the same way along the common edge, the following equivalence\added[id=R2]{s} are enforced between threshold parameters: $\mu_2^{(1)}=\mu_1^{(2)}$, and $\mu_4^{(1)}=\mu_3^{(2)}$.}
\label{fig: GEO construction}
\end{figure}
The discussion so far has focused on the geometric definition of a single cell. To describe a complete structure composed of multiple cells, we now introduce a superscript $(i_c)$ to distinguish all quantities associated with a specific cell, where $i_c \in \{1, \dots, n_c\}$ is the cell index and $n_c$ is the total number of cells. Therefore, the untrimmed and trimmed physical domains of the overall structure are defined, respectively, as:
\begin{equation}
{\Pi} = \mathrm{int}\left(\bigcup_{i_c=1}^{n_c} \bar{\Pi}^{(i_c)}\right), \quad\mathrm{and}\quad{\Omega} = \mathrm{int}\left(\bigcup_{i_c=1}^{n_c} \bar{\Omega}^{(i_c)}\right)\;,
\end{equation}
where $\bar{\bullet}$ denotes the closure of the domain $\bullet$. While this procedure can generate complex topologies and is not restricted to tensor-product grids, the cells must satisfy certain compatibility conditions at both the untrimmed and trimmed levels. 

First, for untrimmed compatibility, the physical mapping of connected edges must coincide. If the $i_{c_1}$-th and $i_{c_2}$-th cells are joined at their respective edges $\partial\hat{\Pi}_i^{(i_{c_1})}$ and $\partial\hat{\Pi}_j^{(i_{c_2})}$, then:
\begin{equation} \label{eq: GEO untrimmed compatibility}
\bm{\mathcal{F}}^{(i_{c_1})}\left(\partial\hat{\Pi}_i^{(i_{c_1})}\right) = \bm{\mathcal{F}}^{(i_{c_2})}\left(\partial\hat{\Pi}_j^{(i_{c_2})}\right).
\end{equation}
where $i,j\in\{1b,1t,2b,2t\}$.

Second, for trimmed compatibility, the trimming function $\phi$ must be continuous across the common interface
\begin{equation} \label{eq: GEO trimmed compatibility}
    \phi^{(i_{c_1})}(\xi)=\phi^{(i_{c_2})}(\xi),
\end{equation}
where $\xi\in[0,1]$ is an auxiliary curvilinear coordinate that parameterizes the interface. This condition, in turn, is satisfied if both the level-set and threshold functions are independently equal: $\phi_0^{(i_{c_1})}(\xi)=\phi_0^{(i_{c_2})}(\xi)$ and $\mu^{(i_{c_1})}(\xi)=\mu^{(i_{c_2})}(\xi)$.
In the particular case of TPMS geometries, the first condition ($\phi_0$) is automatically satisfied by their intrinsic periodicity,
provided that the edge pairs $(i,j)$ in Equation \eqref{eq: GEO untrimmed compatibility} are restricted to the admissible set: ${(1b,1t),(1t,1b),(2b,2t),(2t,2b)}$. The second equality ($\mu$) is satisfied by simply enforcing identical threshold parameters at the common nodes of neighboring cells.

To illustrate, in Figure \ref{fig: GEO construction}, the two adjacent cells $\Pi^{(1)}$ and $\Pi^{(2)}$ are connected at their edges $\partial{\Pi}_{1t}^{(1)}$ and $\partial{\Pi}_{1b}^{(2)}$. Compatibility of the untrimmed interface is ensured because the maps are chosen such that the following relationship holds: $\bm{\mathcal{F}}^{(1)}(1,\xi)=\bm{\mathcal{F}}^{(2)}(0,\xi)$ with $\xi\in[0,1]$. Furthermore, compatibility of the trimmed interface is achieved by imposing  $\mu^{(1)}_{2}=\mu^{(2)}_{1}$, and $\mu^{(1)}_{4}=\mu^{(2)}_{3}$.

Notably, a spline surface can be easily adopted as the untrimmed physical domain, as it naturally satisfies the compatibility condition in Equation \eqref{eq: GEO untrimmed compatibility}. The only preprocessing required is a Bézier extraction for each of its elements, which in turn defines the mapping for each cell in the lattice structure.

\subsection{Model Problem and Discretization}\label{ssec:MODEL}
With the geometric framework for the lattice structure established, we now turn to the mechanical problem defined on the trimmed physical domain $\Omega$. Specifically, we consider the linear elasticity equations, which govern the static response of the structure under prescribed loads and boundary conditions, and describe their discretization using a high-order unfitted finite element method.

\subsubsection{Linear Elasticity}
The variational statement for two-dimensional linear elasticity reads: find the displacement field $\bm{u}\in[\mathcal{H}^1]^2_{\bar{\bm{u}}}$ such that
\begin{equation} \label{eq: ELA bilinear form}
    \mathfrak{L}_{int}(\bm{v},\bm{u}) = \mathfrak{L}_{ext}(\bm{v}) \;,\quad\quad\forall\bm{v}\in[\mathcal{H}^1]^2_{\bm{0}}\;,
\end{equation}
where $[\mathcal{H}^1]^2_{\bar{\bm{u}}}$ and $[\mathcal{H}^1]^2_{\bm{0}}$ represent the subspaces of $[\mathcal{H}^1]^2$ for which Dirichlet boundary conditions are satisfied, in the first case, and where the test function $\bm{v}$ is null on the Dirichlet boundary, in the second one. It is further assumed that only pure Dirichlet or Neumann boundary conditions are considered, meaning that mixed boundary conditions, where Dirichlet is applied to one component of the displacement and Neumann to another, are excluded \added[id=R2]{for simplicity. We note that supporting such conditions is
a straightforward extension of the present framework.}

The bilinear and the linear forms in Equation \eqref{eq: ELA bilinear form} are also referred to as virtual work of the internal and external forces, respectively, and are defined as:
\begin{subequations}
    \begin{align}
         &\mathfrak{L}_{int}(\bm{v},\bm{u}) = \sum_{i_c=1}^{n_c} \mathfrak{L}_{int}^{(i_c)}(\bm{v},\bm{u}),\\
         &\mathfrak{L}_{ext}(\bm{v}) = \sum_{i_c=1}^{n_c} \mathfrak{L}_{ext}^{(i_c)}(\bm{v}) \;,
    \end{align}
\end{subequations}
where the virtual work of the internal and external forces for the $i_c$-th cell are defined as:
\begin{subequations}\label{eq: ELA virtual works}
    \begin{align}
         &\mathfrak{L}_{int}^{(i_c)}(\bm{v},\bm{u}) = \int_{\Omega^{(i_c)}} \bm{\sigma}(\bm{u}):\bm{\varepsilon}(\bm{v})\,\dd\Omega \;,\\
         &\mathfrak{L}_{ext}^{(i_c)}(\bm{v}) = \int_{\Omega^{(i_c)}} \bm{v}\cdot\bm{b}\,\dd\Omega + \int_{\partial\Omega^{(i_c)}_N}\bm{v}\cdot\bm{t}\,\dd\partial\Omega \;,
    \end{align}
\end{subequations}
where $\bm{b}$ and $\bm{t}$ are the distributed domain force and boundary traction, respectively. The components of the strain tensor $\bm{\varepsilon}$ are:
\begin{equation}
    \varepsilon_{ij}(\bm{u}) = \frac{1}{2}\left(\fpd{u_i}{x_j} + \fpd{u_j}{x_i} \right) \;,
\end{equation}
whereas the components of the stress tensor are obtained through the constitutive relationship:
\begin{equation}
    \sigma_{ij} = C_{ijkl} \, \varepsilon_{kl} \;,
\end{equation}
being $C_{ijkl}$ the elasticity tensor that for isotropic materials reads as:
\begin{equation}
    C_{ijkl} = \lambda \, \delta_{ij}\delta_{kl} + \mu \, (\delta_{ik}\delta_{jl} + \delta_{il}\delta_{jk}),
\end{equation}
where $\lambda$ and $\mu$ are the Lamé parameters of the material. The linear elasticity variational statement can be reformulated through an adequate transformation directly in the parametric domain of the cell, which will be useful later in this article for the fast assembly procedure. To do so, one shall notice that
\begin{equation}
    \fpd{u_i}{x_j} = \fpd{u_i}{\xi_k}\fpd{\xi_k}{x_j}\;,\quad \dd\Omega=j_\Omega\dd\hat{\Omega}\;,\quad \dd\partial\Omega=j_{\partial\Omega}\dd\partial\hat{\Omega}\;,\quad\mathrm{and}\quad \varepsilon_{ij}C_{ijkl} \varepsilon_{kl}=\fpd{u_i}{x_j}C_{ijkl} \fpd{u_k}{x_l}\,
\end{equation}
where the last equivalence comes from the symmetries of the elasticity tensor, $j_{\Omega}$ is the determinant of the Jacobian of the map in Equation \eqref{eq: GEO map}, and $j_{\partial\Omega}$ is the curve Jacobian associated with the boundary. Substituting into Equation \eqref{eq: ELA virtual works}, the virtual works of the internal and external forces become:
\begin{subequations}\label{eq: ELA virtual works parametric}
    \begin{align}
         &\mathfrak{L}_{int}^{(i_c)}(\bm{v},\bm{u}) = \int_{\hat{\Omega}^{(i_c)}} \fpd{u_i}{\xi_j}  \hat{C}_{ijkl} \fpd{u_k}{\xi_l}\,\dd\hat{\Omega} \;,\\
         &\mathfrak{L}_{ext}^{(i_c)}(\bm{v}) = \int_{\hat{\Omega}^{(i_c)}} \bm{v}\cdot\hat{\bm{b}}\,\dd\hat{\Omega} + \int_{\partial\hat{\Omega}^{(i_c)}_N}\bm{v}\cdot\hat{\bm{t}}\,\dd\partial\hat{\Omega} \;,
    \end{align}
\end{subequations}
where the elasticity tensor and external forces have been modified as:
\begin{subequations}\label{eq: ELA constitutive}
\begin{align}
    &\hat{C}_{ijkl} = j_{\Omega}\;\fpd{\xi_j}{x_m}C_{imkn}\fpd{\xi_l}{x_n}\;,\\
    &\hat{\bm{b}} = j_{\Omega}\;\bm{b}\;, \\
    &\hat{\bm{t}} = j_{\partial\Omega}\; \bm{t}\;.
\end{align}
\end{subequations}
It shall be mentioned that in the previous equations, the term relative to the differential geometry of the map should be enriched with a superscript $(i_c)$ since this is specific to the cell. However, to enhance readability and simplify the notation the superscript has not been used here.

\subsubsection{p-FEM finite element method}

\begin{figure}[ht]
\centering
\includegraphics[width=0.35\textwidth]{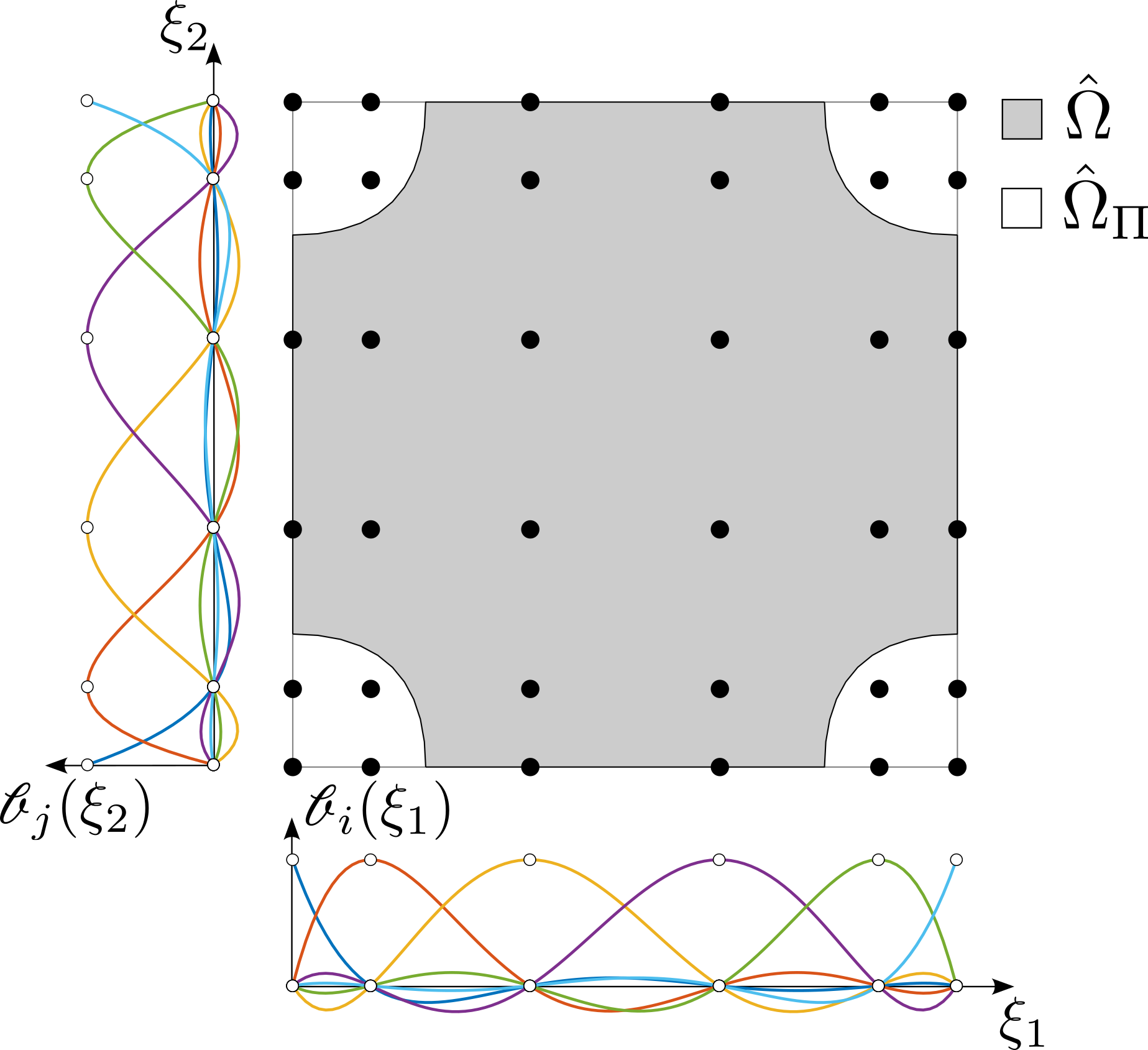}
\caption{Position of the Gauss-Lobatto-Legendre nodes of the Lagrangian basis superimposed to the trimmed parametric domain $\hat{\Omega}$ and its complement $\hat{\Omega}_\Pi$. The associated \replaced[id=R2]{univariate polynomials}{monovariate polynomial} in $\xi_1$ and $\xi_2$ are also shown.}\label{fig: MET GLL}
\end{figure}

The approximation space employed in this work is based on Lagrange polynomials defined on Gauss-Lobatto-Legendre (GLL) nodes, constructed over the untrimmed parametric domain $\hat{\Pi}^{(i_c)}$. Accordingly, the polynomials are expressed as functions of the curvilinear coordinates $\xi_1$ and $\xi_2$. Figure \ref{fig: MET GLL} illustrates the univariate polynomials of degree 6 in both directions, along with the associated tensor-product GLL nodes. The corresponding bivariate polynomials are defined as:
\begin{equation}
\mathscr{B}_{ij}(\xi_1,\xi_2) = \mathscr{b}_i(\xi_1) \cdot \mathscr{b}_j(\xi_2), \quad \mathrm{with} \quad i,j \in \{0,\dots,p\},
\end{equation}
where $\mathscr{b}_i(\xi_1)$ and $\mathscr{b}_j(\xi_2)$ denote the $i$-th and $j$-th univariate GLL polynomials in the $\xi_1$ and $\xi_2$ directions, respectively, and $p$ represents the degree of the approximation space. Formally, the approximation space for the $i_c$-th cell is given by
\begin{equation}
\mathcal{S}_h=\mathrm{span}\{\mathscr{B}_{ij}\circ\bm{\mathcal{F}}^{-1}: i\in\{0,\dots,p\} ,\; j\in\{0,\dots,p\}\},    
\end{equation}
It is worth mentioning that in the p-FEM context, the mesh is not refined beyond the level of the individual cell, thereby making the approximation element and the geometric cell equivalent. Despite this seemingly coarse discretization of the individual cell, the high polynomial order ensures substantial approximation accuracy as it will be shown in Section \ref{sec: RES}. 

Moreover, unlike other high-order bases, the GLL basis features nodes located at the domain boundaries, as illustrated in Figure \ref{fig: MET GLL} for $p=6$. Consequently, degrees of freedom at cell interfaces could, in principle, be strongly enforced. However, as will be detailed in Section \ref{sec: DD}, the domain decomposition strategy adopted in this work does not impose a strong coupling of all interface degrees of freedom, but still exploits the coincidence of their position.

\subsubsection{Unfitted FEM} \label{sec: MOD unfitted}
The polynomial basis and the location of the nodes are initially defined on the untrimmed parametric domain $\hat{\Pi}^{(i_c)}$. After the trimming operation, however, the support of each basis function is restricted to its active portion corresponding to the trimmed parametric domain of the cell, $\hat{\Omega}^{(i_c)}$. Importantly, regardless of the shape of the trimming, each basis function remains active, even if the associated node lies outside the active domain, leading to fully dense stiffness matrices. 

One might argue that, for a comparable approximation accuracy, using lower-degree polynomials with mesh refinement could provide the advantage of increased sparsity. Nevertheless, the property that the basis functions remain active irrespective of the trimming geometry is particularly advantageous for building a surrogate model of the stiffness, as discussed in Section \ref{sec: SPEED}, and is a primary motivation for the adoption of p-FEM elements.

Unfitted discretizations naturally introduce three main challenges: quadrature, application of Dirichlet boundary conditions, and conditioning.

\begin{figure}[ht]\centering
    \begin{subfigure}[b]{0.42\textwidth}
        \centering
        \begin{tikzpicture}
            \node[anchor=south west,inner sep=0] (img)
                {\includegraphics[width=\textwidth]{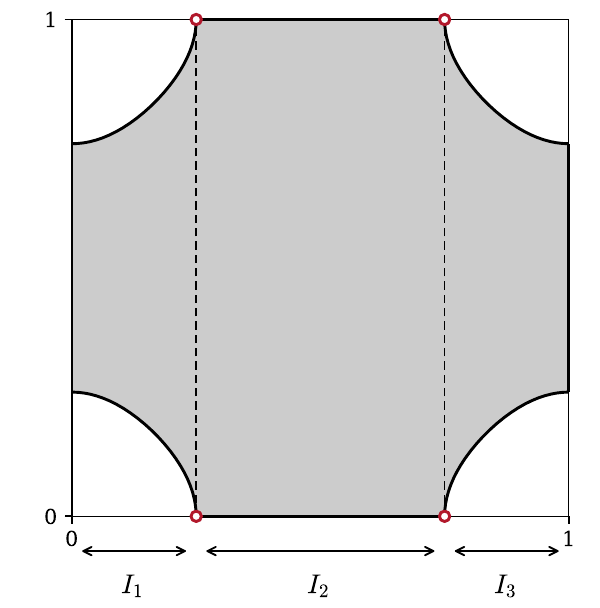}};
            \begin{scope}[x={(img.south east)},y={(img.north west)}]
                \node at (0.534,0.559) {$\hat{\Omega}$};
            \end{scope}
        \end{tikzpicture}
        \caption{}
        \label{fig: MOD quadrature a}
    \end{subfigure}\hfill
    \begin{subfigure}[b]{0.42\textwidth}
        \centering
        \includegraphics[width=\textwidth]{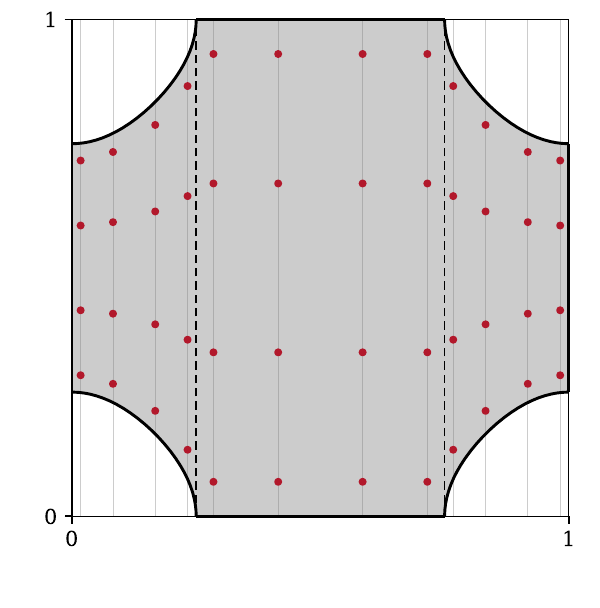}
        \caption{}
        \label{fig: MOD quadrature b}
    \end{subfigure}
\caption{\added[id=R1]{Height-function quadrature of QUGaR/algoim over the active domain $\hat{\Omega}$ of a Schwarz Primitive cell, with $\xi_2$ selected as the integration direction. (a) The base direction $\xi_1$ is split at the critical points of the level-set (red circles) into subintervals $I_1,I_2,I_3$ over which the height function is smooth. (b) A tensor-product Gauss rule then places the quadrature points within each subinterval, clustering them in the active domain.}}\label{fig: MOD quadrature}
\end{figure}

\paragraph{Quadrature algorithm}
\added[id=R1]{Quadrature over the complex trimmed domains is handled by the recently released QUGaR library \cite{qugar}, which provides a high-order quadrature algorithm for level-set-defined geometries within the FEniCSx framework \cite{dolfinx2023preprint}, building on the algoim library \cite{algoim,saye2015high,saye2022high}. The algorithm recasts the integral over the implicitly defined domain as a recursive sequence of one-dimensional integrations, thereby avoiding any explicit reconstruction of the trimmed boundary. Once an integration direction has been selected, the height function induced by the level-set, that is, the location of its roots along that direction, provides the integration bounds in that coordinate. The remaining base direction is then partitioned into subintervals delimited by the critical points of the level-set, where the height function loses smoothness or the number of roots changes. In the interior of each subinterval the integrand is infinitely smooth, so that a standard quadrature rule (e.g. Gauss-Legendre or tanh-sinh), applied in a tensor-product fashion, evaluates the integral to high accuracy. The procedure is illustrated in Figure \figref{fig: MOD quadrature} for one of the level-set functions listed in Table \ref{tab: GEO tpms}. This methodology has proven to yield high-order accuracy and spectral convergence of the integration error. Consequently, the consistency error associated with the numerical integration can be safely assumed to be negligible. Furthermore, since changes in the threshold functions alter the underlying geometry, a distinct set of quadrature nodes and weights must be generated for each set of threshold parameters.}

\paragraph{Application of boundary conditions}
\deleted[id=R2]{The application of Dirichlet boundary conditions on trimmed portions of the boundary $\partial\hat{\Omega}_t$ is not considered here, as such conditions are uncommon in lattice structures. Dirichlet conditions on conformal portions of the boundary $\partial\hat{\Omega}_c$, however, can be applied in a strong sense by acting directly on the degrees of freedom that lie on $\partial\hat{\Pi}$.}
\added[id=R2]{We restrict our formulation to boundary conditions applied on the conformal portions of the boundary, $\partial\hat{\Omega}_c$. Dirichlet conditions on $\partial\hat{\Omega}_c$ can be applied directly in a strong sense by acting on the degrees of freedom lying on $\partial\hat{\Pi}$. Conversely, we do not consider Dirichlet boundary conditions on trimmed boundary portions, $\partial\hat{\Omega}_t$. This is a deliberate choice since, in an unfitted setting, the degrees of freedom do not generally align with the cut boundary. As such, enforcing Dirichlet conditions there requires weak imposition techniques. Integrating such methods into the BDDC preconditioner is mathematically non-trivial and is deferred to future work. Non-homogeneous Neumann conditions are supported on the conformal portions of the boundary, including their active parts. On trimmed boundary portions, such conditions pose no mathematical difficulty, as they only involve an integral over the cut boundary, and they are indeed employed in the single-cell study of Section~\ref{sssec: RES single cell}; however, they are not covered by the ROM acceleration of Section~\ref{sec: SPEED} and are therefore not considered in the remainder of this work. While boundary conditions on trimmed boundaries might be relevant for local fixtures, local loads, or contact zones, they are uncommon for the interior surfaces of lattice structures, which are typically left unloaded.}

\paragraph{Stabilization}
Depending on the configuration of the trimming, poorly conditioned stiffness matrices may arise even with a spectral basis. To address this issue, we employ the $\alpha$-stabilization proposed in \cite{dauge2015theoretical}. This non-consistent term is defined as:
\begin{equation} \label{eq: MOD stabilization}
    \mathfrak{L}_{sta}(\bm{v},\bm{u}) = \rho\sum_{i_c=1}^{n_c}\int_{\hat{\Omega}^{(i_c)}_\Pi} \fpd{u_i}{\xi_j}  \hat{C}_{ijkl} \fpd{u_k}{\xi_l}\,\dd\hat{\Omega},
\end{equation}
where $\rho$ is a constant that balances the condition number with the consistency error. The stabilization is applied over the complement of the active domain, $\hat{\Omega}^{(i_c)}_\Pi = \hat{\Pi}^{(i_c)} \setminus \hat{\Omega}^{(i_c)}$. Physically, this stabilization can be interpreted as immersing the domain in a softer material, which, when $\rho$ is small, does not significantly contribute to the overall stiffness of the structure.

 \section{Domain Decomposition Method} \label{sec: DD}
This section presents the Balancing Domain Decomposition by Constraints (BDDC) method \added[id=R2]{\cite{mathew2008domain,dohrmann2003preconditioner}}, a non-overlapping domain decomposition technique used to solve the linear system arising from the assembly of the lattice structures. In the present framework, each cell is associated with a corresponding subdomain in the domain decomposition setting. The global Schur complement system, obtained after condensing the internal degrees of freedom (DoFs) of each subdomain, is solved iteratively using the preconditioned conjugate gradient (PCG) method. The BDDC preconditioner decomposes the correction of the residual into two complementary parts: (i) independent local Neumann-type solves in each subdomain, made invertible by constraining the correction to vanish at a selected set of coarse DoFs, and (ii) a global coarse solve that balances the local corrections through a basis spanning the coarse DoF space. Since the resulting correction is only guaranteed to be continuous at the coarse DoFs, a weighted averaging step using a discrete partition of unity restores full continuity across the entire interface skeleton.

Two important aspects of the problem formulation must be noted already at this stage. First, for the 2D elasticity problems considered, there is no one-to-one correspondence between degrees of freedom\deleted[id=R1]{ (DoFs)} and nodes; each node is associated with two DoFs, corresponding to the components of the displacement vector. Second, eventual Dirichlet boundary conditions are enforced strongly along entire edges even though only part of it is effectively active. Furthermore, only pure Dirichlet and Neumann boundary conditions are considered, while mixed boundary conditions are neglected.

To improve readability, quantities associated with an individual cell will be denoted by a superscript $(i)$, which differs from the notation used in the previous section. Throughout this section, the following notational convention is adopted: a tilde denotes the redundant (vertically stacked) version of a quantity, obtained by collecting contributions from all cells, e.g., $\tilde{u}$; the absence of a tilde denotes the assembled (non-redundant) counterpart, e.g., $u$; and a superscript $(i)$ denotes the local version associated with $i$-th cell, e.g., $u^{(i)}$. The same convention applies to operators such as $\tilde{\mathcal{R}}_U$ and $\mathcal{R}_U^{(i)}$, and to DoF counts such as $\tilde{n}_U$, $n_U$, and $n_U^{(i)}$. Furthermore, the following \emph{vertical stacking} and \emph{diagonal stacking} operators are introduced here and adopted throughout this section:
\begin{equation*}
    \mathrm{vstack}(v_1\dots,v_n)=\begin{bmatrix} v_1 \\v_2\\ \\\vdots\\v_n \end{bmatrix},\quad\mathrm{and}\quad
    \mathrm{blockdiag}(M_1\dots,M_n)=\begin{bmatrix}M_1&0&\cdots&0\\
    \vdots&\vdots&\ddots&\vdots\\0&0&\cdots& M_n\end{bmatrix}.
\end{equation*}

\subsection{Preliminary decompositions}
\begin{figure}[ht]
\centering
\includegraphics[width=0.7\textwidth]{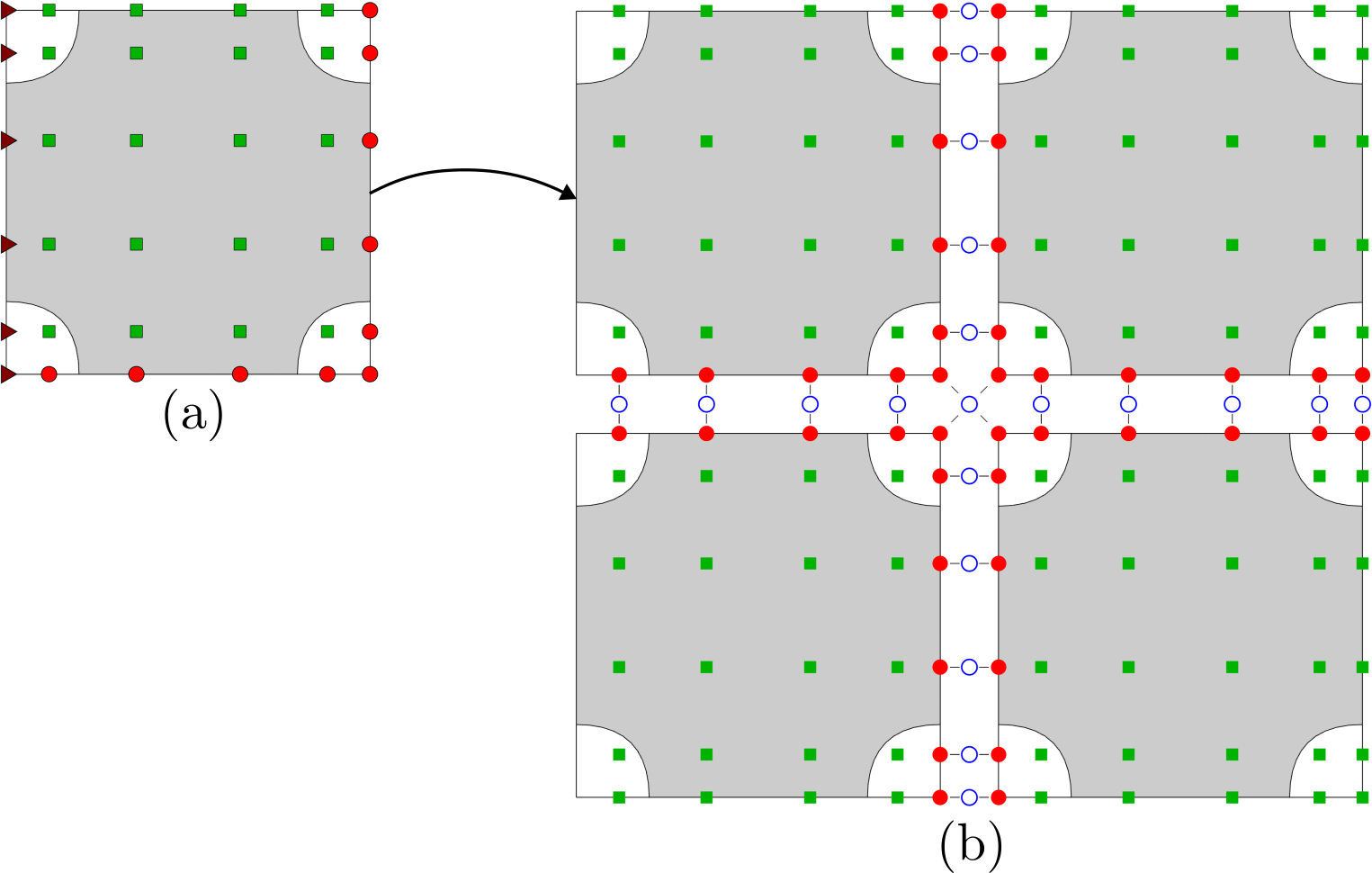}
\caption{
    Classification of the DoFs associated with a single cell (a) and with a two-by-two cells structure (b). For the single cell in subfigure (a): 
    \includegraphics[height=.7em]{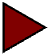} denotes the local DoFs associated with the Dirichlet boundary conditions collected in $u_D^{(i)}$; 
    \includegraphics[height=.7em]{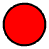} denotes the local DoFs whose associated nodes lie on the internal skeleton, collected in $u^{(i)}$; 
    \includegraphics[height=.7em]{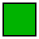} denotes the local internal DoFs collected in $u_I^{(i)}$. 
    Regarding the entire structure in subfigure (b): 
    \includegraphics[height=.7em]{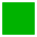} denotes the global internal DoFs collected in $\tilde{u}_I$; 
    \includegraphics[height=.7em]{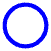} denotes the non-redundant skeleton DoFs collected in $u$; 
    \includegraphics[height=.7em]{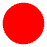} denotes the redundant global skeleton DoFs collected in $\tilde{u}$. 
    The relationship between $\tilde{u}$ and $u$, visualized by dashed lines in (b), is enforced by the equation $\tilde{u}=\tilde{\mathcal{R}}_U u$.
    }
    \label{fig: BDDC DoFs classification}
\end{figure}

Figure \figref{fig: BDDC DoFs classification}{a} illustrates a single cell together with the nodes associated with the\deleted[id=R1]{degrees of freedom} DoFs of its spectral approximation space. These DoFs are grouped into three distinct vectors: $u_D^{(i)}$, $u^{(i)}$, and $u_I^{(i)}$. The first vector, $u_D^{(i)}$ (triangles), contains the local Dirichlet DoFs, i.e., those at which Dirichlet boundary conditions are imposed in strong form. The second vector, $u^{(i)}$ (squares), contains the local skeleton DoFs, defined as those that do not belong to $u_D^{(i)}$ and lie on an interface with another cell. Finally, the third vector, $u_I^{(i)}$ (solid circles), contains the remaining internal DoFs. The sizes of such vectors are denoted as $n_D^{(i)}$, $n_U^{(i)}$, and $n_I^{(i)}$, respectively. For the assembled structure, the internal and skeleton DoFs are stacked into their global counterparts, denoted by $\tilde{u}_I$ and $\tilde{u}$, respectively, as follows:
\begin{equation}
\tilde{u}_I=\mathrm{vstack} \left(u_I^{(1)}, \dots, u_I^{(n_c)}\right),\quad \mathrm{and}\quad
\tilde{u} =\mathrm{vstack} \left(u^{(1)}, \dots, u^{(n_c)}\right)\;.
\end{equation}
It is worth recalling that $n_c$ denotes the number of cells and coincides with the number of subdomains in the decomposition. The dimensions of these vectors are given by $\tilde{n}_I=\sum_{i=1}^{n_c} n_I^{(i)}$ and $\tilde{n}_U=\sum_{i=1}^{n_c} n_U^{(i)}$. While the local internal DoFs $u_I^{(i)}$ vectors are mutually independent, the skeleton DoFs are not. Consequently, the vector $\tilde{u}$ contains redundant DoFs. In contrast, we denote by $u$ the vector of independent non-redundant skeleton DoFs, and by $n_U$ its dimension. The relationship between $\tilde{n}_U$ and $n_U$ depends on the specific geometry of the mesh. After introducing a global numbering of the non-redundant skeleton DoFs, we introduce the extension/restriction operators $\mathcal{R}_U^{(i)}\in\mathbb{R}^{n_U^{(i)}\times n_U}$ associated to the relation
\begin{equation}
u^{(i)}=\mathcal{R}_U^{(i)} u,
\end{equation}
that extract the vector of local skeleton DoFs from its global non-redundant version. $\mathcal{R}_U^{(i)}$ is therefore composed of zeros and ones accordingly to the global and local DoFs orderings. The corresponding global expression reads
\begin{equation}
\tilde{u}=\tilde{\mathcal{R}}_U u,
\end{equation}
with $\tilde{\mathcal{R}}_U\in\mathbb{R}^{\tilde{n}_U\times n_U}$
\begin{equation}
\tilde{\mathcal{R}}_U =\mathrm{vstack}\left(\mathcal{R}_U^{(1)},\dots,\mathcal{R}_U^{(n_c)}\right).
\end{equation}

\subsection{\added[id=R2]{Assembly of the stiffness matrix}}
The standard assembly of the bilinear forms introduced in Section \ref{ssec:MODEL}, for the $i$-th cell, yields a local stiffness matrix that can be partitioned into submatrices denoted by $K_{AB}^{(i)}$. Here, the subscripts $A$ and $B$ take values in $\{I,U,D\}$, depending on whether the corresponding rows and columns are associated with $u_I^{(i)}$, $u^{(i)}$, and $u_D^{(i)}$, respectively. Similarly, the assembly of the linear forms in Section \ref{ssec:MODEL} yields the local external force vector, which is partitioned into $f_I^{*(i)}$ and $f_U^{*(i)}$. These local force vectors are further adjusted by condensing the Dirichlet DoFs as
\begin{align}
    &f_I^{(i)}=f_I^{*(i)}-K_{ID}^{(i)} u_D^{(i)} \;,\\
    &f_U^{(i)}=f_U^{*(i)}-K_{UD}^{(i)} u_D^{(i)} \;.
\end{align}
Assembling the local stiffness matrices and force vectors over the entire structure, while preserving the separation between internal and skeleton DoFs in the global ordering, yields the following linear system:
\begin{equation} \label{eq: BDDC full system}
    \begin{bmatrix} \tilde{K}_{II} & \tilde{K}_{IU}\tilde{\mathcal{R}}_U \\
    \tilde{\mathcal{R}}_U^\top \tilde{K}_{IU}^\top & \tilde{\mathcal{R}}_U^\top \tilde{K}_{UU}\tilde{\mathcal{R}}_U \end{bmatrix}
    \begin{bmatrix} \tilde{u}_I\\u \end{bmatrix} = \begin{bmatrix}\tilde{f}_I\\ \tilde{\mathcal{R}}_U^\top \tilde{f}_{U}\end{bmatrix}
\end{equation}
Here, 
\begin{equation}
    \tilde{K}_{AB}=\mathrm{blockdiag}\left(K_{AB}^{(1)},\dots,K_{AB}^{(n_c)}\right)\;,\quad \mathrm{and}\quad \tilde{f}_{A}=\mathrm{vstack}\left(f_{A}^{(1)},\dots,f_{A}^{(n_c)}\right)\;,
\end{equation}
where, the subscripts $A$ and $B$ take values in $\{I,U\}$. Then, condensing the internal DoFs within each cell leads to the local Schur complement $S^{(i)}$ and the corresponding condensed force vector $f^{(i)}$, defined by
\begin{align*}
    &S^{(i)} = K_{UU}^{(i)} - K_{UI}^{(i)}K_{II}^{(i) -1}K_{IU}^{(i)}\;,\\
    &f^{(i)} = f_U^{(i)} - K_{UI}^{(i)}K_{II}^{(i) -1}f_{I}^{(i)}\;.
\end{align*}
Note that $K_{II}^{(i)}$ is invertible, since neither skeleton nor external boundary DoFs are included in its definition. In this work, the local Schur complements are assembled explicitly, although the inverse $K_{II}^{(i)-1}$ is never formed. Instead, the action of $K_{II}^{(i)-1}$ on $K_{IU}^{(i)}$ is computed by solving linear systems of the form $K_{II}^{(i)} X = K_{IU}^{(i)}$, i.e., by applying a local solve to each column of $K_{IU}^{(i)}$. The resulting global system reads
\begin{equation*}
    \tilde{\mathcal{R}}_U^\top \tilde{S} \tilde{\mathcal{R}}_U \,u = \tilde{\mathcal{R}}_U^\top \tilde{f}
\end{equation*}
where 
\begin{equation*}
    \tilde{S} = \mathrm{blockdiag}\left(S^{(1)},\dots,S^{(n_c)}\right)\;,\quad\mathrm{and}\quad\tilde{f} = \mathrm{vstack}\left(f^{(1)},\dots,f^{(n_c)}\right)\;.
\end{equation*}
Equivalently,
\begin{equation}\label{eq: BDDC system}
    S u = f
\end{equation}
with $S = \tilde{\mathcal{R}}_U^\top \tilde{S} \tilde{\mathcal{R}}_U$, and $f = \tilde{\mathcal{R}}_U^\top \tilde{f}$, where $S$ is the global assembled Schur complement.

\subsection{The BDDC method}
In the BDDC method \cite{badia2018robust,mathew2008domain}, the system in Equation \eqref{eq: BDDC system} is solved using an iterative method, here the preconditioned conjugate gradient (PCG) method, where the construction of the preconditioner $S_{BDDC}^{-1}$ is described in the following subsection. Although the details are omitted here for brevity, the PCG algorithm requires the multiplication of search directions by $S$ and of residuals by $S_{BDDC}^{-1}$. 

In this contribution, neither of these quantities are explicitly assembled. Only local Schur complements $S^{(i)}$ are assembled and stored. Therefore, the action of $S$ on a vector is performed through sequences of matrix-vector multiplications that can be efficiently distributed in parallel. In contrast, the application of $S_{BDDC}^{-1}$ also involves local solves, which are likewise readily parallelizable, with only one coarse global solve performed using a parallel Cholesky factorization (although alternative solvers for the coarse correction could also be employed).

The BDDC preconditioner consists of two components: a coarse correction, which addresses errors associated with global low-frequency modes affecting the entire structure, and a fine correction, which targets errors in local high-frequency behaviors. To properly separate these two contributions, a set of coarse DoFs is selected.

\subsubsection{Coarse degrees of freedom}
\begin{figure}[ht]
\centering
\includegraphics[width=0.7\textwidth]{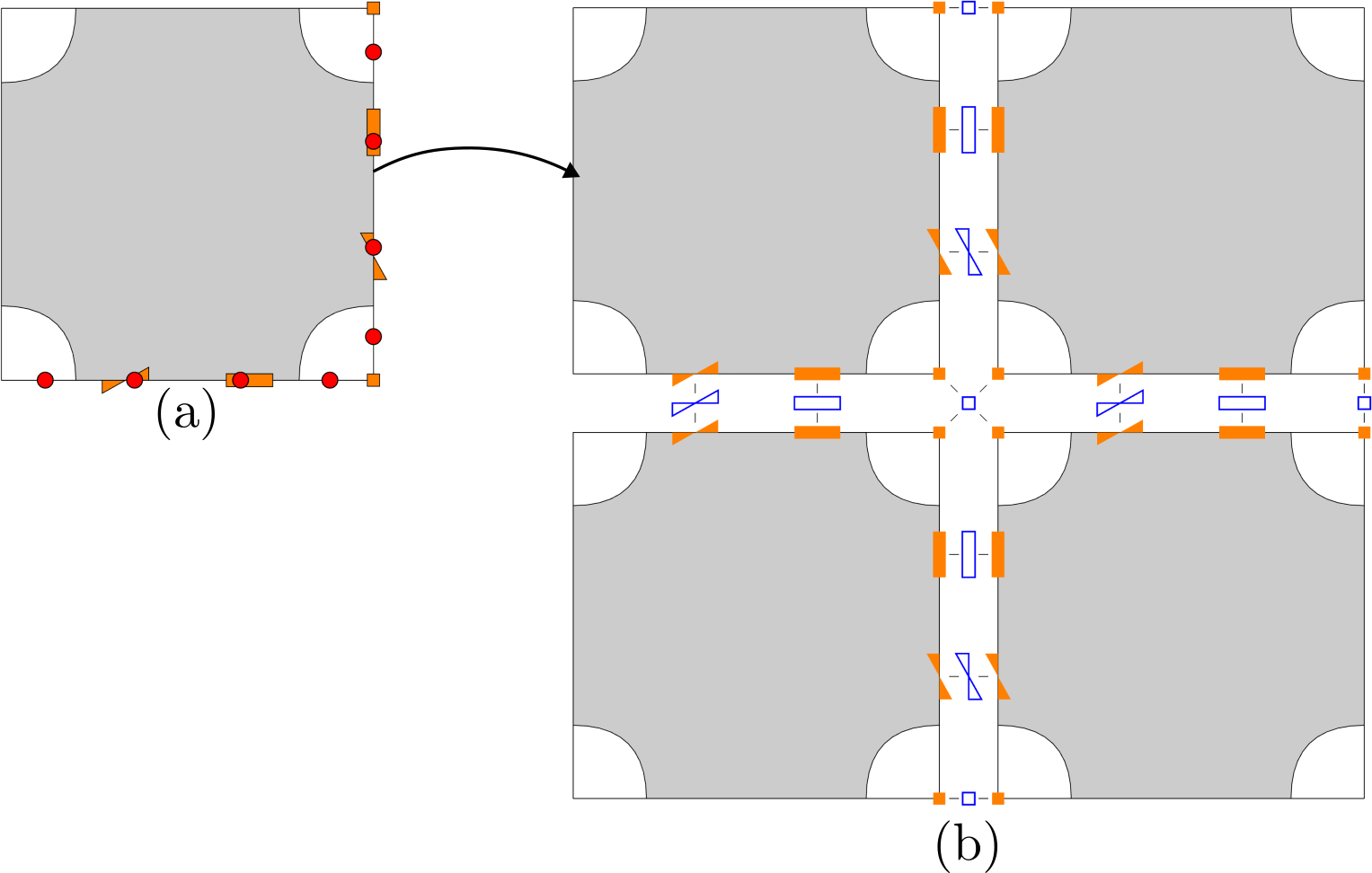}
\caption{
    Classification of the coarse DoFs associated with a single cell (a) and with a two-by-two cells structure (b). 
    For the single cell in subfigure (a): 
    \includegraphics[height=.7em]{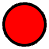} local DoFs of the internal skeleton used to compute edge averages and moments; 
    \includegraphics[height=.7em]{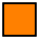} local DoFs at crosspoints that are also coarse DoFs; 
    \includegraphics[height=.7em]{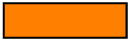} local averages of DoFs on edges lying on the internal skeleton; 
    \includegraphics[height=1em]{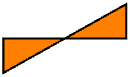} local moments of DoFs on edges lying on the internal skeleton. 
    In particular, these last three categories are collected in the vector of local coarse DoFs $c^{(i)}$.
    Regarding the entire structure in subfigure (b): 
    \includegraphics[height=.7em]{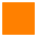} global redundant DoFs at crosspoints that are also coarse DoFs; 
    \includegraphics[height=.7em]{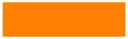} global redundant averages of DoFs on internal skeleton edges; 
    \includegraphics[height=1em]{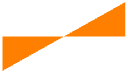} global redundant moments of DoFs on internal skeleton edges. 
    These last three categories are collected in the global vector of redundant coarse DoFs $\tilde{c}$;  
    \includegraphics[height=.7em]{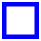} global non-redundant DoFs at crosspoints that are also coarse DoFs; 
    \includegraphics[height=.7em]{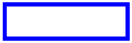} global non-redundant averages of DoFs on internal skeleton edges; 
    \includegraphics[height=1em]{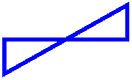} global non-redundant moments of DoFs on internal skeleton edges. 
    These last three categories are collected in the global vector of non-redundant coarse DoFs $c$. The relationship between $\tilde{c}$ and $c$, visualized by dashed lines in (b), is enforced by the equation $\tilde{c}=\tilde{\mathcal{R}}_C c$.}
    \label{fig: Coarse DoFs classification}
\end{figure}
Coarse\deleted[id=R1]{degrees of freedom} DoFs are selected skeleton DoFs, or combinations thereof, used to capture the global behavior of the structure. In this work, three types of coarse DoFs are considered: cross-point values, edge averages, and edge moments (see Figure \figref{fig: Coarse DoFs classification}). This choice is inspired by \added[id=R2]{\cite{heinlein2019frugal, mandel2005algebraic}}, although not identical. Cross points are nodes located at the intersection of mesh lines. The coarse DoFs associated with the $i_x$-th cross point are defined as

\begin{subequations}
\begin{align}
    c_{i_x}^{u_1}&=u_{1}|_{\boldsymbol{x}_{i_x}} \,,\\
    c_{i_x}^{u_2}&=u_{2}|_{\boldsymbol{x}_{i_x}} \,,
\end{align}
\end{subequations}
where $u_{1}$ and $u_{2}$ denote the components of the displacement vector, and $\boldsymbol{x}_{i_x}$ are the coordinates of the $i_x$-th cross point. The remaining coarse DoFs are associated to the $i_e$-th edge. Two average coarse DoFs are introduced, one for each displacement component:
\begin{subequations}
\begin{align}
    c_{i_e}^{a_1}&=\int_{\partial\Omega_{i_e}} u_{1} \, \text{d}\partial\Omega \,,\\
    c_{i_e}^{a_2}&=\int_{\partial\Omega_{i_e}} u_{2} \, \text{d}\partial\Omega \,,
\end{align}
\end{subequations}
where $\partial\Omega_{i_e}$ denotes the domain of the $i_e$-th edge. In addition, one moment-based coarse DoF is defined as
\begin{equation}
c_{i_e}^{m} = \int_{\partial\Omega_{i_e}} \boldsymbol{f} \cdot \bm{u} \, \text{d}\partial\Omega \,,
\end{equation}
where $\boldsymbol{f} = \begin{pmatrix} x_2 - \bar{x}_2 \\ \bar{x}_1 - x_1 \end{pmatrix}$ and $(\bar{x}_1, \bar{x}_2)$ are the coordinates of the center of the edge. Collecting all coarse DoFs into the vector $c$ of size $n_C$, the following relationship holds:
\begin{equation}
    c = Q u \;,
\end{equation}
where $Q \in \mathbb{R}^{n_C \times n_U}$ is the matrix that enforces the above definitions through adequate integration and evaluation of basis functions. Figure \figref{fig: Coarse DoFs classification} shows the different types of coarse DoFs for a single cell and for a two-by-two cells structure. As for the skeleton DoFs, the coarse DoFs vector on the $i$-th cell, denoted by $c^{(i)}$ and of size $n_C^{(i)}$, is obtained through an extension/restriction operator $\mathcal{R}_C^{(i)}\in\mathbb{R}^{n_C^{(i)} \times n_C}$,
\begin{equation}
    c^{(i)} = \mathcal{R}_C^{(i)} c.
\end{equation}
Similarly to $\mathcal{R}_U^{(i)}$, the matrix $\mathcal{R}_C^{(i)}$ consists of properly placed zeros and ones according to the global and local orderings of the coarse DoFs. By stacking the local coarse DoFs, the vector of redundant coarse DoFs $\tilde{c}$ of size $\tilde{n}_C=\sum_{i=1}^{n_c} n_C^{(i)}$ is defined as
\begin{equation}
    \tilde{c} = \mathrm{vstack}\left(c^{(1)}, \dots, c^{(n_c)}\right)\;,
\end{equation}
and the associated matrix $\tilde{\mathcal{R}}_C\in\mathbb{R}^{\tilde{n}_C\times n_C}$
\begin{equation*}
    \tilde{\mathcal{R}}_C = \mathrm{vstack}\left(\mathcal{R}_C^{(1)}, \dots, \mathcal{R}_C^{(n_c)}\right)
\end{equation*}
performs the operation
\begin{equation*}
    \tilde{c} = \tilde{\mathcal{R}}_C \, c \;.
\end{equation*}
Additionally, the constraint matrices $C^{(i)}\in\mathbb{R}^{n_C^{(i)}\times n_U^{(i)}}$ are introduced to extract the coarse DoFs of cell $i$ from its skeleton DoFs:
\begin{equation*}
    C^{(i)} = \mathcal{R}_C^{(i)} Q \mathcal{R}_U^{(i)\top}\;,
\end{equation*}
associated to the relation
\begin{equation*}
    c^{(i)}=C^{(i)}u^{(i)} \;.
\end{equation*}
The global version of the precedent relation is expressed as:
\begin{equation*}
    \tilde{c}=\tilde{C}\tilde{u} \;,
\end{equation*}
where matrix $\tilde{C}\in\mathbb{R}^{\tilde{n}_C\times \tilde{n}_U}$ is defined as
\begin{equation*}
    \tilde{C} = \mathrm{blockdiag}\left(C^{(1)}, \dots, C^{(n_c)}\right) \;.
\end{equation*}

\subsubsection{Coarse basis functions}
Associated with the coarse DoFs $c$, the matrix $\tilde{\Psi}$ of size $\tilde{n}_U\times n_C$ is introduced, whose columns represent a coarse basis for the global correction. The matrix $\tilde{\Psi}$ is defined through
\begin{equation*}
    \begin{bmatrix}\tilde{S}&\tilde{C}^\top\\\tilde{C}&0\end{bmatrix}
    \begin{bmatrix}\tilde{\Psi}\\\tilde{\Lambda}\end{bmatrix} = 
    \begin{bmatrix}0\\\mathcal{\tilde{R}}_C\end{bmatrix}\;,
\end{equation*}
which implies that the $i$-th column of the matrix $\tilde{\Psi}$ is obtained as the solution associated with the $i$-th column of the matrix $\tilde{\mathcal{R}}_C$. The matrix $\tilde{\Lambda}$ contains the Lagrange multipliers associated with the coarse DoF constraints. Additionally, due to the block structures of $\tilde{S}$, $\tilde{C}$, and $\tilde{\mathcal{R}}_C$, the previous problem can be decomposed, and therefore parallelized, into local problems as
\begin{equation}\label{eq: BDDC local basis}
    \begin{bmatrix}S^{(i)}&C^{(i)\top}\\C^{(i)}&0\end{bmatrix}
    \begin{bmatrix}\Psi^{(i)}\\\Lambda^{(i)}\end{bmatrix} = 
    \begin{bmatrix}0\\\mathcal{R}_C^{(i)}\end{bmatrix}\;,
\end{equation}
where $\Lambda^{(i)}$ contains the Lagrange multipliers associated to coarse DoFs at the cell level. Furthermore, it should be noted that while $\tilde{\mathcal{R}}_C$ has at least two nonzero elements in every column, this is not necessarily the case for the local matrix $\mathcal{R}_C^{(i)}$ in cell $i$. Therefore, only a subset of the problems \eqref{eq: BDDC local basis} needs to be solved for each cell. The coarse basis functions $\tilde{\Psi}$ are then assembled as
\begin{equation*}
    \tilde{\Psi} =\mathrm{vstack}\left(\Psi^{(1)},\dots,\Psi^{(n_c)}\right)\;.
\end{equation*}
From a physical point of view, the $i$-th basis function (i.e., the $i$-th column of $\tilde{\Psi}$) represents the displacement field obtained by solving in each cell a local Dirichlet problem where, in addition to homogeneous Dirichlet conditions on the external boundary DoFs $u_D^{(i)}$, one coarse DoF is set equal to one while all the others are set to zero.

\subsubsection{Discrete partition of unity matrix}
In agreement with the PCG scheme, the input to the BDDC preconditioner is the residual force vector on the skeleton. As shown in the next subsection, the first step of the preconditioner is to distribute this vector to each cell by means of the extension/restriction operator $\tilde{\mathcal{R}}_U$. In doing so, it is important that the sum of the extended contributions equals one and that the forces are distributed proportionally to the local stiffness associated with their DoFs. For this purpose, the diagonal matrix $\tilde{D}$ is introduced, constructed as
\begin{equation}
    \tilde{D} = \mathrm{blockdiag}\left(D^{(1)},\dots,D^{(n_c)}\right),
\end{equation}
where, for the $i$-th cell, the $j\,j$-th entry of the diagonal matrix $D^{(i)}$, corresponding to the node $X_l$, is defined as
\begin{equation}
    D^{(i)}_{[j\,j]} = \frac{K_{[j\,j]}^{(i)}}{K_{l}^{tot}},
\end{equation}
where $K_{l}^{tot}$ denotes the sum of all stiffness contributions associated with node $X_l$ from the cells sharing that node.

\subsubsection{The preconditioner}
The BDDC preconditioner, denoted by $S_{BDDC}^{-1}$, provides an approximation of the inverse of the global Schur complement. Its action decomposes the correction of the residual into two contributions acting on distinct $\tilde{S}$-orthogonal subspaces. The first contribution is a \emph{fine correction}: independent local Neumann-type problems are solved in each subdomain. Because such local Neumann problems are in general not invertible, the correction is computed in the subspace of functions that vanish at the coarse DoFs (i.e., $\ker(C^{(i)})$), which is enforced through local saddle-point systems (see Algorithm \ref{alg: Local solves a}). The second contribution is a \emph{coarse correction}: the residual is projected onto the subspace spanned by the coarse basis functions $\tilde{\Psi}$, and a global coarse system involving $S_C = \tilde{\Psi}^\top \tilde{S} \tilde{\Psi}^\top$ is solved. This coarse solve balances the local fine corrections, hence the name of the method. The constraints thus serve a dual purpose: they ensure the invertibility of local problems and define a coarse basis that allows the solver to scale effectively as the problem size grows. While the exact solution is globally continuous, the approximate correction obtained from the direct sum of these two subspaces is only guaranteed to be continuous at the coarse DoFs. To restore global continuity across the entire interface skeleton, the partially discontinuous correction is projected back onto the continuous space using the restriction operator $\tilde{\mathcal{R}}_U$ and the diagonal scaling matrix $\tilde{D}$. This weighting is essential for ensuring the robustness and convergence of the method.

The formal expression of the preconditioner reads
\begin{equation} \label{eq: BDDC preconditioner}
    S_{BDDC}^{-1} = \tilde{\mathcal{R}}_U^\top \tilde{D}^\top \tilde{\Psi} S_C^{-1} \tilde{\Psi}^\top \tilde{D} \tilde{\mathcal{R}}_U + \sum_{i=1}^{n_c}
    \begin{bmatrix} {D^{(i)}}{\mathcal{R}}_U^{(i)}\\0\end{bmatrix}^\top
    \begin{bmatrix} S^{(i)} & C^{(i)^\top} \\C^{(i)} & 0 \end{bmatrix}^{-1}
    \begin{bmatrix} {D^{(i)}}{\mathcal{R}}_U^{(i)}\\0  \end{bmatrix} \;.
\end{equation}
In this expression, the first term corresponds to the coarse correction and the second term to the fine correction. The application of $S_{BDDC}^{-1}$ is carried out through a sequence of matrix-vector multiplications and local solves, as detailed in Algorithm \ref{alg: SBDDC multiplication}. While the local problems are fully decoupled and can be solved in parallel, the coarse problem requires the solution of a global system. Since $S_C$ is explicitly assembled, the corresponding linear systems are solved using a direct method based on a parallel Cholesky factorization \added[id=R2]{(specifically, the MUMPS solver \cite{amestoy2001mumps}, accessed through PETSc)}.

\section{Solver accelerations}\label{sec: SPEED}
The BDDC method described in the previous section requires assembling the local stiffness matrix for each cell. However, as demonstrated in Section \ref{sec: RES}, this assembly step is relatively slow and represents the main computational bottleneck of the solver. This section outlines the strategies adopted in this work to accelerate the assembly process.

First, in Section \ref{ssec: Fast Assembly} we employ the fast assembly technique introduced in \cite{hirschler2022}. This approach separates the contribution of the cell mapping, which determines the integrands of Equations \eqref{eq: ELA virtual works}, from the contribution arising from the implicit description of the trimmed domain, which affects only the quadrature rule. As a result, the integrand functions obtained no longer depend on the geometric map. Second, in Section \ref{ssec: ROM} a reduced order model (ROM) is constructed for these integrand functions to avoid their computationally expensive evaluation. This process yields a surrogate model that, while specific to a given level-set, is valid for any geometric mapping \added[id=R2]{for which the pulled-back coefficients satisfy the smoothness conditions required by the fast assembly technique}. The final assembly process for the stiffness matrix is then summarized in Algorithm \ref{alg:local_stiffness}.

\subsection{Fast assembly} \label{ssec: Fast Assembly}

Equation \eqref{eq: ELA constitutive}, introduces the terms that combine the constitutive equation and the external forces with the differential geometry quantities associated to the map. \replaced[id=R2]{Such quantities, namely, the equivalent constitutive tensor $\hat{C}_{i_1j_1i_2j_2}$, the equivalent body force $\hat{\bm{b}}$, and the equivalent traction $\hat{\bm{t}}$ are smooth on each cell, provided that the geometric map is regular (non-vanishing Jacobian) on the closure of the cell and that the material coefficients and loads are smooth cell-wise. Under these conditions, they can be interpolated at the cell level with high accuracy, resulting in}{Assuming sufficient smoothness, these quantities can be interpolated at the cell level, resulting in}
\begin{subequations}\label{eq: FAST interpolation}
\begin{align}
    \hat{C}_{i_1 j_1 i_2 j_2}(\xi_1,\xi_2) &\approx \mathscr{L}_k(\xi_1,\xi_2) \mathrm{C}_{i_1 j_1 i_2 j_2 k}\;, \\
    \hat{b}_{i}(\xi_1,\xi_2) &\approx \mathscr{L}_k(\xi_1,\xi_2) \mathrm{B}_{ik}\;, \label{eq: FAST b interpolation} \\
    \hat{t}_{i}(\xi_1,\xi_2) &\approx \mathscr{L}_k(\xi_1,\xi_2) \mathrm{T}_{ik}\;, \label{eq: FAST t interpolation}
\end{align}
\end{subequations}
being $\hat{b}_{i}$ and $\hat{t}_{i}$ the components of $\hat{\bm{b}}$ and $\hat{\bm{t}}$, respectively, $\mathscr{L}_k(\xi_1,\xi_2)$ the Lagrangian polynomials of degree $q$ associated with the interpolation nodes, selected also here as the Gauss-Lobatto-Legendre (GLL) tensor-product points, and $\mathrm{C}_{i_1 j_1 i_2 j_2 k}$, $\mathrm{B}_{ik}$, and $\mathrm{T}_{ik}$ are the values of $\hat{C}_{i_1 j_1 i_2 j_2}$, $\hat{b}_{i}$, and $\hat{t}_{i}$, respectively, at the interpolation nodes. The components of the test ($v_{i_1}$) and trial ($u_{i_2}$) functions are first approximated using the basis functions $\mathscr{B}_k$:
\begin{subequations}
\begin{align}
     v_{i_1}(\xi_1,\xi_2) = \mathscr{B}_{k_1}(\xi_1,\xi_2) \mathrm{V}_{i_1 k_1} \\
     u_{i_2}(\xi_1,\xi_2) = \mathscr{B}_{k_2}(\xi_1,\xi_2) \mathrm{U}_{i_2 k_2}
\end{align}    
\end{subequations}
where $\mathrm{V}_{ik}$ and $\mathrm{U}_{ik}$ are the degrees of freedom for the test and trial functions, respectively. Substituting this approximation into Equation \eqref{eq: ELA virtual works parametric}, the expressions for the virtual work of the internal and external forces at the cell level are modified as follows:
\begin{subequations}
\begin{align}
     &\mathfrak{L}_{int}^{(i_c)}\approx  \mathrm{V}_{i_1 k_1} \left[\left( \int_{\hat{\Omega}^{(c)}} \fpd{\mathscr{B}_{k_1}}{\xi_{j_1}}
     \fpd{\mathscr{B}_{k_2}}{\xi_{j_2}}\mathscr{L}_{k_3}
     \,\dd\hat{\Omega}  \right) \mathrm{C}_{i_1 j_1 i_2 j_2 k_3}\right] \mathrm{U}_{i_2 k_2}\;,\label{eq: SPEED L int}\\
     &\mathfrak{L}_{ext}^{(i_c)}\approx \mathrm{V}_{i_1 k_1} \left[\left( \int_{\hat{\Omega}^{(i_c)}} \mathscr{B}_{k_1}\mathscr{L}_{k_3}\,\dd\hat{\Omega}\right) \mathrm{B}_{i_1k_3} + \left( \int_{\partial\hat{\Omega}^{(i_c)}} \mathscr{B}_{k_1}\mathscr{L}_{k_3}\,\dd\hat{\Omega}\right)\mathrm{T}_{i_1k_3} \right]\label{eq: SPEED L ext}
     \;.
\end{align}
\end{subequations}
The quantities in the square brackets are the entries of the stiffness matrix and the external force vector associated to the cell. However, only the quantities inside the integrals actually depend on the curvilinear coordinates $\xi_1$ and $\xi_2$. If the value of the integral were known, the assembly procedure would be equivalent to computing the tensorial contraction along $j_1$, $j_2$, and $k_3$, which is computationally cheap.

In a nutshell, if the constitutive relationship, the external forces, and the differential geometry, can be combined and accurately approximated as polynomials, then the integration could be precomputed with some efficient technique and the actual stiffness matrix and external forces vectors obtained by a fast tensorial product with the coefficients of such polynomials. This procedure, introduced in \cite{hirschler2022fast}, is referred to as the fast assembly technique.

It should be noted that the interpolation in Equations \eqref{eq: FAST b interpolation} and \eqref{eq: FAST t interpolation} is not always applicable. \added[id=R2]{The fast assembly technique requires the geometric map to be regular, i.e., with non-vanishing Jacobian, on the closure of each cell, and the material coefficients, the body forces, and the tractions to be smooth on each cell. Therefore, discontinuities aligned with cell boundaries are admissible. Under these conditions, and in particular for the spline maps and constant material coefficients employed in this work, the pulled-back quantities $\hat{C}_{i_1j_1i_2j_2}$, $\hat{\bm{b}}$, and $\hat{\bm{t}}$ are analytic within each cell, and their polynomial interpolation of degree $q$ converges exponentially. Furthermore, the interpolation error can be evaluated a priori, by sampling the difference between the exact and interpolated quantities, so that $q$ can be selected adaptively to meet a prescribed tolerance; the resulting consistency error is controlled through Strang's first lemma \cite{hirschler2022fast,mantzaflaris2015}.}

\added[id=R2]{When these regularity requirements are violated, for instance, when the geometric map approaches degeneracy or in the presence of discontinuous material fields or non-smooth loads, the interpolation error ceases to decrease at high order, and a full quadrature of the affected term is required. We remark, however, that the impact of non-smooth data is not limited to the fast assembly step. In fact, spectral elements themselves lose their high-order approximation properties when the solution lacks sufficient regularity.} \deleted[id=R2]{and requires a sufficient degree of smoothness in the functions representing the forces to ensure accuracy. In critical cases, a full assembly of the force term may be required. However, this remains computationally less expensive than assembling the contribution of the cell to the stiffness matrix through full integration.}

Regarding the stabilization term in Equation \eqref{eq: MOD stabilization}, it can also be assembled using a fast assembly procedure. In particular,
\begin{equation}
    \mathfrak{L}_{sta} \approx  \mathrm{V}_{i_1 k_1} \left[\left( \int_{\hat{\Omega}^{(i_c)}_\Pi} \fpd{\mathscr{B}_{k_1}}{\xi_{j_1}}
     \fpd{\mathscr{B}_{k_2}}{\xi_{j_2}}\mathscr{L}_{k_3}
     \,\dd\hat{\Omega}  \right) \rho\mathrm{C}_{i_1 j_1 i_2 j_2 k_3}\right] \mathrm{U}_{i_2 k_2}\;.
\end{equation}
Exploiting the identity \begin{equation}\int_{\hat{\Omega}^{(i_c)}_\Pi}\left(\bullet\dd\hat{\Omega}\right)=\int_{\hat{\Pi}^{(i_c)}}\left(\bullet\dd\hat{\Omega}\right)-\int_{\hat{\Omega}^{(i_c)}}\left(\bullet\dd\hat{\Omega}\right),
\end{equation}
the above expression can be rewritten as
\begin{equation}
    \mathfrak{L}_{sta} \approx  \mathrm{V}_{i_1 k_1} \left[\left( \int_{\hat{\Pi}^{(i_c)}} \fpd{\mathscr{B}_{k_1}}{\xi_{j_1}}
     \fpd{\mathscr{B}_{k_2}}{\xi_{j_2}}\mathscr{L}_{k_3}
     \,\dd\hat{\Omega} - \int_{\hat{\Omega}^{(i_c)}} \fpd{\mathscr{B}_{k_1}}{\xi_{j_1}}
     \fpd{\mathscr{B}_{k_2}}{\xi_{j_2}}\mathscr{L}_{k_3}
     \,\dd\hat{\Omega}  \right) \rho\mathrm{C}_{i_1 j_1 i_2 j_2 k_3}\right] \mathrm{U}_{i_2 k_2}\;.
\end{equation}
The first integral in this expression can be precomputed and stored, as it does not depend on the trimming configuration. The second integral is formally equivalent to the one appearing in Equation \eqref{eq: SPEED L int}, which can therefore be reused here.

\subsection{Reduced order modeling} \label{ssec: ROM}
In this section, we describe how to construct the ROM in order to avoid the need for full numerical quadrature. In this regard, to make it adaptable to different mappings, the models are constructed on the fast assembly integrals introduced in Section \ref{ssec: Fast Assembly}. The procedure is presented here for the tensor associated with the stiffness term, as in Equation \eqref{eq: SPEED L int}, but it is adopted also for the integrals related to the external forces as in Equation \eqref{eq: SPEED L ext}.
Let us denote as $\mathcal{I}$ the tensor whose components are
\begin{equation}
\mathcal{I}^{k_1k_2k_3}_{j_1j_2} = \fpd{\mathscr{B}_{k_1}}{\xi_{j_1}}
     \fpd{\mathscr{B}_{k_2}}{\xi_{j_2}}\mathscr{L}_{k_3}.
\end{equation}
The vector containing the integrals of the components of $\mathcal{I}$ is denoted as:
\begin{equation}\label{eq: MDEIM I}
    I(\bm{\mu}) = \int_{\hat{\Omega}(\bm{\mu})}\mathrm{vec}(\mathcal{I})
     \,\dd\hat{\Omega}
\end{equation}
where $\mathrm{vec}(\bullet)$ stands for the vectorization operation of the tensor $\bullet$.
Here, the dependence of $I$ from the level-set threshold parameters $\bm{\mu}= \{\mu_{1}, \mu_{2}, \mu_{3}, \mu_{4}\}$, arises from the variation of the trimmed domain $\hat{\Omega}$. This section revolves around the strategy to avoid explicitly computing the integrals in Equation \eqref{eq: MDEIM I} by means of an offline built reduced order model. 

In particular, the goal is to construct an approximation of the form:
\begin{equation} \label{eq: MDEIM I approximation}
    I(\bm{\mu}) \approx U_r I_r(\bm{\mu})\;,
\end{equation}
where $U_r\in\mathbb{R}^{n_i\times n_r}$, and $I_r\in\mathbb{R}^{n_r}$. The columns of $U_r$ constitute the basis vectors of the reduced order model, and the reduced order vector $I_r$ contains the components of $I$ with respect to such basis. The length of the vector $I$ is denoted as $n_i$, while the number of basis vector\added[id=R2]{s} is denoted as $n_r$.

\subsubsection{Clustering of the threshold parameters space}
The construction of the ROM approximation begins with the definition of the threshold parameter space, $\boldsymbol{\mu} \in M$. This space is defined as the hypercube $M=[\mu_\mathrm{min}, \mu_\mathrm{max}]^d$, where $\mu_\mathrm{min}$ and $\mu_\mathrm{max}$ are the respective lower and upper bounds, and $d$ the number of parameters, four in this case. This domain is then partitioned into clusters by subdividing each interval into $n_k$ subintervals, resulting in a total of $n_K = (n_k)^4$ cluster domains denoted as $M_{i_k}$. The approximation described in Equation \eqref{eq: MDEIM I approximation} is then performed independently within each cluster. In other words, a distinct ROM is constructed for each cluster, following the procedure outlined below. 

To keep the notation concise, no additional index is introduced to explicitly indicate the cluster to which each quantity belongs, as this is not expected to generate ambiguity. It should be noted, however, that when evaluating functions of $\bm{\mu}$, the first step consists in identifying the cluster to which the corresponding threshold parameter vector belongs.

\subsubsection{Reduced order basis}
The basis $U_r$ is constructed in an offline phase. First, a set of $n_s$ representative level-set threshold parameter vectors, ${\bm{\mu}^{{i_s}}}$, is generated within the associated cluster domain. To ensure a representative sample with a limited number of elements, these parameter sets are selected using the Latin hypercube sampling method \cite{mckay2000}. For each of these samples, the corresponding snapshot integrand, $I(\bm{\mu}^{{i_s}})$, is then computed exactly using the accurate quadrature procedure described in Section \ref{sec: MOD unfitted}. These snapshots are then collected in the matrix
\begin{equation}
\mathfrak{I} = \left[I(\bm{\mu}^{\{1\}}),\; I(\bm{\mu}^{\{2\}}),\; \cdots,\; I(\bm{\mu}^{\{n_s\}})\right]
\in \mathbb{R}^{n_i \times n_s}.
\end{equation}
A randomized singular values decomposition is then applied on this matrix and truncated to $n_v$ singular values, with $n_v>n_r$ to account for the loss of accuracy on the trailing singular vectors, leading to the following approximated expression
\begin{equation}
    \mathfrak{J} \approx U_v \Sigma_v V_v^\Tr,
\end{equation}
where $U_v\in\mathbb{R}^{n_i\times n_v}$ is the matrix containing as columns the left singular vectors, $\Sigma_v\in\mathbb{R}^{n_v\times n_v}$ contains on the main diagonal the singular values, and $V_v\in\mathbb{R}^{n_s\times n_v}$ contains as columns the right singular vectors. The vectors of the basis $U_r$ in Equation \eqref{eq: MDEIM I approximation} are selected as the first $n_r$ columns of $U_v$.

\subsubsection{Coefficients of the reduced model} \label{sssec: MDEIM}
With the basis for the reduced-order model constructed, the non-linear function $I_r(\bm{\mu})$ is then approximated using the Matrix Discrete Empirical Interpolation Method (MDEIM) \cite{negri2015}.

Firstly, the so-called \emph{magic points} are selected \cite{barrault2004}. The number of magic points is equal to the size of the reduced basis $n_r$, and correspond to positions within the basis vectors. As such, they are integers in the set $m_{i_r}\in\{1,\dots,n_i\}$, with $i_r=1,\dots,n_r$. How the magic points are computed is detailed in Algorithm \ref{alg: MDEIM magic points}. In a nutshell, the algorithm works by iterating over the basis vectors; at each iteration, the entry that is most orthogonal to the previous vectors is selected.

The magic points are used to compute the reduced vector $I_r$ from the full vector $I$. In particular, $I_r$ contains coefficients chosen so that the reduced approximation matches the entries of $I$ exactly at the magic points. This condition is expressed as
\begin{equation}
    U_{r\,i}^{[m_j]} \, I_r^{[j]} = I^{[m_j]}, \qquad i,j = 1,\dots,n_r ,
\end{equation}
which corresponds to solving an $n_r \times n_r$ linear system, where $U_{r\,i}^{[m_j]}$ denote the $m_j$-th component of the $i$-th vector of the reduced-order basis, and $I_r^{[j]}$ and $I^{[m_j]}$ are the $j$-th and the $m_j$-th components of $I_r$ and $I$, respectively.

However, computing $I_r$ exactly requires access to the magic-points entries of the non-reduced vector $I$, which in turn would require a full integration process, which the proposed approach aims to avoid in the online phase. To overcome this issue, the vector 
$I_r(\boldsymbol{\mu})$ is computed exactly in an offline phase for a set of interpolatory 
snapshots. This sampling, which does not coincide with that used for the 
singular value decomposition, has a tensor-product structure within the associated cluster domain $M_{i_k}$, and follows the distribution of GLL points. 

The function $I_r(\boldsymbol{\mu})$ is then obtained by interpolating these values, \added[id=R2]{using Lagrangian polynomial of degree $r$}.  In this work, the interpolation is based on standard Lagrange polynomials defined over GLL nodes, which are preferred over radial basis functions \cite{powell1992theory}. This choice is motivated by the aim of minimizing the impact of the interpolation error on the overall accuracy of the MDEIM model. For problems with a larger number of parameters, Lagrangian interpolation may become too expensive in terms of both memory and computational cost. Alternative approaches, such as sparse grids \cite{bugartz2004},  might be used in such cases.

\subsection{\added[id=R2]{Offline cost and speedup}}\label{ssec: ROM cost}
\added[id=R2]{This section quantifies the offline and online costs of the ROM-based acceleration strategies introduced above. The offline phase involves training a surrogate model for each level-set family. The online phase consists in evaluating the surrogate model to assemble each cell stiffness matrix, replacing the expensive numerical quadrature.}

\textbf{\added[id=R2]{Offline training:}}
\added[id=R2]{The training is performed independently for each level-set geometry. The first step is to define the interval within which the threshold parameters are allowed to vary. This interval is chosen to exclude geometrically degenerate configurations, such as disconnected cells. For the Schwarz Diamond, two models are trained over the parameter spaces $\mathbb{P}_1^\mathrm{SD} = [0.1, 0.9]^4$ and $\mathbb{P}_2^\mathrm{SD} = [0.1, 1.0]^4$. For the Schoen IWP, the two corresponding spaces are $\mathbb{P}_1^\mathrm{IWP} = [-2.5, 2.5]^4$ and $\mathbb{P}_2^\mathrm{IWP} = [-2.5, 3.0]^4$. For both geometries, the latter space is defined to include the fully solid cell, which however comes at the cost of a slight reduction in accuracy (see Section \ref{sssec: RES multiple}).}

\added[id=R2]{Each model is built for a polynomial degree $p=8$ and a fast assembly interpolation degree $q=2$. The parameter space is partitioned into $n_K = n_k^4 = 16$ clusters, $n_k=2$ for each parametric direction. Within each cluster, $n_s = 100$ fully integrated snapshots are generated via Latin hypercube sampling, and a randomized SVD is performed, retaining $n_v = 50$ singular vectors of which the first $n_r = 40$ form the reduced basis.}

\added[id=R2]{The MDEIM coefficient fitting requires, per cluster, $(r+1)^4 = 7^4 = 2{,}401$ additional fully integrated evaluations on a tensor-product GLL grid, for a total of $n_K \times 7^4 = 38{,}416$ evaluations across all clusters. These parameter values were tuned based on the numerical results presented in Section \ref{sec: RES}.
For the ROM corresponding to $\mathbb{P}_1^\mathrm{IWP}$, the total training time is approximately 9 hours (using 30 cores on a workstation with 4 Intel Xeon Gold 6148 processors and 1 TB of memory), and the total disk-storage footprint of the trained ROM data is 2.16 GB. Comparable training times were observed for the other geometries considered in this work.   These figures include the model for the stiffness, the model for the mass matrix (not utilized in this contribution), and the model used to impose Neumann boundary conditions. The difference in the speedup factor arises from the higher geometric complexity of the Schoen IWP-based geometries, resulting in a larger number of quadrature points.}

\begin{table}[ht]\centering
\caption{\added[id=R2]{Online assembly cost per cell, with and without the ROM, for two level-set families. Timings are averaged over 1000 cells.}}
\label{tab: ROM online}
\renewcommand{\arraystretch}{1.2}
\begin{tabular}{lll}\hline
\textbf{Quantity} & \textbf{Schwarz Diamond} & \textbf{Schoen IWP} \\ \hline
Full quadrature assembly & 499 ms & 870 ms \\
Assembly time with ROM   & 2.3 ms & 2.5 ms \\
Speedup factor           & 218$\times$ & 342$\times$ \\ \hline
\end{tabular}
\renewcommand{\arraystretch}{1.0}
\end{table}

\textbf{\added[id=R2]{Online assembly and amortization threshold:}}
\added[id=R2]{Once the ROM is available, the assembly cost per cell is reduced to a Lagrangian interpolation of the MDEIM reduced vector, followed by a multiplication with the reduced basis. The assembly times, with and without the ROM, for two level-set families and the settings described in Section~\ref{sssec: RES single cell}, are reported in Table~\ref{tab: ROM online}.}

\added[id=R2]{The amortization threshold is defined as the total number of cell assembly evaluations required for the cumulative online savings to offset the offline training cost. Denoting the training time as $T_{\mathrm{off}}$, the per-cell assembly time without ROM as $T_{\mathrm{full}}$, and with ROM as $T_{\mathrm{rom}}$, the threshold is}
\begin{equation}
    n_{\mathrm{amm}} = \left\lceil \frac{T_{\mathrm{off}}}{T_{\mathrm{full}} - T_{\mathrm{rom}}} \right\rceil.
\end{equation}
\added[id=R2]{As an example, for the ROM model $\mathbb{P}_1^\mathrm{IWP}$, the break-even point is reached after approximately $38 000$ cell evaluations, beyond which the offline training cost is amortized. This threshold is largely independent of the geometry, as the offline training cost is primarily determined by the fixed number of fully integrated cells used to construct the model.}

\added[id=R2]{On this regard, it is worth noting that the number of TPMS geometries used as lattice micro-structures is limited. The training therefore represents a one-time community investment rather than a per-user burden. Consistent with this philosophy, the precomputed ROM models used in this work are publicly available on Zenodo~\cite{flash_rom_zenodo}. These models are directly reusable in any simulation involving the same level-set families: indeed, they were computed once and then reused, without any retraining, in all the analyses of Section~\ref{sec: RES}.}

\subsection{Summary of the computational workflow}\label{ssec: workflow}
\begin{figure}[htbp]
\centering
\begin{subfigure}[b]{0.48\textwidth}\centering\includegraphics[width=\textwidth]{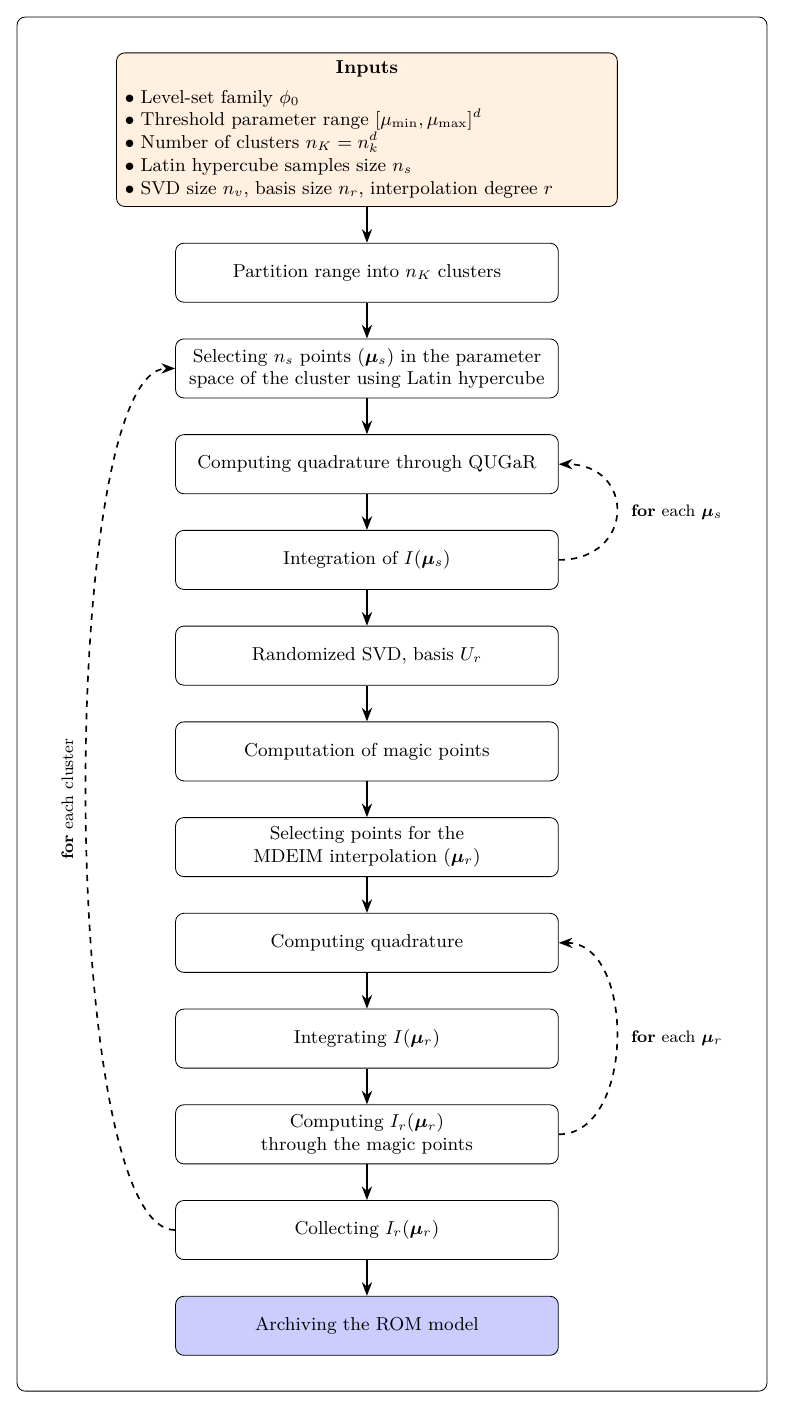}\caption{}\label{fig:workflow_off}\end{subfigure}\hfill
\begin{subfigure}[b]{0.48\textwidth}\centering\includegraphics[width=\textwidth]{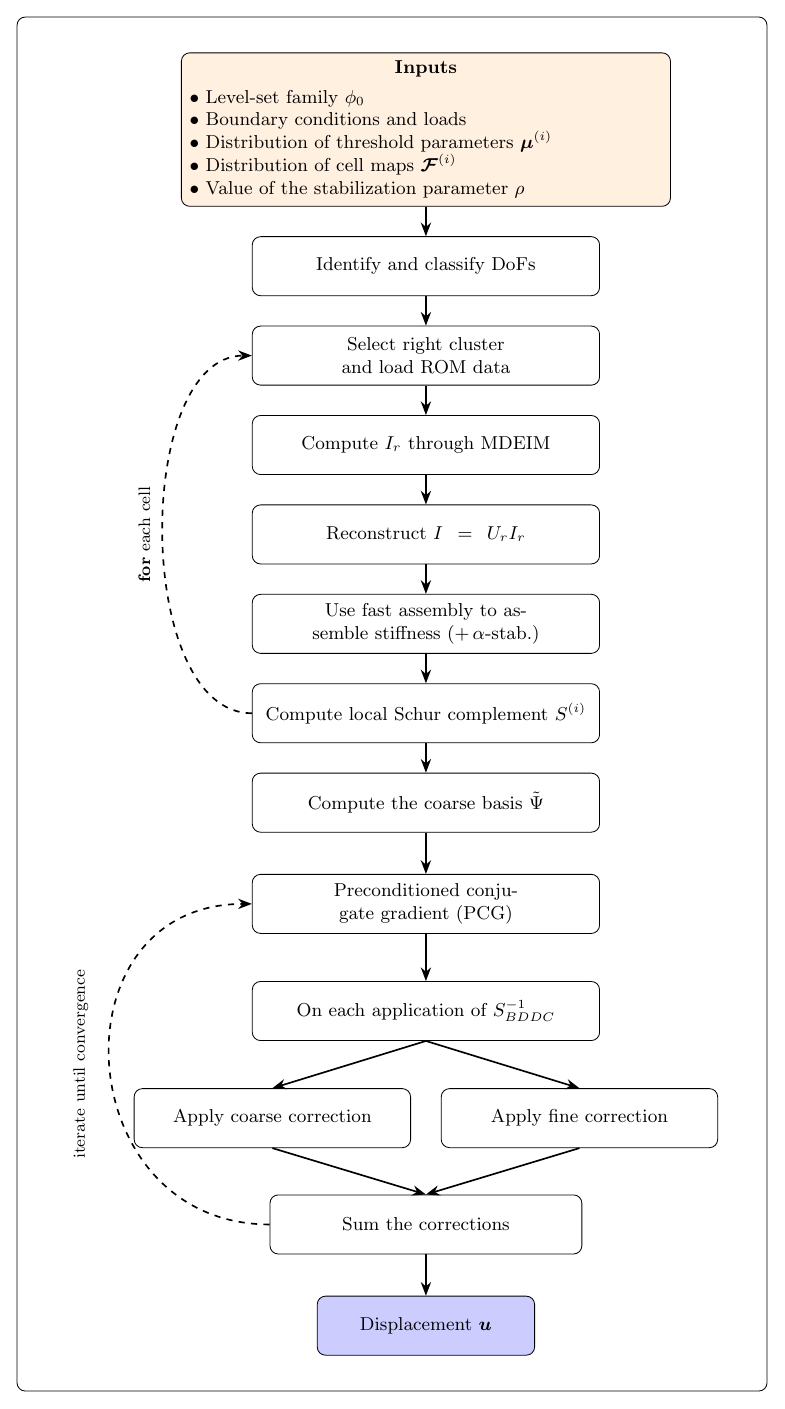}\caption{}\label{fig:workflow_on}\end{subfigure}
\caption{\added[id=R1]{Computational workflow of the ROM-based BDDC solver, read top to bottom from the input to the output box. (a) Offline ROM-training pipeline, built once per level-set family. (b) Online per-cell assembly and BDDC solve, performed for every analysis. The individual steps, with the corresponding algorithms, are detailed in the text.}}
\label{fig:workflow}
\end{figure}

\added[id=R1]{Having introduced the individual ingredients and quantified their offline and online cost, we now summarize how they combine into a single offline--online pipeline, illustrated in Figure~\ref{fig:workflow}. The two panels show the respective workflows from the input box to the output box: Figure~\figref{fig:workflow}{a} traces the one-time offline training of the ROM, while Figure~\figref{fig:workflow}{b} traces the online assembly and the subsequent BDDC solve. The two stages are deliberately decoupled, and trimmed quadrature is confined entirely to the offline stage.}

\textbf{Offline phase:} Following Figure \figref{fig:workflow}{a}, the range of the threshold parameters is partitioned into the $n_K=(n_k)^4$ clusters, and every subsequent step is repeated independently for each cluster. Within a cluster, snapshots of the fast-assembly tensor $I(\bm{\mu})$ are computed by full trimmed quadrature (Section~\ref{sec: MOD unfitted}) at the $n_s$ parameter samples $\bm{\mu}_s$ obtained by Latin hypercube sampling, and a truncated randomized SVD of these snapshots yields the reduced basis $U_r$. The MDEIM magic points are then extracted (Algorithm~\ref{alg: MDEIM magic points}), and the interpolation coefficients are fitted from a second, tensor-product set of samples $\bm{\mu}_r$, which likewise require full quadrature. Once all clusters have been processed, the trained model is archived for reuse (Section~\ref{ssec: ROM cost}). We note once again that the resulting trimming-based ROM is geometry-specific but mapping-independent.

\textbf{Online phase:} The online stage corresponds to Figure~\figref{fig:workflow}{b}. After the degrees of freedom have been identified and classified as either Dirichlet, internal, or skeleton DoFs, the solver loops over the cells: given the threshold parameters $\bm{\mu}^{(i)}$ of a cell, the relevant cluster is selected and the reduced vector $I_r$ is obtained by MDEIM interpolation, from which the fast-assembly tensor is reconstructed as $I\approx U_r I_r$ (Algorithm~\ref{alg: mdeim}). This tensor is contracted with the mapping coefficients corresponding to $\bm{\mathcal{F}}^{(i)}$ to form the local stiffness matrix, to which the $\alpha$-stabilization is added (Algorithm~\ref{alg:local_stiffness}), the local Schur complement $S^{(i)}$ is then formed. Once every cell has been processed, it is computed the coarse basis $\tilde{\Psi}$ and the preconditioned conjugate gradient is initiated. When requested by the iterative solver scheme, application of preconditioner $S_{BDDC}^{-1}$ is carried out by superimposition of a coarse and a fine correction (Algorithm~\ref{alg: SBDDC multiplication}), the first relying on a global direct solve, and the latter relying on the local saddle-point solves of Algorithm~\ref{alg: Local solves a}, and the iteration proceeds until the prescribed tolerance is met, returning the displacement field $\bm{u}$.

\added[id=R1]{Read together, the two panels make explicit the migration of the dominant cost. In the unaccelerated baseline, the per-cell trimmed quadrature dominates the setup phase (Figure~\figref{fig:test_4}{b}); the acceleration layers move this cost to the one-time offline training of Figure~\figref{fig:workflow}{a}, so that the online assembly of Figure~\figref{fig:workflow}{b} becomes negligible and the BDDC setup and solve times become comparable (Figure~\ref{fig:test_1}). The ingredients remain simultaneously effective provided that the pulled-back integrands are smooth, as required by the fast assembly (Section~\ref{ssec: Fast Assembly}); that the threshold parameters lie within the admissible range used to train the ROM and to exclude degenerate cells; and that the $\alpha$-stabilization keeps the conditioning of strongly trimmed cells under control, which in turn supports the robustness of the BDDC coarse space. }

\section{Numerical Results}

\label{sec: RES}
In this section, several numerical examples are presented to assess the performance of the proposed computational framework. The section is divided into two parts. Initially, tests are conducted on simple geometries to evaluate the impact of the  various approximations introduced by the method, namely the $\alpha$-stabilization, the fast assembly, and the reduced-order modeling, and to guide the selection of the relevant parameters. The performance is also assessed in terms of iteration count and total analysis time. In the second part, the proposed approach is applied to more realistic geometries, namely a sandwich wing and a lattice wrench. 

It is further recalled that, within the presented BDDC framework, there is a one-to-one correspondence between cells and subdomains. Accordingly, the two terms are used interchangeably throughout the section.

The method presented in this paper is implemented in FLAS$_h$, an open-source Python library \cite{flash_github}. \added[id=R2]{In particular, the repository includes dedicated example scripts that reproduce each figure and timing table presented in this section, along with the corresponding input configurations and solver settings}. Table~\ref{tab: RES reproducibility} lists, for each figure of this section, the script(s) of the FLAS$_h$ repository reproducing it, all referring to version v0.2.0. FLAS$_h$ \deleted[id=R2]{that} relies on PETSc4py \cite{petsc-user-ref,petsc4py} for the parallel solver, NumPy and SciPy for numerical computations, and QUGaR \cite{qugar} for the unfitted discretizations. Precomputed ROM models used in this work are publicly available on Zenodo \cite{flash_rom_zenodo} and are automatically downloaded during the library installation. All the numerical results presented in this section are fully reproducible using FLAS$_h$. The reported tests were performed on a \added[id=R2]{Apple} MacBook Pro \added[id=R2]{with an Apple M4 Pro chip (14-core CPU: 10 performance + 4 efficiency cores; 48 GB unified memory).} \deleted[id=R2]{(14 cores) with 48 GB of RAM} The simulations were run in parallel using 8 \added[id=R2]{performance} cores.

\begin{table}[ht]\centering
\caption{\added[id=R2]{Scripts of the FLAS$_h$ repository (v0.2.0) reproducing the figures of this section.}}
\label{tab: RES reproducibility}
\renewcommand{\arraystretch}{1.2}
\begin{tabular}{cl>{\raggedright\arraybackslash}p{7.5cm}}\hline
\textbf{Figure} & \textbf{Section} & \textbf{Script(s)} \\ \hline
9  & 5.1.1 & \texttt{test\_convergence.py} \\
10 & 5.1.1 & \texttt{test\_stabilization\_1.py} \\
11 & 5.1.1 & \texttt{test\_stabilization\_2.py} (a), \texttt{test\_stabilization\_3.py} (b) \\
13 & 5.1.2 & \texttt{test\_fast\_assembly\_accuracy.py} \\
14 & 5.1.2 & \texttt{test\_rom\_basis.py} \\
15 & 5.1.2 & \texttt{test\_rom\_accuracy.py} \\
16 & 5.1.3 & \texttt{test\_solver\_comparison.py} \\
17 & 5.1.3 & \texttt{test\_acceleration\_efficiency.py} \\
18 & 5.1.3 & \texttt{test\_scalability.py} \\
19 & 5.2.1 & \texttt{example\_wing.py} \\
20 & 5.2.2 & \texttt{example\_wrench\_coarse.py} (a,b), \texttt{wrench\_refined.py} (c,d) \\
21 & 5.2.2 & \texttt{wrench\_refined.py}, \texttt{wrench\_validation.py}, \texttt{wrench\_paper\_summary.py}, \texttt{plot\_wrench\_h\_study.py} \\ \hline
\end{tabular}
\renewcommand{\arraystretch}{1.0}
\end{table}

\subsection{Method validation and benchmarking}
This section \deleted[id=R2]{first} justifies the choice of discretization strategy, \deleted[id=R2]{It then} evaluates the accuracy of each approximation component (fast assembly, ROM, and stabilization), and assesses the computational performance of the overall method.
Throughout this section, the term \emph{baseline} refers to the BDDC method \added[id=R1]{assembled with full integration} without applying any acceleration technique \deleted[id=R2]{n}or stabilization. Additionally, all errors reported in this section are normalized with respect to the corresponding reference quantities.

\subsubsection{\replaced[id=R2]{Accuracy assessment on a single-cell structure}{Discretization strategy}}
\label{sssec: RES single cell}
As a starting point, a simple single-cell geometry is considered.\deleted[id=R2]{The size of the cell is taken to be unitary in both dimensions, such that the mapping of the cell reduces to the identity.}  The threshold parameter $\mu(\boldsymbol{\xi})$ is set \added[id=R2]{to be uniform, with value specified in the tests}. \deleted[id=R2]{uniformly to 0 (same as in Figure \figref{fig: GEO TPMS geometries}(c))}. \added[id=R2]{The geometry is subsequently transformed through the map.}
\added[id=R2]{
\begin{equation*}
    \bm{\mathcal{F}}=\begin{bmatrix} \frac{1}{2}+\left(\xi_1-\frac{1}{2}\right)(1+(a-1)(2\xi_2-\xi_2^2))\\ \xi_2\end{bmatrix}\;,
\end{equation*}
where $a$ is a parameter controlling the severity of the geometric distortion.} The body forces, prescribed displacements, and boundary tractions for this test are derived from the reference manufactured solution:
\begin{equation}\label{eq:manufactured}
\boldsymbol{u}_{ex} = \begin{bmatrix} u_0(2x_2-1)\sin^2{(\pi x_1)}\sin^2{(\pi x_2)}\\
u_0(1-2x_1)\sin^2{(\pi x_1)}\sin^2{(\pi x_2)}
\end{bmatrix}\;.
\end{equation}
It is important to note that, here, Neumann boundary conditions are also applied on the trimmed boundary $\partial\Omega_t$ to enable a consistent comparison with the exact solution, as permitted by the manufactured nature of the problem settings. The Lamé constants of the adopted material are both set unitary. Here, no stabilization or acceleration techniques are employed in the assembly of the stiffness matrices.

\begin{figure}[ht]
\centering
\begin{subfigure}[b]{0.98\textwidth}\centering
\includegraphics[width=\textwidth]{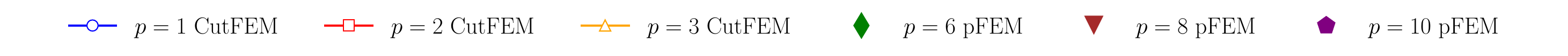}
\end{subfigure}\\
\begin{subfigure}[b]{0.32\textwidth}\centering
\includegraphics[width=\textwidth]{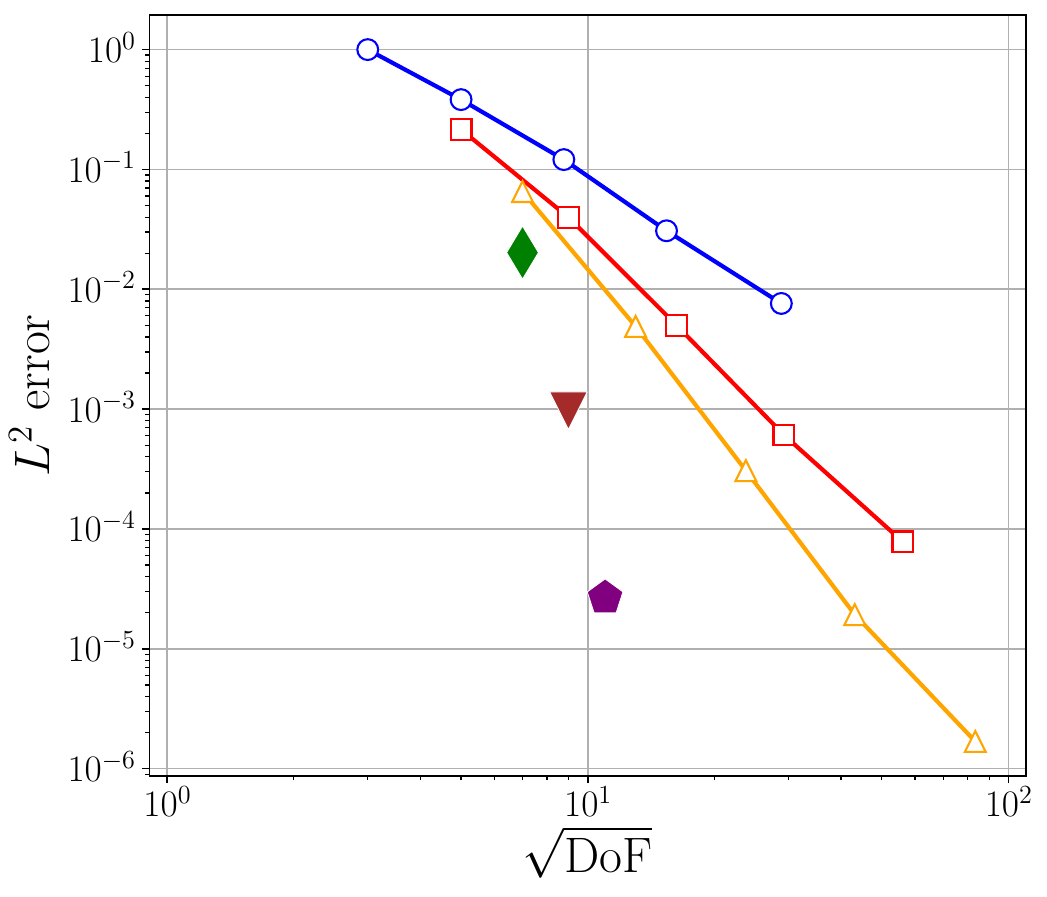}
\caption{}
\end{subfigure}
\begin{subfigure}[b]{0.32\textwidth}\centering
\includegraphics[width=\textwidth]{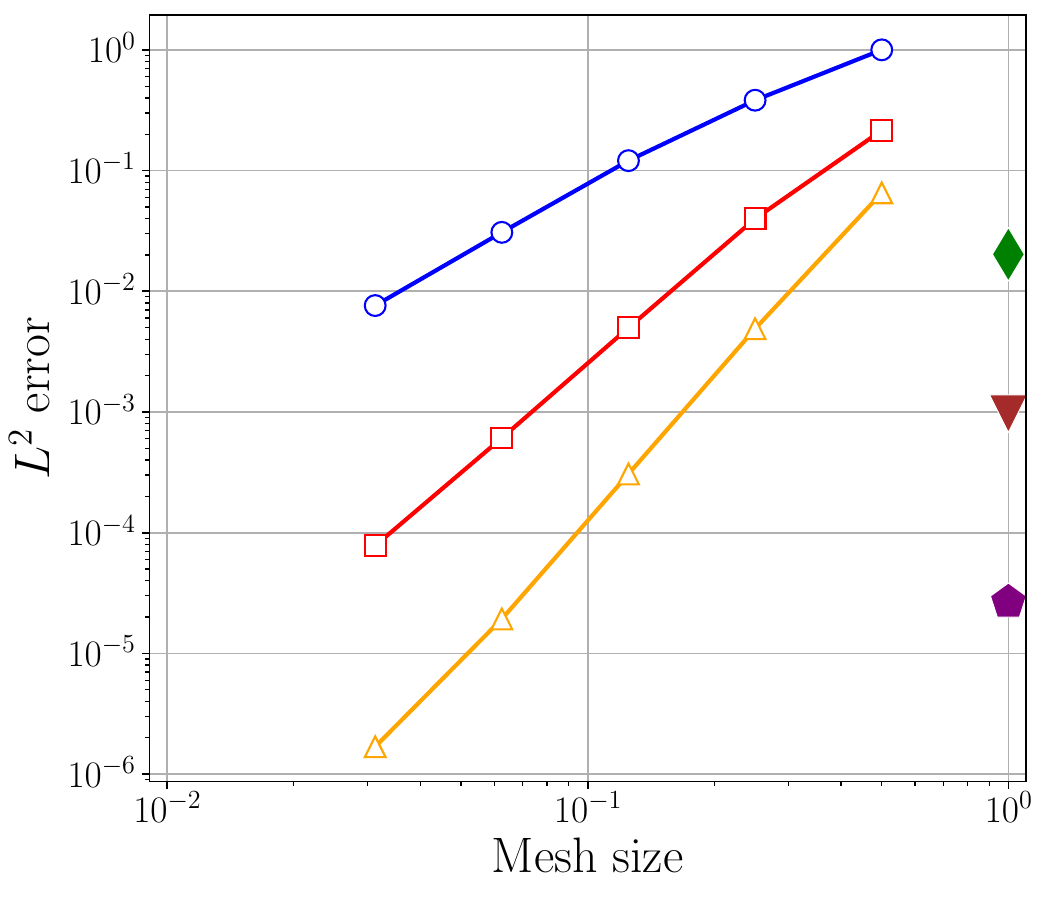}
\caption{}
\end{subfigure}
\begin{subfigure}[b]{0.32\textwidth}\centering
\includegraphics[width=\textwidth]{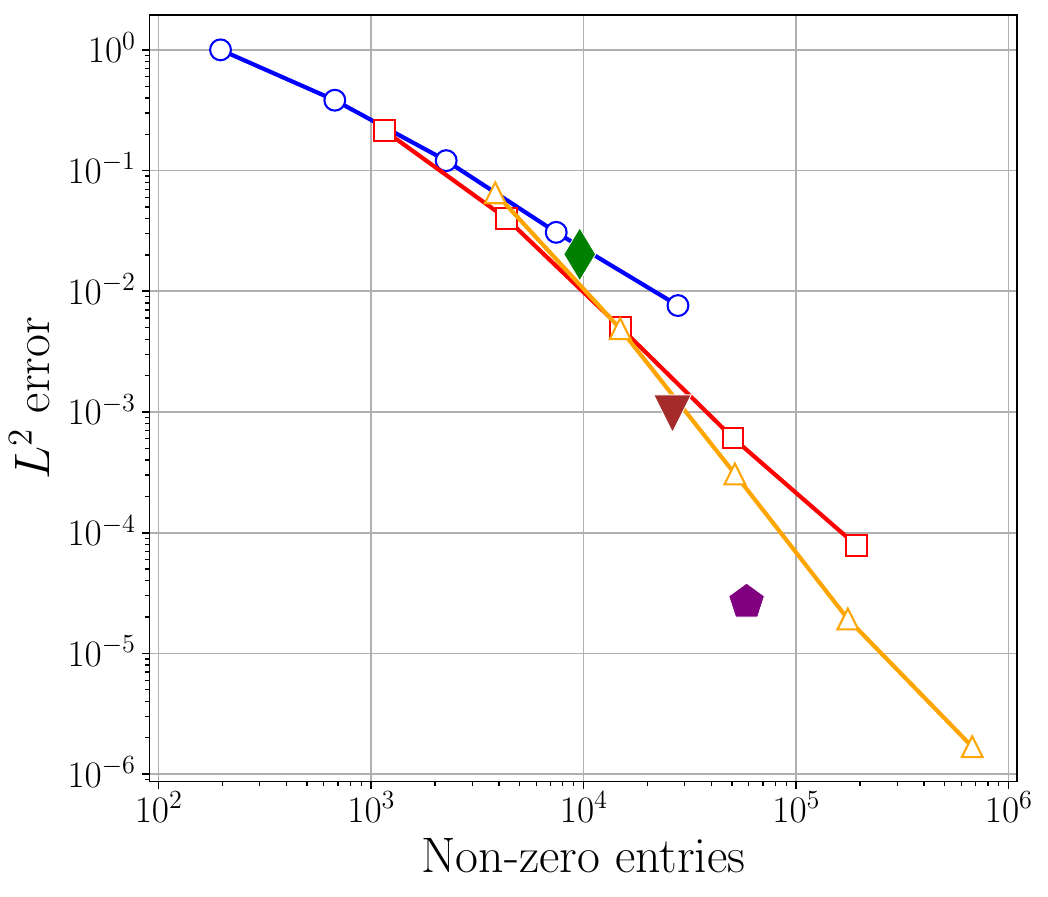}
\caption{}
\end{subfigure}
\caption{Comparison between unfitted p-FEM element method and CutFEM for the manufactured problem described in Section \ref{sssec: RES single cell}. In the left picture the $L^2$ error is plotted against the square root of the number of DoFs. In the center, the $L^2$ error is plotted against the mesh size. Here p-FEM data are vertically aligned since the mesh size in all such cases is unitary. On the right, the $L^2$ error is reported as a function of non-zero entries of the stiffness matrix.}
\label{fig: RES spectral}
\end{figure}

\paragraph{\added[id=R2]{Accuracy of the p-FEM discretization}}
\added[id=R2]{In this test, the adopted level-set is the Schoen FRD shown in Figure~\figref{fig: GEO TPMS geometries}{c}. The threshold parameter is set to zero, and the transformation is taken as the identity ($a=1$). No stabilization is adopted in this test.} Two different discretizations are taken into account. Specifically, either a single unfitted p-FEM element is adopted, or a grid of rectangular elements is employed within a CutFEM framework \cite{burman2015,burman2025cut}. Figure \figref{fig: RES spectral} shows the convergence curves of the $L^2$ error in the displacement field for both discretizations, with respect to the square root of the degrees of freedom, the mesh size, and the number of non-zero entries of the stiffness matrix. 
From these results, it is evident that a single high-degree p-FEM element can achieve the same level of accuracy with a significantly lower number of degrees of freedom. This is particularly apparent for $p = 10$, which outperforms the CutFEM refinement of degree $p = 3$ by requiring more than an order of magnitude fewer degrees of freedom. This advantage becomes less pronounced when considering the number of nonzero entries of the stiffness matrix. In this case, the p-FEM with $p = 10$ requires fewer than four times the number of nonzero entries compared to the corresponding CutFEM curve with $p = 3$. This comparison becomes even more relevant when examining the $p = 6$ and $p = 8$ points for the p-FEM, where CutFEM refinements of order $p = 2$ and $p = 3$ can be more efficient.

However, when adopting a CutFEM approach, shape functions may switch between active and inactive states, which can introduce difficulties in the construction of the ROM. This issue does not arise in the p-FEM case. Consequently, the stiffness matrix to be approximated within the ROM framework depends continuously on the threshold parameters, therefore making the surrogate model more accurate. For this reason, in this work p-FEM elements are adopted to approximate the solution within each cell. Additionally, unless otherwise specified, the polynomial degree is set to $p = 8$, as it provides a suitable compromise between accuracy and computational cost.

\paragraph{\added[id=R2]{Accuracy of the stabilization}}
\begin{figure}[ht]
\centering
\begin{subfigure}[b]{0.75\textwidth}\centering
\includegraphics[width=\textwidth]{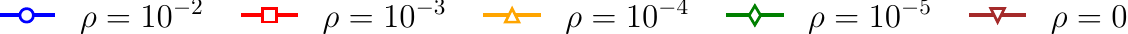}
\end{subfigure}\\
\begin{subfigure}[b]{0.44\textwidth}\centering
\includegraphics[width=\textwidth]{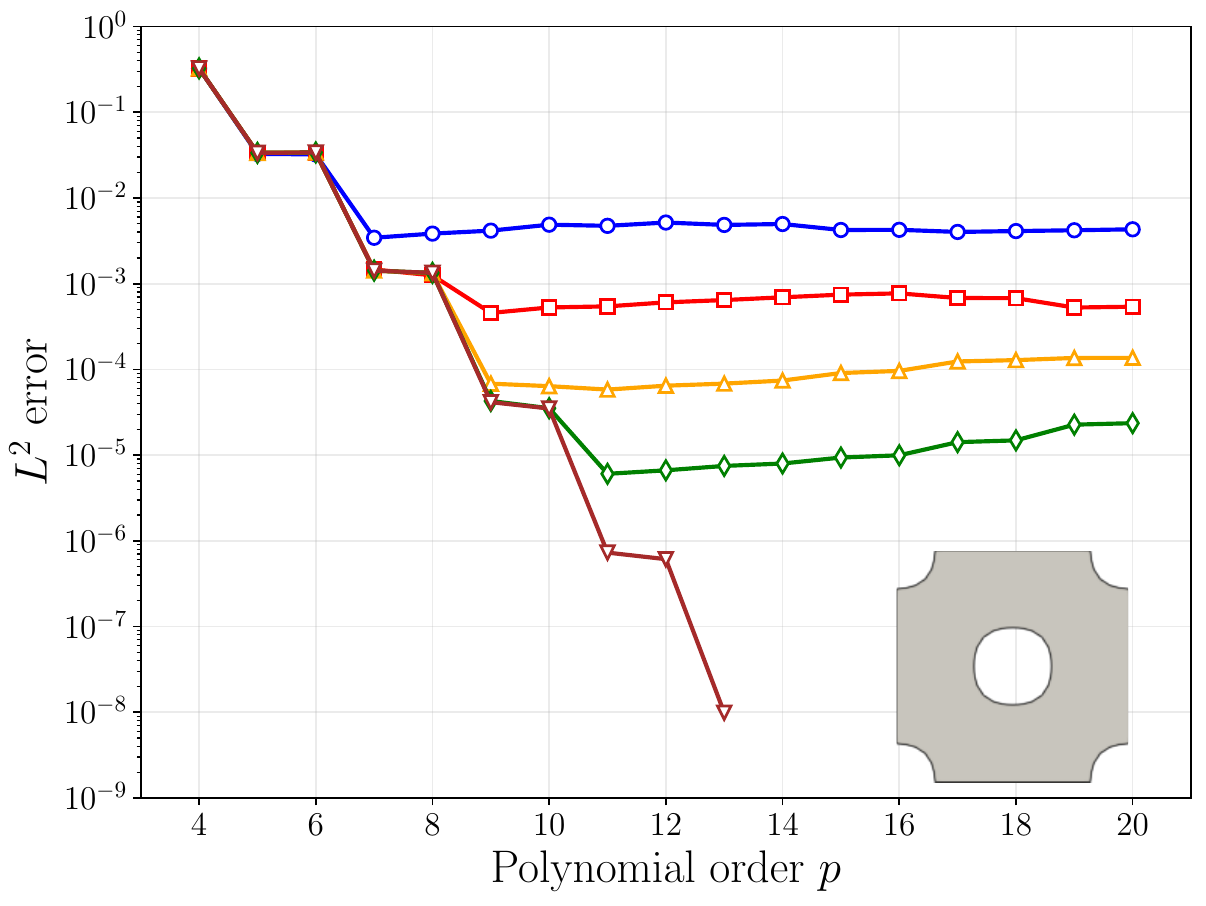}
\caption{}
\end{subfigure}
\begin{subfigure}[b]{0.44\textwidth}\centering
\includegraphics[width=\textwidth]{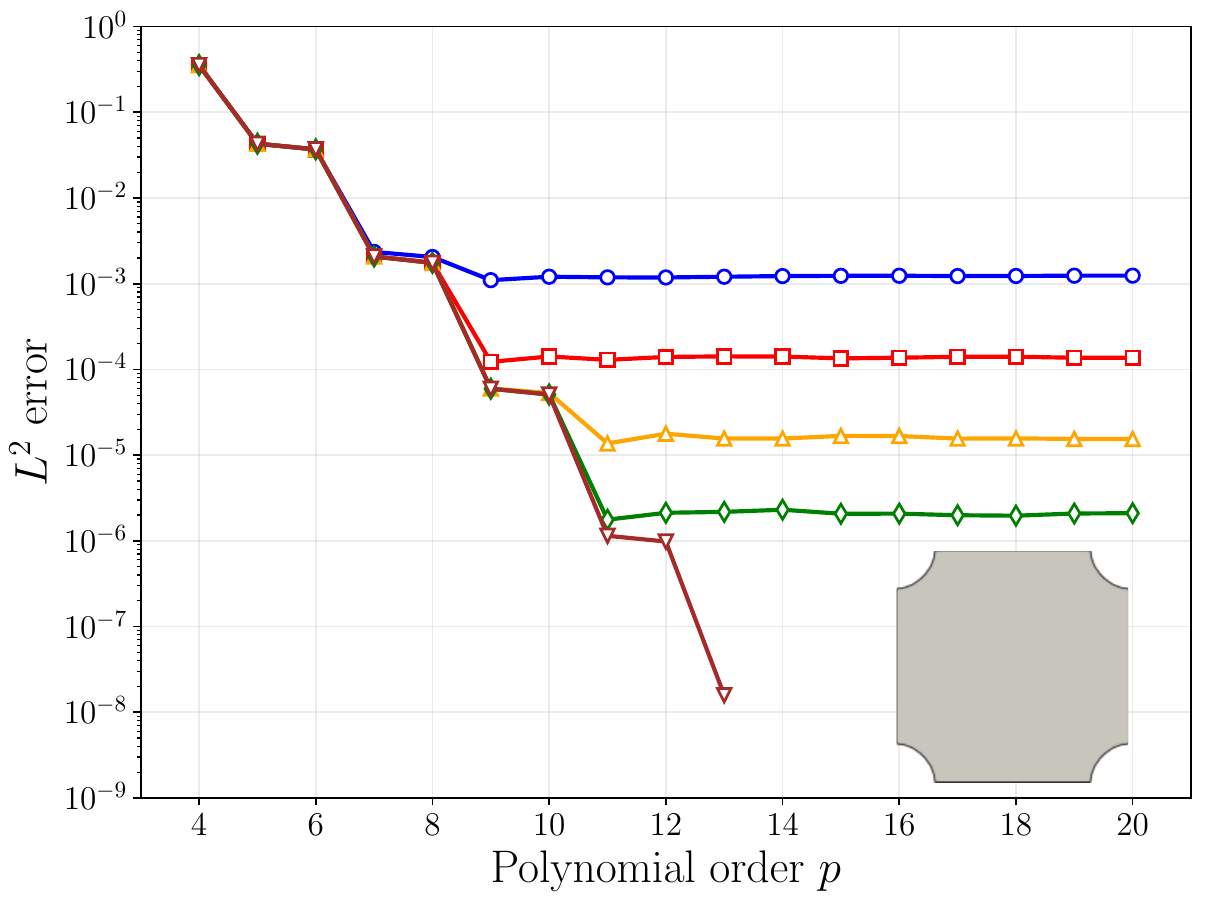}
\caption{}
\end{subfigure}\\
\begin{subfigure}[b]{0.44\textwidth}\centering
\includegraphics[width=\textwidth]{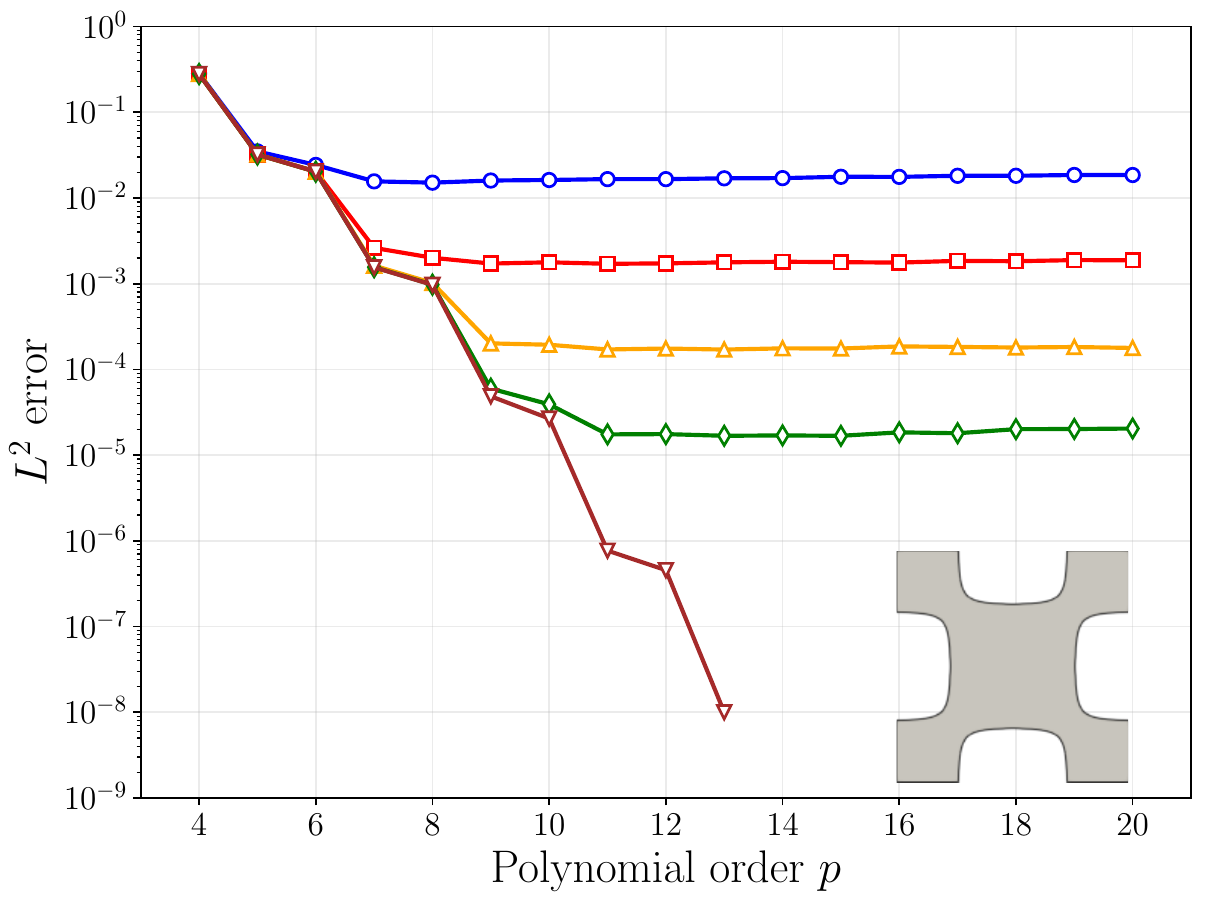}
\caption{}
\end{subfigure}
\begin{subfigure}[b]{0.44\textwidth}\centering
\includegraphics[width=\textwidth]{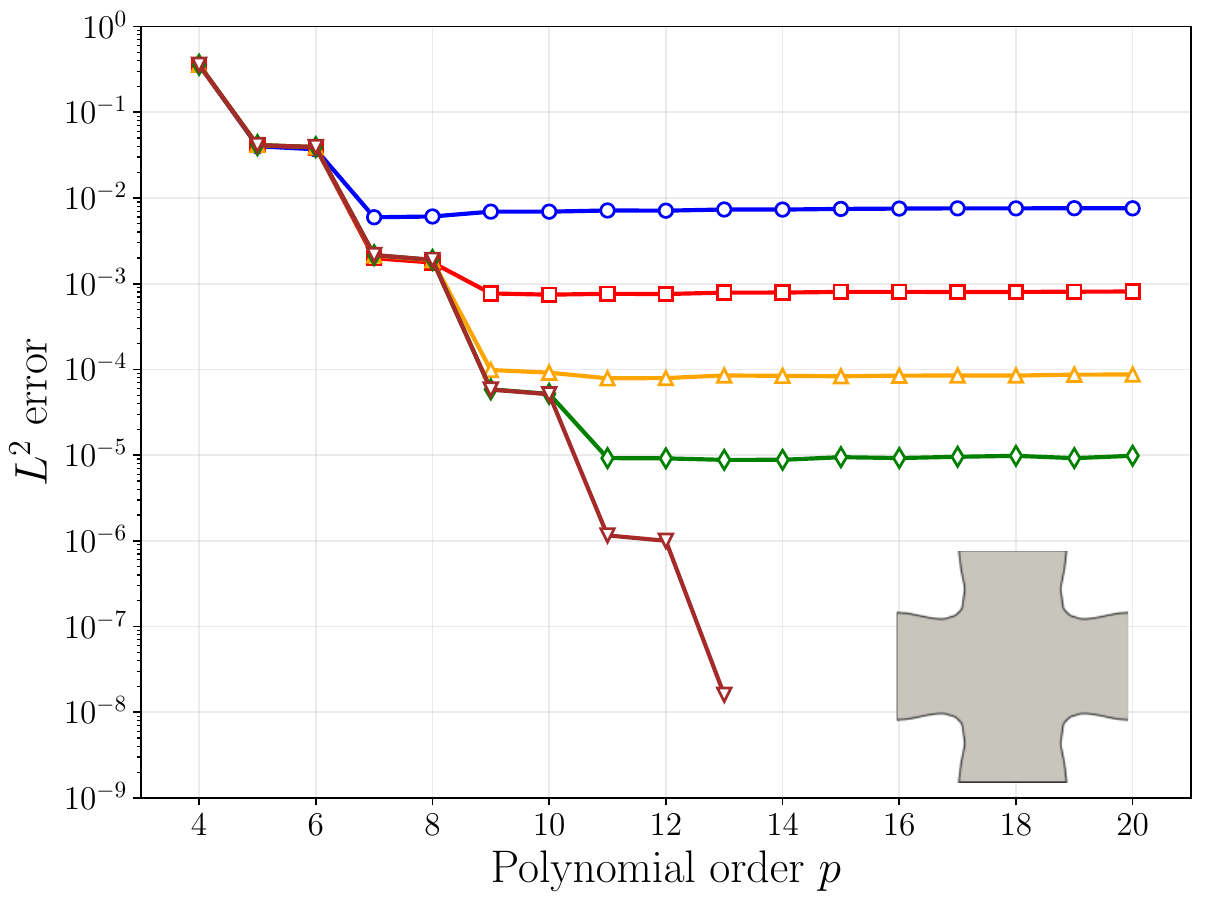}
\caption{}
\end{subfigure}
\caption{\added[id=R2]{$L^2$ error as a function of the approximation degree $p$ for different values of the stabilization parameter $\rho$, computed on a single cell with the manufactured solution described in Section \ref{sssec: RES single cell}. Results are shown for the four TPMS geometries introduced in Figure \ref{fig: GEO TPMS geometries}: (a) Schwarz Diamond, (b) Schwarz Primitive, (c) Schoen FRD, and (d) Schoen IWP.}}
\label{fig: RES stabilization pconv}
\end{figure}

\begin{figure}[ht]
\centering
\begin{subfigure}[b]{0.75\textwidth}\centering
\includegraphics[width=\textwidth]{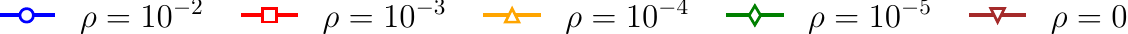}
\end{subfigure}\\
\begin{subfigure}[b]{0.44\textwidth}\centering
\includegraphics[width=\textwidth]{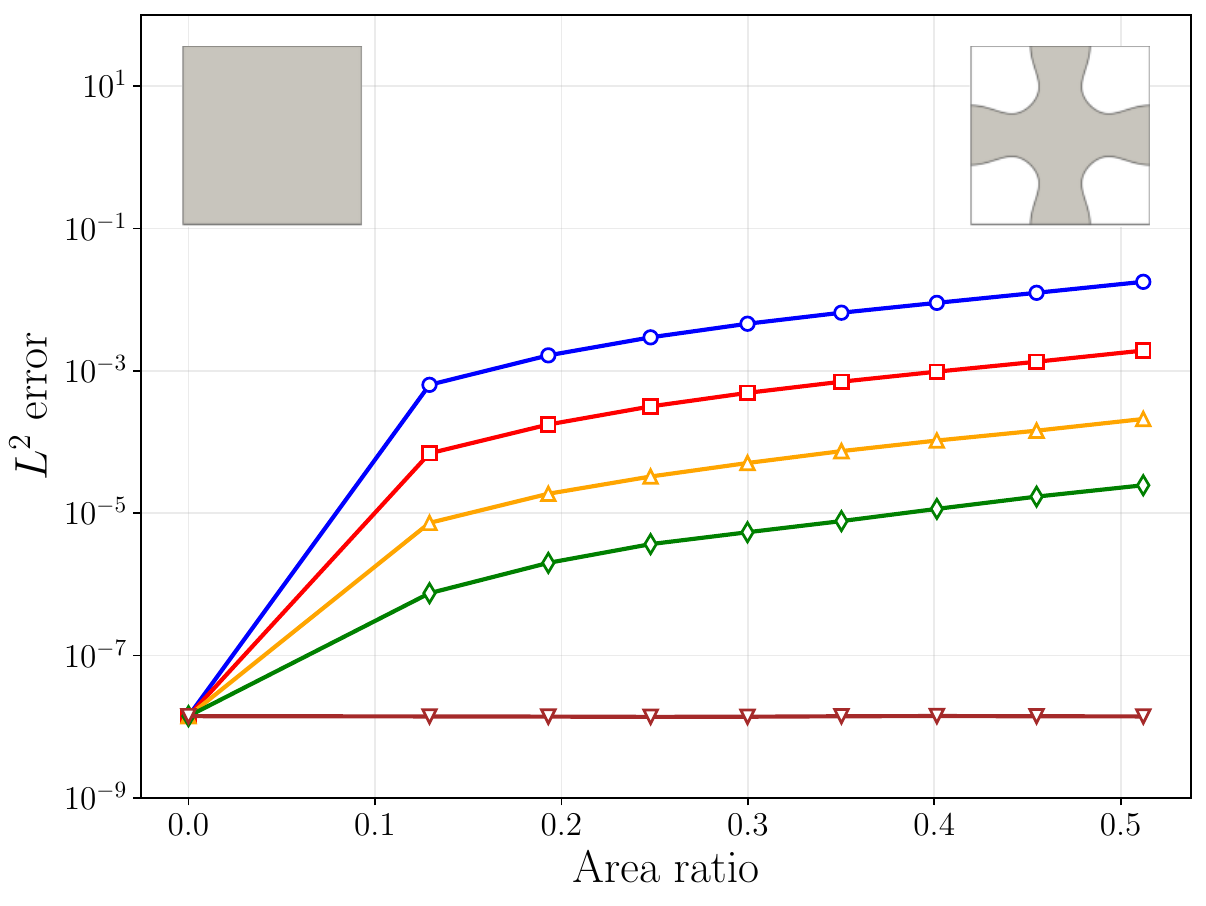}
\caption{}
\label{fig: RES stabilization trimming mapstudy a}
\end{subfigure}
\begin{subfigure}[b]{0.44\textwidth}\centering
\includegraphics[width=\textwidth]{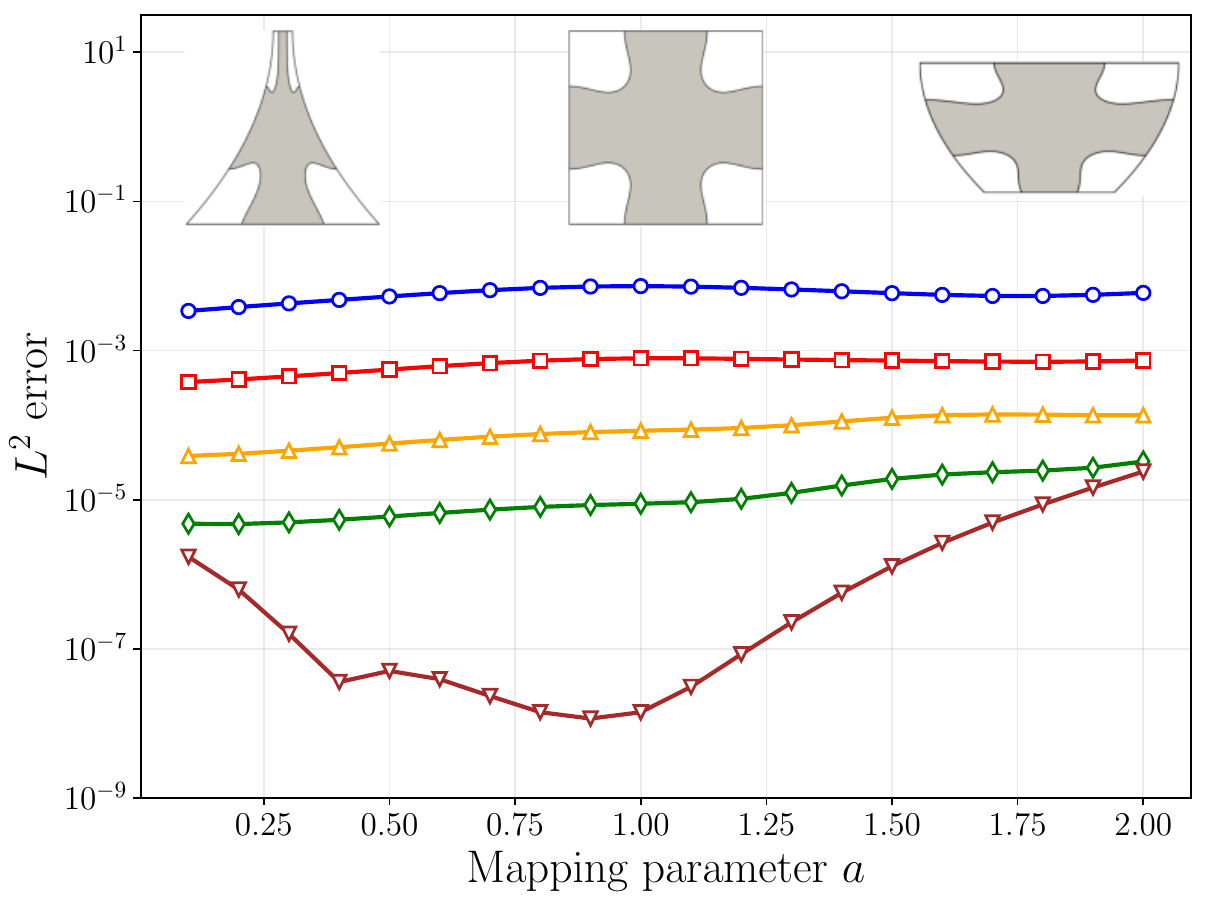}
\caption{}
\label{fig: RES stabilization trimming mapstudy b}
\end{subfigure}
\caption{\added[id=R2]{$L^2$ error for the Schoen IWP geometry at fixed polynomial degree $p=14$, for different values of the stabilization parameter $\rho$. (a) Error as a function of the trimmed area ratio, controlled by a uniform additive shift of the threshold parameter, ranging from no trimming to increasingly trimmed configurations. (b) Error as a function of the mapping distortion parameter $a$, where $a=1$ corresponds to the identity map.}}
\label{fig: RES stabilization trimming mapstudy}
\end{figure}

\added[id=R2]{The accuracy impact of the stabilization parameter $\rho$ is assessed using the single-cell manufactured solution described above, considering five values $\rho \in \{10^{-2}, 10^{-3}, 10^{-4}, 10^{-5}, 0\}$. Three aspects are investigated: the effect on $p$-convergence, the sensitivity to trimming severity, and the robustness with respect to mapping distortion.}

\added[id=R2]{Figure~\ref{fig: RES stabilization pconv} reports the $L^2$ error as a function of the polynomial degree $p$ for all four TPMS geometries, with the identity map ($a=1$) and a uniform threshold parameter (0.1 for the Schwarz Diamond and 0 for the remaining cases). Without stabilization ($\rho = 0$), spectral convergence is achieved in all cases. The staircase pattern for $\rho=0$ results from the higher ability of odd-order polynomials to capture the reference solution during $p$-convergence. As $\rho$ increases, the error eventually plateaus at a level approximately equal to the stabilization parameter itself, with minor differences among geometries. This confirms that the stabilization introduces a consistency error whose magnitude scales linearly with $\rho$, limiting high-order accuracy for large stabilization values.}

\added[id=R2]{Figure~\figref{fig: RES stabilization trimming mapstudy}{a} shows the $L^2$ error at fixed $p = 14$ as a function of the trimmed area ratio. Such a high polynomial degree, compared to the degree $p=8$ employed in the rest of this work, is deliberately chosen to drive the approximation error well below the error floor induced by the stabilization, thereby isolating the consistency error introduced by $\rho$ from the approximation error. Without stabilization, the error remains at the level dictated by the element’s approximation accuracy, regardless of how severely the cell is trimmed, demonstrating that the base method is inherently robust to trimming. With stabilization, the error increases monotonically with the trimmed area, since the penalty term acts over an increasingly large fictitious extension of the domain. Larger values of $\rho$ lead to proportionally higher errors across all trimming levels.}

\added[id=R2]{Figure~\figref{fig: RES stabilization trimming mapstudy}{b} reports the $L^2$ error at fixed $p = 14$ as a function of the mapping distortion parameter $a$, for a homogeneously zero threshold parameter. For $\rho = 0$, the error is minimized near $a = 1$ (identity map), corresponding to a trimmed area ratio of 0.67, and grows for strongly compressed ($a=0.1$) or stretched ($a=2$) configurations, reflecting the reduced approximation accuracy of the polynomial basis defined in the parametric domain. With stabilization, the error remains bounded across the full range of mapping parameters, with the floor determined by $\rho$.}

\added[id=R2]{These results characterize how the interaction between $p$, $\rho$, $\mu$, and $\bm{\mathcal{F}}$ affects the accuracy of the analysis. Extending this characterization to a broader range of multi-cell geometries and loading conditions is a natural continuation of this study.}

\subsubsection{Accuracy assessment \added[id=R2]{on a multi-cells structure}} \label{sssec: RES multiple}

\begin{figure}[ht]
\centering
\begin{subfigure}[b]{0.30\textwidth}
    \centering
    \includegraphics[width=\textwidth]{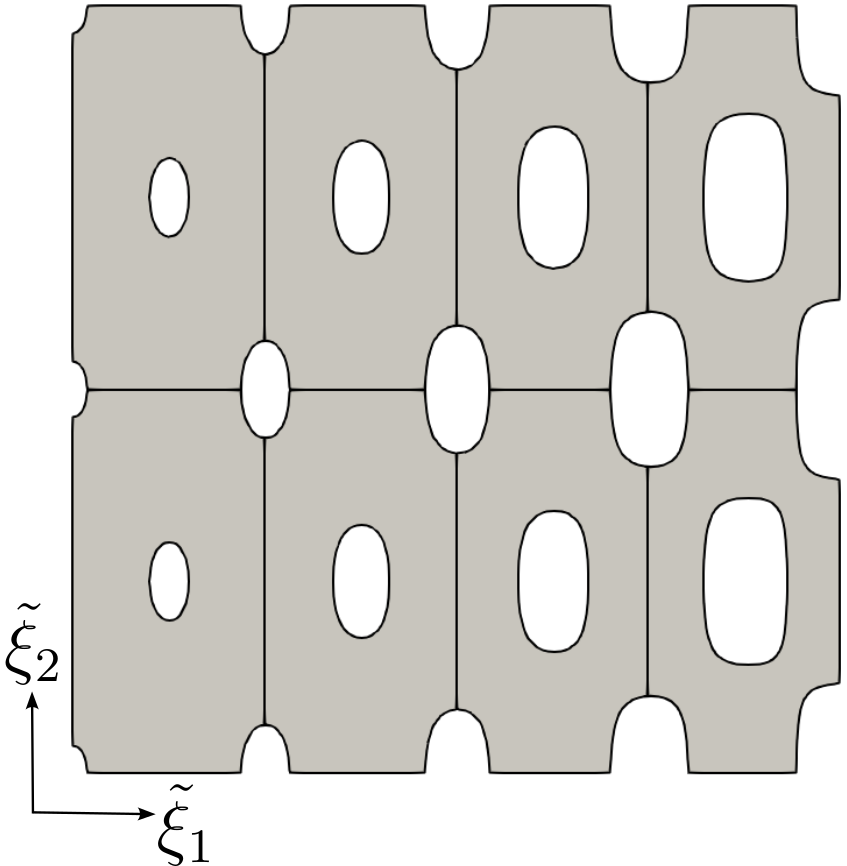}
    \caption{}
\end{subfigure}\hfill
\begin{subfigure}[b]{0.32\textwidth}
    \centering
    \includegraphics[width=\textwidth]{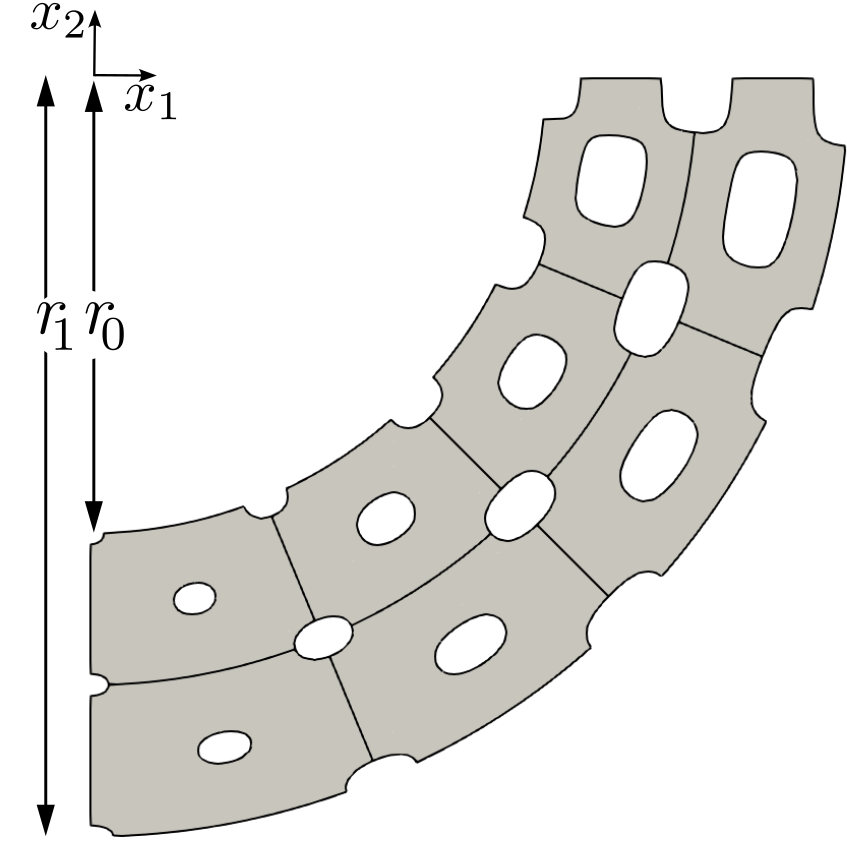}
    \caption{}
\end{subfigure}\hfill
\begin{subfigure}[b]{0.32\textwidth}
    \centering
    \includegraphics[width=\textwidth]{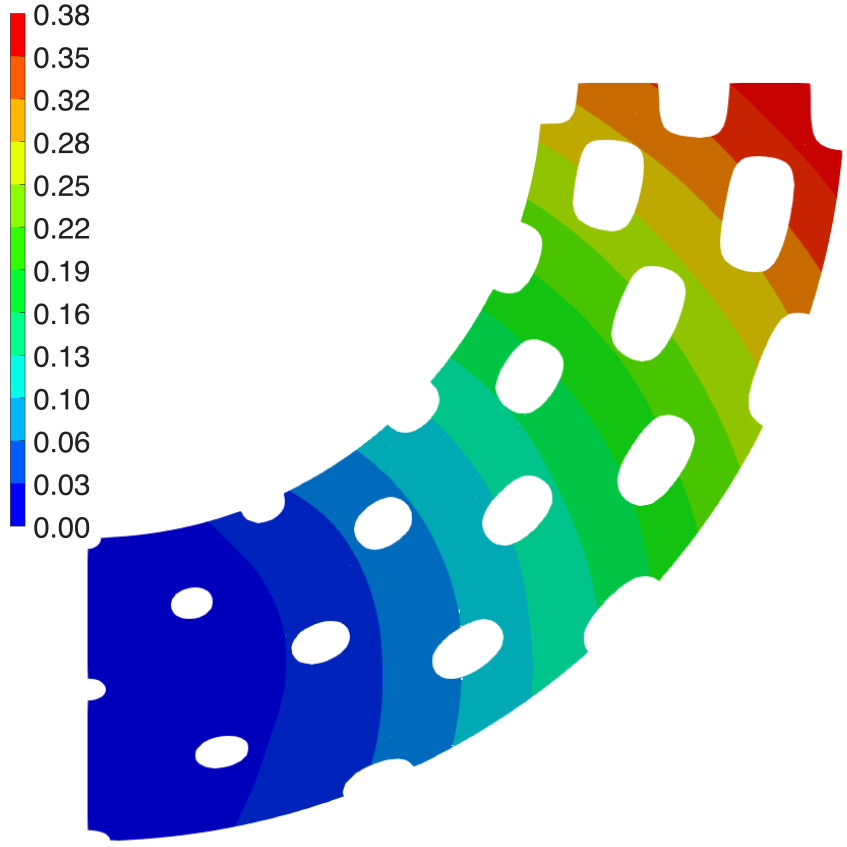}
    \caption{}
\end{subfigure}
\caption{Benchmark geometry for the tests in Sections \ref{sssec: RES multiple} and \ref{sssec: RES performance}. (a) Auxiliary parametric domain after trimming operation. (b) Physical domain resulting from the map given in Equation \eqref{eq: RES multi map}. (c) Contour of the magnitude of the displacements expressed in m.}
\label{fig:test_description}
\end{figure}
The remaining examples in this subsection share a common benchmark geometry, namely the lattice structure depicted in Figure \figref{fig:test_description}. The geometry under consideration is constructed \added[id=R2]{by} adopting an intermediate step with respect to the procedure described in Section \ref{ssec:geometry}. Here, it is introduced an auxiliary parametric domain $\tilde{\Pi} = [0,1]^2$, with associated curvilinear coordinates $\tilde{\xi}_1$ and $\tilde{\xi}_2$. The untrimmed physical domain $\Pi$ is then defined through the auxiliary mapping
\begin{equation}\label{eq: RES multi map}
    \tilde{\mathcal{F}}(\tilde{\xi}_1,\tilde{\xi}_2) =
    \begin{bmatrix}
    +( r_0 + (r_1 - r_0)\tilde{\xi}_2) \sin(\frac{\pi}{2}\tilde{\xi}_1) \\
    -( r_0 + (r_1 - r_0)\tilde{\xi}_2) \cos(\frac{\pi}{2}\tilde{\xi}_1)
    \end{bmatrix},
\end{equation}
where $r_0 = 0.6$ m and $r_1 = 1$ m, corresponding to a quarter annulus.

The untrimmed cells are then constructed with the aid of $\tilde{\Pi}$ by partitioning it into a tensor-product rectangular mesh. The number of cells is chosen such that the number of cells in the $\tilde{\xi}_1$ direction is twice that in the $\tilde{\xi}_2$ direction. This mesh naturally induces a partitioning of $\Pi$ into curvilinear quadrilateral cells. Starting from this geometric construction, it is straightforward to define, for each cell, a mapping from its local parametric domain $\hat{\Pi}^{(i_c)} = [0,1]^2$ to the corresponding physical untrimmed domain in the assembled structure, thereby recovering the procedure described in Section \ref{sssec:single cell map}. Such an intermediate auxiliary parameter domain $\tilde{\Pi}$, although not strictly necessary, is particularly useful when the arrangement of the lattice cells follows a tensor-product structure, allowing the mapping to be reduced to a common global one.

With regard to the trimming procedure described in Section \ref{sssec:single cell trimming}, the level-set function selected here is the Schwarz diamond, as defined in Table \ref{tab: GEO tpms}. The threshold parameters are chosen according to the position of the corresponding node in the auxiliary parametric domain $\tilde{\Pi}$. In particular, the threshold parameters are defined as
\begin{equation*}
\mu_i(\tilde{\bm{\xi}}) = 0.9 - 0.8\,\tilde{\xi}_1\,.
\end{equation*}
It is important to remark that the introduction of the auxiliary parametric domain $\tilde{\Pi}$ is particularly convenient for lattice structures whose cell distribution follows a tensor-product pattern, as in the case of Figure \ref{fig:test_description}. However, this assumption is not strictly required, and structures that do not follow such \added[id=R2]{a} cells organization can also be analyzed within the proposed framework, as will be demonstrated in Section \ref{ssec: RES numerical test}.

The material properties are defined by a Young's modulus $E = 5$ N/m and a Poisson's ratio $\nu = 0.25$. No body forces are applied within the domain. The boundary corresponding to $\tilde{\xi}_1 = 0$ is clamped, while a traction $f = 0.1$ N/m, directed downward in the physical space, is applied on the boundary corresponding to $\tilde{\xi}_1 = 1$. \added[id=R2]{The} deformation magnitude is depicted in Figure \figref{fig:test_description}c.

Throughout the remainder of this subsection, the parameters under investigation are the stabilization parameter $\rho$, the degree $q$ of the fast assembly technique, and the number of basis vectors $n_r$ and clusters $n_k$ for the ROM, as described in Sections \ref{sec: MOD unfitted}, \ref{ssec: Fast Assembly}, and \ref{ssec: ROM}, respectively. A polynomial basis of degree $8$ is employed in all tests. However, when not explicitly varied, or differently specified, the parameters are set to $\rho = 5 \cdot 10^{-4}$, $q = 2$.\deleted[id=R2]{, $n_r = 40$, and $n_k = 2$.}

\paragraph{Accuracy of the fast assembly technique} 
\begin{figure}[ht]
\centering
\begin{subfigure}[b]{0.70\textwidth}\centering\includegraphics[width=\textwidth]{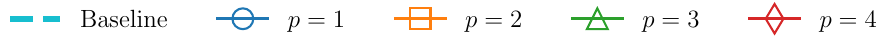}\end{subfigure}\\
\begin{subfigure}[b]{0.44\textwidth}\centering\includegraphics[width=\textwidth]{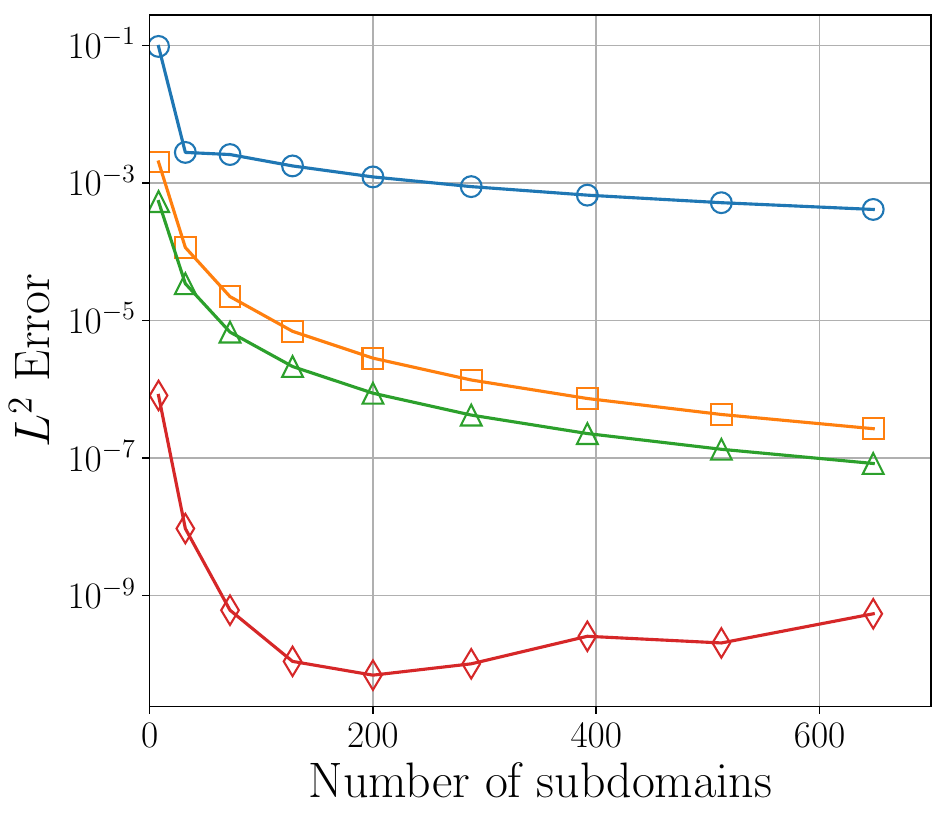}\caption{}\end{subfigure}\hfill
\begin{subfigure}[b]{0.44\textwidth}\centering\includegraphics[width=\textwidth]{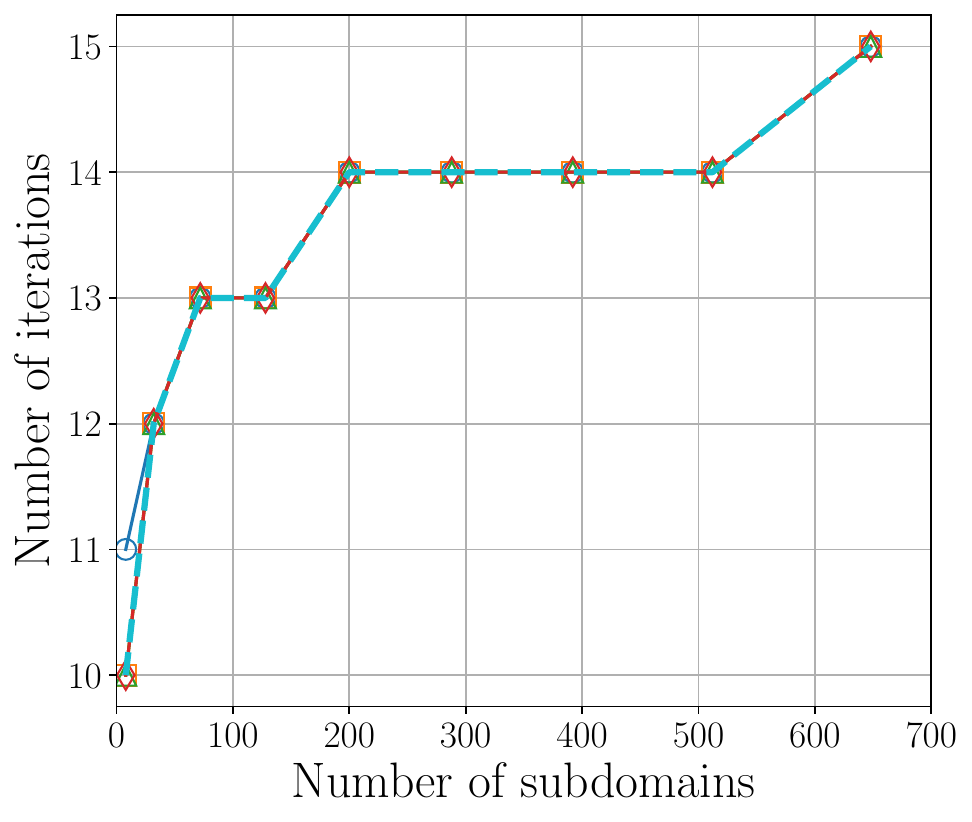}\caption{}\end{subfigure}
\caption{(a) Accuracy of the fast assembly interpolation measured in the $L^2$ norm of the displacements error. (b) Number of iteration\added[id=R2]{s} required by the solver when the fast assembly technique is adopted.}\label{fig:test_6}
\end{figure}

We begin by assessing the capability of the fast assembly technique described in Section \ref{ssec: Fast Assembly} to accurately interpolate the terms that combine the constitutive parameters and external forces with the differential geometry contributions associated with the mapping. Here, the stabilization parameter $\rho$ is set to zero and no ROM acceleration is applied.

To evaluate the accuracy of the interpolation in Equation \eqref{eq: FAST interpolation}, Figure \figref{fig:test_6}{a} shows the $L^2$ error on the displacement vector obtained with the fast assembly technique compared with a standard integration. It can be observed that increasing the interpolation degree reduces the error as expected. Furthermore, the interpolation becomes increasingly accurate as the number of subdomains increases, since the mapping for a small cell approaches a simpler bilinear one.

Additionally, Figure \figref{fig:test_6}{b} reports the number of iteration required by the solver, varying the degree $q$ of the fast assembly interpolation and the number of subdomains. The plots demonstrate that the interpolation degree has no significant influence on the number of iterations. In the remainder of the tests, the interpolation degree for the fast assembly technique is set to two, as it provides a good compromise between accuracy and memory consumption. In fact, a higher interpolation degree leads to larger tensors $I$ and, consequently, a larger ROM. It should be noticed that for $q=2$, even considering a single subdomain, the error in Figure \figref{fig:test_6}{a} is still of the order of magnitude of $10^{-3}$.

\paragraph{Accuracy of the reduced order model} 
\begin{figure}[ht]
    \centering
    \begin{subfigure}
    {0.80\textwidth}\centering\includegraphics[width=\textwidth]{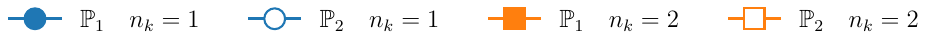}\end{subfigure}\\
    \begin{subfigure}[b]{0.45\textwidth}
        \centering
        \includegraphics[width=\textwidth]{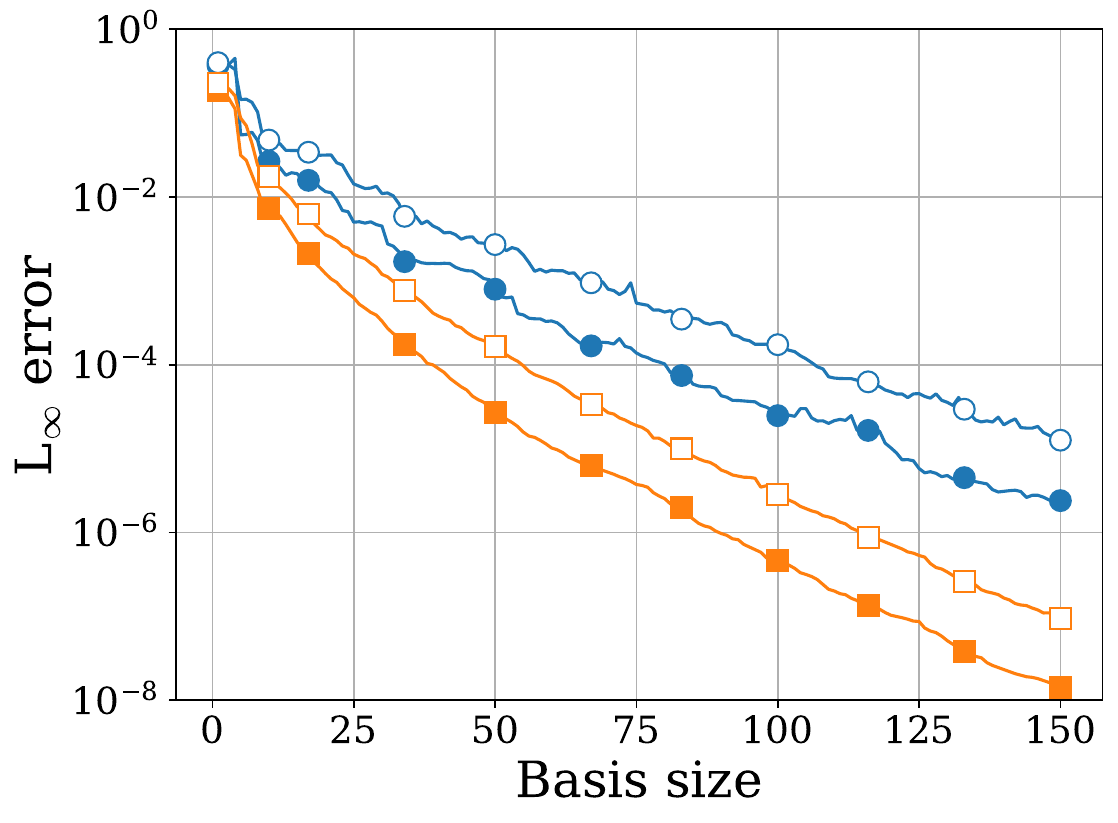}
        \caption{}
    \end{subfigure}
    \hfill
    \begin{subfigure}[b]{0.45\textwidth}
        \centering
        \includegraphics[width=\textwidth]{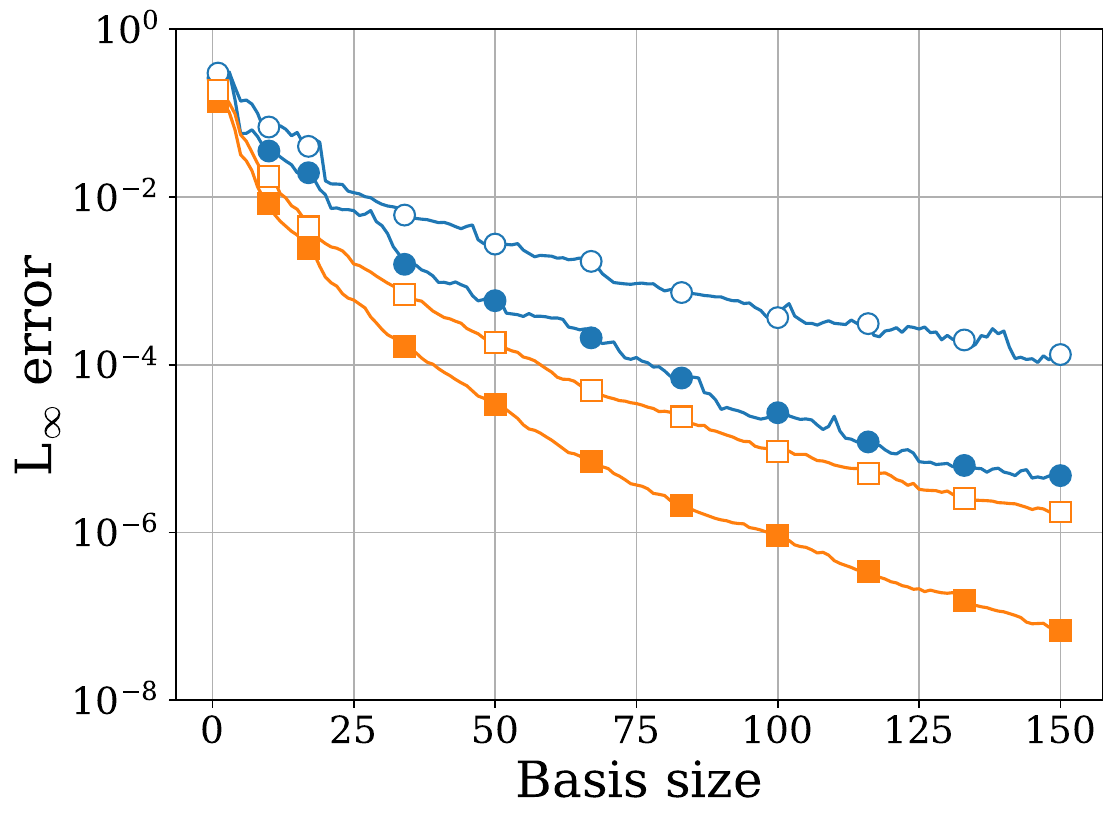}
        \caption{}
    \end{subfigure}
    \\
    \begin{subfigure}[b]{0.45\textwidth}
        \centering
        \includegraphics[width=\textwidth]{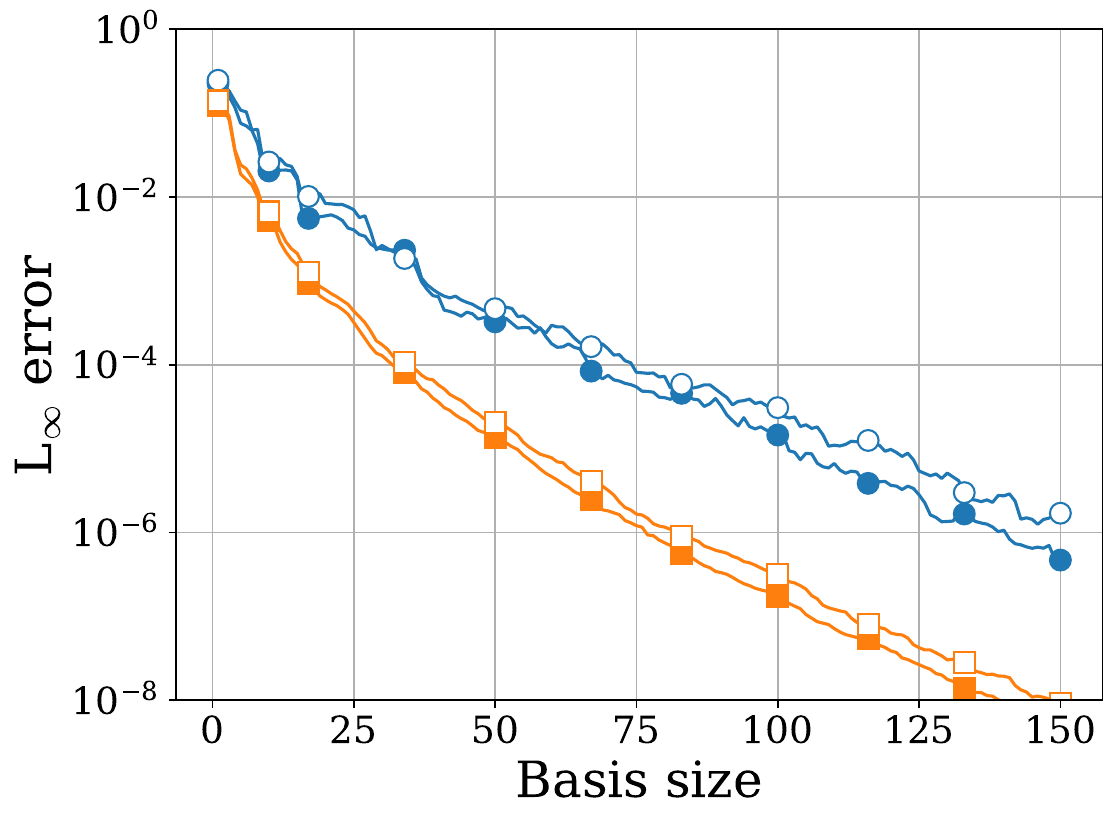}
        \caption{}
    \end{subfigure}
    \hfill
    \begin{subfigure}[b]{0.45\textwidth}
        \centering
        \includegraphics[width=\textwidth]{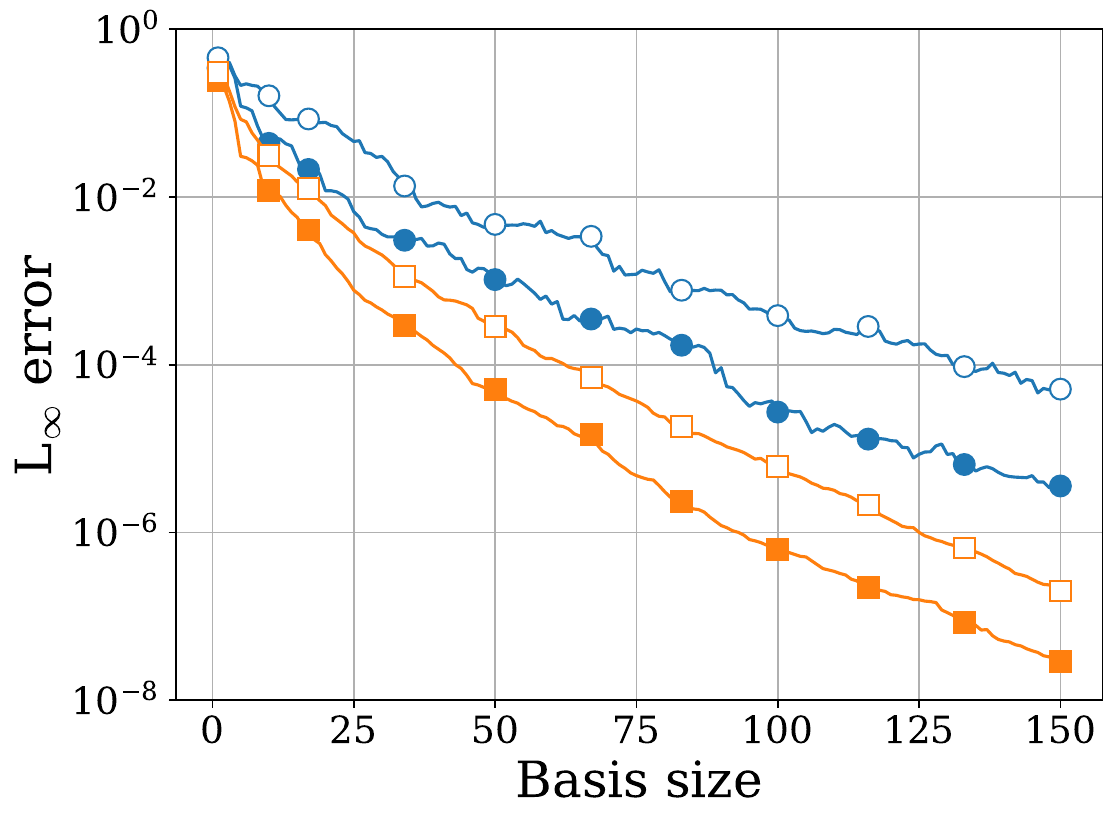}
        \caption{}
    \end{subfigure}
    \caption{Error of the fast assembly tensor depending on the basis size $n_r$ and for different number of cluster\added[id=R2]{s} $n_k$. (a) Schwarz diamond, with space of the threshold parameters $\mathbb{P}_1^\mathrm{SD}$\deleted[id=R1]{$ = \left[0.1, 0.9\right]^{4}$} and $\mathbb{P}_2^\mathrm{SD}$\deleted[id=R1]{$ = \left[0.1, 1.0\right]^{4}$}, (b) Schoen IWP, with \replaced[id=R2]{threshold-parameter space}{space of the threshold parameters} $\mathbb{P}_1^\mathrm{IWP}$\deleted[id=R1]{$ = \left[-2.5, 2.5\right]^{4}$} and $\mathbb{P}_2^\mathrm{IWP}$\deleted[id=R1]{$ = \left[-2.5, 3.0\right]^{4}$}\added[id=R2]{, (c) Schoen FRD, with space of the threshold parameters $\mathbb{P}_1^\mathrm{FRD}$ and $\mathbb{P}_2^\mathrm{FRD}$, (d) Schwarz Primitive, with space of the threshold parameters $\mathbb{P}_1^\mathrm{SP}$ and $\mathbb{P}_2^\mathrm{SP}$}}
    \label{fig:test_rom}
\end{figure}

Next, we study how the performance of the method depends on the number of basis functions $n_r$ and the number of clusters $n_k$ in the ROM, seeking a suitable compromise between accuracy and efficiency. \replaced[id=R2]{Four level-set geometries are considered, namely the Schwarz diamond, the Schoen IWP, the Schoen FRD, and the Schwarz Primitive}{Two level-set geometries are considered, namely the Schwarz diamond and the Schoen IWP} (see Figure \figref{fig: GEO TPMS geometries}). For each geometry, the parameter domain over which the threshold parameters are allowed to vary is specified for the construction of the ROM. In particular, two options are considered for each case, \added[id=R2]{as defined in Section~\ref{ssec: ROM cost}}. For the Schwarz diamond, the threshold parameters are taken within either $\mathbb{P}_1^\mathrm{SD} = [0.1, 0.9]^4$ or $\mathbb{P}_2^\mathrm{SD} = [0.1, 1.0]^4$; for the Schoen IWP\deleted[id=R2]{, the threshold parameters are taken} within either $\mathbb{P}_1^\mathrm{IWP} = \left[-2.5, 2.5\right]^{4}$ or $\mathbb{P}_2^\mathrm{IWP} = \left[-2.5, 3.0\right]^{4}$\added[id=R2]{; for the Schoen FRD, within either $\mathbb{P}_1^\mathrm{FRD} = \left[-6.5, 0.5\right]^4$ or $\mathbb{P}_2^\mathrm{FRD} = \left[-7.0, 0.5\right]^4$; and for the Schwarz Primitive, within either $\mathbb{P}_1^\mathrm{SP} = \left[-0.5, 0.8\right]^4$ or $\mathbb{P}_2^\mathrm{SP} = \left[-0.5, 1.0\right]^4$}. The use of two different parameter spaces for each geometry is intended to illustrate the performance of the ROM when the parameter domain used for its construction is progressively enlarged. Specifically, the sets $\mathbb{P}_2^\mathrm{SD}$, $\mathbb{P}_2^\mathrm{IWP}$\added[id=R2]{, $\mathbb{P}_2^\mathrm{FRD}$, and $\mathbb{P}_2^\mathrm{SP}$} contain the combination of the threshold parameters that reproduces a full cell. \replaced[id=R2]{We recall}{It is reminded} that the ROM models are built on the fast assembly tensor $I$ as described in Sections \ref{ssec: Fast Assembly} and \ref{ssec: ROM}.
The ROM coefficients are computed via Lagra\added[id=R2]{n}gian interpolation (see Section \ref{sssec: MDEIM}), with interpolation order 6 used throughout the results presented below.

Figure \figref{fig:test_rom} reports the results of this study. \replaced[id=R2]{Figure \figref{fig:test_rom}{a} refers to the Schwarz diamond level-set, Figure \figref{fig:test_rom}{b} to the Schoen IWP, Figure \figref{fig:test_rom}{c} to the Schoen FRD, and Figure \figref{fig:test_rom}{d} to the Schwarz Primitive.}{Figure \figref{fig:test_rom}{a} refers to the Schwarz diamond level-set while Figure \figref{fig:test_rom}{b} refers to the Schoen IWP.} The vertical axis represents a measure of the ROM error, computed as follows. For 50 different randomly generated snapshots (for each cluster) that are not included in the singular value decomposition training set, the L$_\infty$ norm of the difference between the vectorized tensor $I$ obtained through full integration and its projection in the reduced-order basis is evaluated and normalized by the L$_\infty$ norm of $I$. The final ROM error is then obtained as the arithmetic average of these values.

As expected, the error decreases as both the number of basis vectors and the number of clusters increase, reflecting the improved accuracy of the ROM. \replaced[id=R2]{Furthermore, as the parameter domain expands from $\mathbb{P}_1$ to $\mathbb{P}_2$, a deterioration in accuracy is observed for all the tested geometries. Similar trends are evident in all panels of Figure \figref{fig:test_rom}; however, the error level is geometry-dependent, with higher values observed for the Schoen IWP and a milder but still visible deterioration also for the newly included Schoen FRD and Schwarz Primitive cases, indicating that the ROM reliability is sensitive to both the geometry family and the extent of the threshold-parameter space.}{Furthermore, as the parameter domain expands from $\mathbb{P}_1^\mathrm{SD}$ to $\mathbb{P}_2^\mathrm{SD}$, a deterioration in accuracy is observed. Similar trends are evident in Figure \figref{fig:test_rom}{b}; however, the error is generally higher for this level-set function, indicating that the ROM accuracy remains dependent on the underlying geometry.}

For the remainder of this paper, unless otherwise specified, the ROM is constructed, \added[id=R2]{as described in Section~\ref{ssec: ROM cost}}, using two clusters per threshold parameter (resulting in a total of 16 clusters) and 40 basis vectors per cluster. Furthermore, the threshold parameters are taken in the sets $\mathbb{P}_1$. This choice ensures an accuracy comparable to that of the fast assembly procedure, so that the different components of the method exhibit a balanced level of approximation.

\added[id=R2]{We remark that, since in the online phase the surrogate is only queried inside the threshold-parameter box sampled during training, the study of Figure~\figref{fig:test_rom} constitutes a practical acceptance test that can be carried out entirely at the ROM training stage, at a marginal cost compared to the training itself (50 additional full-quadrature evaluations per cluster): the user selects the basis size and the number of clusters so that the indicator meets a prescribed tolerance, or otherwise rejects the model and retrains it with more clusters or a smaller parameter range. The reported projection error relates to the actual interpolation error through a computable amplification factor associated with the magic-point rows of the basis \cite{chaturantabut2010nonlinear}. Its transfer to the solution follows a first-order perturbation argument: the entries of the cell operators depend linearly on the tensor $I$, so a relative error $\varepsilon$ on $I$ induces an operator perturbation of the same order, and the ensuing relative solution error is bounded, to first order, by $\varepsilon$ times a constant depending on the conditioning of the stabilized problem. For energy-type quantities of interest, such as the compliance, the first-order error is bounded by the energy-relative operator perturbation, with no conditioning amplification. These transfer constants are moderate in practice: as will be shown in the wrench example of Section~\ref{sssec: RES wrench}, replacing the ROM by full integration changes the monitored quantities of interest by $1.4\cdot10^{-4}$, of the same order as the offline indicator for the settings employed here. A practical guideline is thus to select, at the ROM training stage, the ROM configuration such that the offline indicator lies below the target discretization error, as done in this work.}

\paragraph{Effect of the stabilization}
\begin{figure}[ht]
\centering
\includegraphics[width=0.58\textwidth]{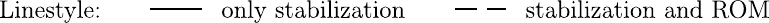}\\
\includegraphics[width=0.98\textwidth]{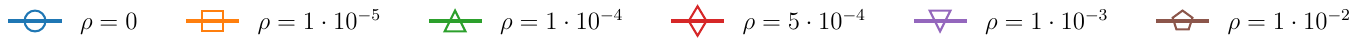}\\
\begin{subfigure}[b]{0.48\textwidth}\centering\includegraphics[width=\textwidth]{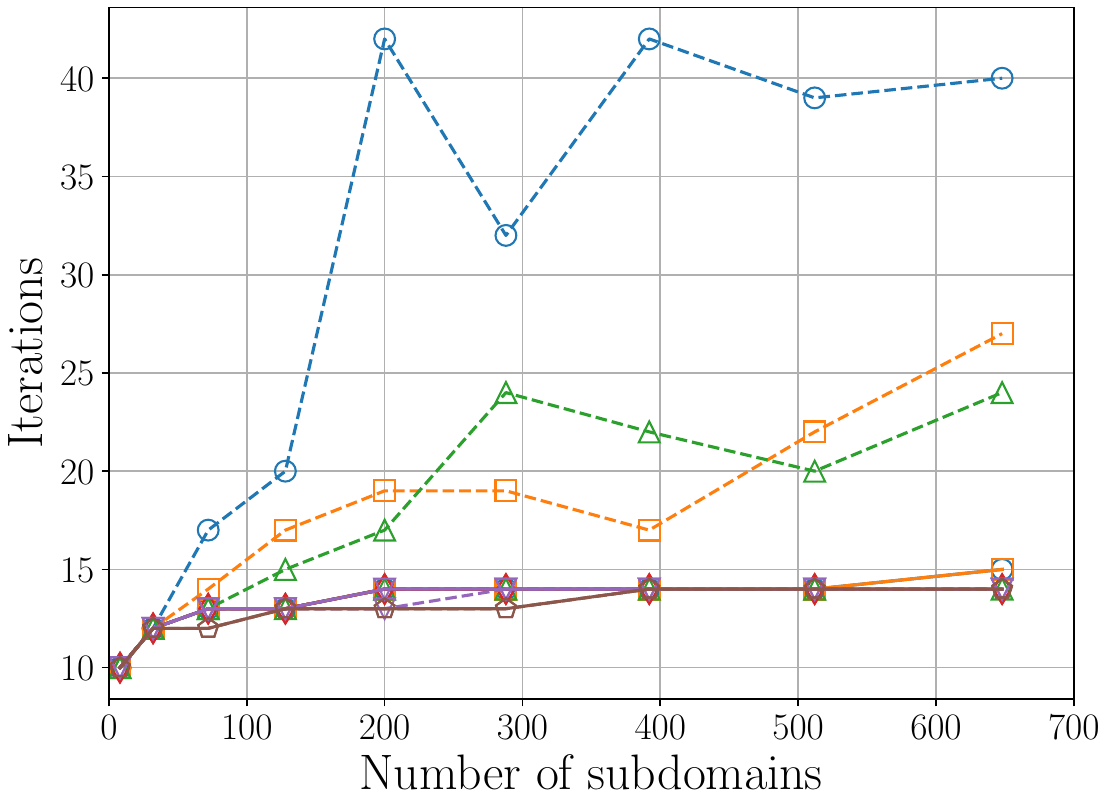}\caption{}\end{subfigure}\hfill
\begin{subfigure}[b]{0.48\textwidth}\centering\includegraphics[width=\textwidth]{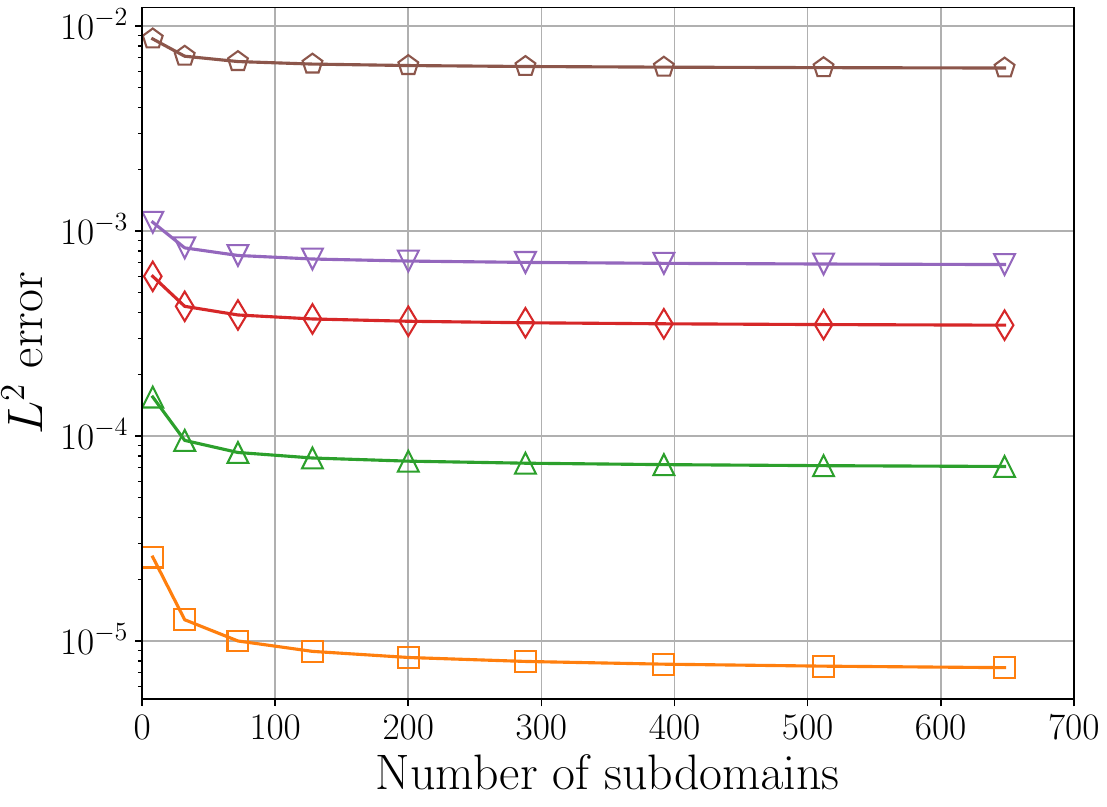}\caption{}\end{subfigure}
\begin{subfigure}[b]{0.48\textwidth}\centering\includegraphics[width=\textwidth]{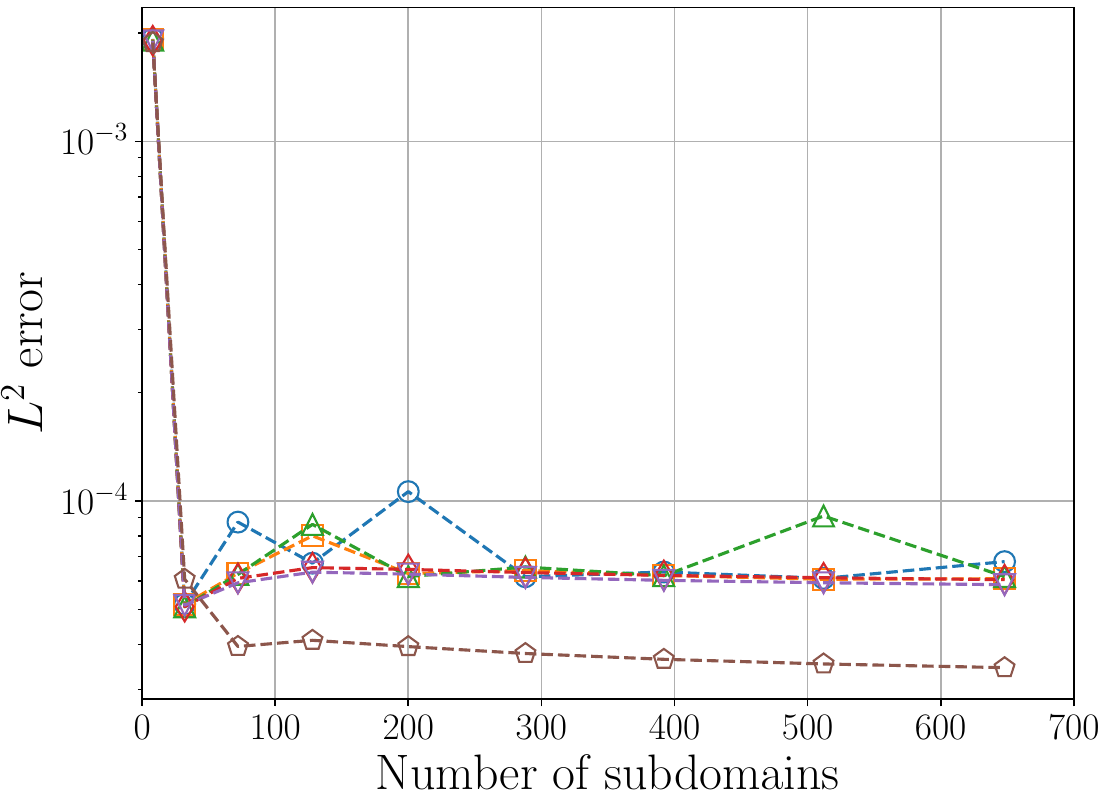}\caption{}\end{subfigure}\hfill
\begin{subfigure}[b]{0.48\textwidth}\centering\includegraphics[width=\textwidth]{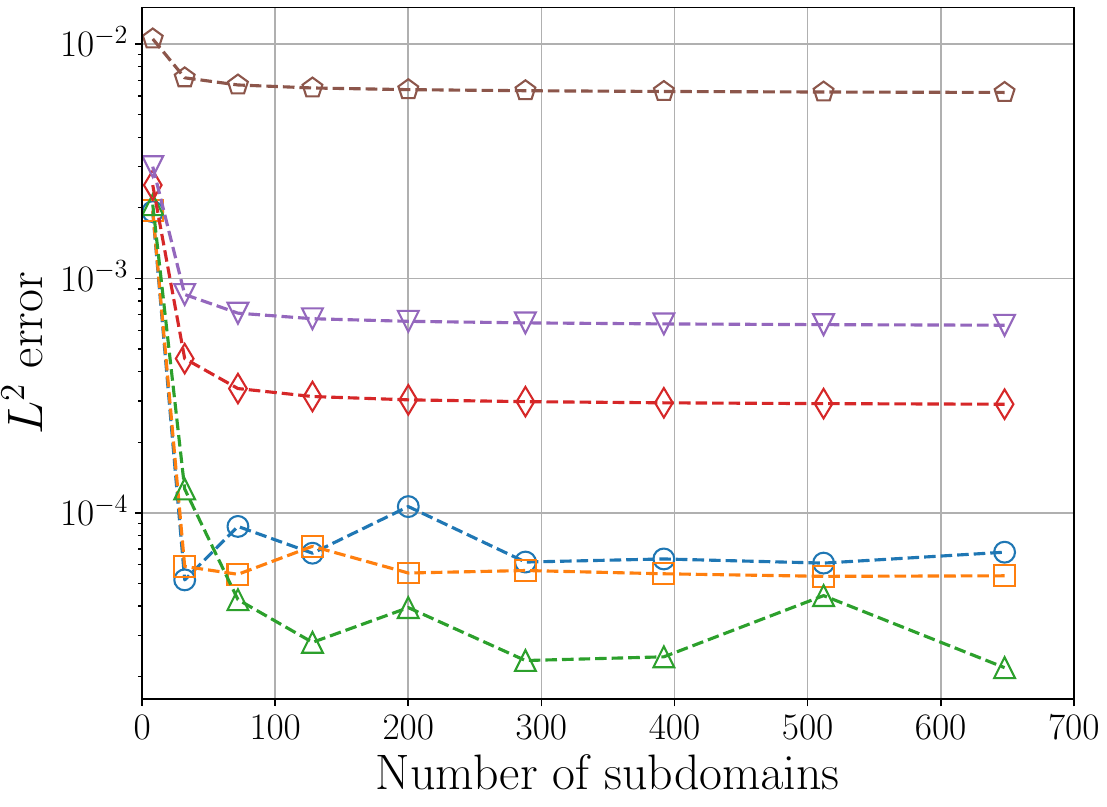}\caption{}\end{subfigure}
    \caption{Effect of the stabilization \replaced[id=R2]{on}{in} the solver iterations and the solution error. (a) Iterations with and without ROM for different stabilizations. (b)  Error between the original solution and stabilized solution without using ROM. (c) Error between the solution using ROM and without ROM (both with the same stabilization). (d) Error between the original solution and the solution using ROM + stabilization}
    \label{fig:test_2}
\end{figure}

The following analysis is aimed at demonstrating the effect of stabilization on the different components of the proposed method. Figure \figref{fig:test_2} summarizes the results and illustrates its interaction with the ROM. It is worth noting that, since the ROM is constructed from the fast assembly tensors, its use is inherently coupled with the fast assembly technique. In contrast, when the ROM is not employed, standard numerical integration is used for matrix assembly.

Figure \figref{fig:test_2}{a} shows the number of iterations required by the BDDC preconditioner when employing both the ROM and stabilization, and when using stabilization alone. It can be observed that stabilization has a negligible effect on the number of iterations when standard integration is adopted (solid lines). However, when the matrices are assembled through the ROM (dashed lines), increasing the stabilization parameter helps control the number of iterations, bringing it closer to the baseline case, while reducing it can cause the iterations to grow unbounded. This highlights that the stability of the solver can be delicate. Although not reported here, for different microgeometries, the solver performed well when stabilization was applied. However, for some geometries, it did not scale properly in the absence of stabilization, even when no ROM was used.

Figure \figref{fig:test_2}{b} assesses the consistency error associated with the stabilization by measuring the $L^2$ norm of the difference between the solutions obtained with and without stabilization. As expected, increasing the stabilization parameter also increases the error. Moreover, for sufficiently large numbers of cells, the error no longer varies as a function of this quantity. \added[id=R2]{In addition, the error levels are approximately proportional to $\rho$, consistently with the single-cell study of Figure~\figref{fig: RES stabilization pconv}: the $O(\rho)$ scaling of the consistency error carries over from the single cell to multi-cell structures.}

Figure \figref{fig:test_2}{c} compares the solution obtained with both ROM and stabilization to the solution obtained with stabilization alone, in order to isolate the contribution of the ROM. It can be observed that, for values below $10^{-3}$, this error becomes independent of the stabilization, indicating that it is primarily due to the ROM approximation, whose accuracy is independent from $\rho$.

Finally, Figure \figref{fig:test_2}{d} depicts the error obtained by comparing the solution without ROM or stabilization to the solution computed with both ROM and stabilization. It can be observed that, initially, the accuracy is primarily governed by the stabilization, and decreasing $\rho$ leads to a reduction in the error. However, once the accuracy limit of the ROM is reached, the error saturates at approximately the same value, and further reduction of $\rho$ does not yield any additional improvement.

The value chosen for the stabilization parameter in the remainder of this section is $\rho = 5 \cdot 10^{-4}$, \added[id=R2]{which as demonstrated by these results and those of the related paragraph of Section \ref{sssec: RES single cell},} is sufficient to control the number of solver iterations while maintaining a satisfactory level of accuracy.
 
\added[id=R2]{More generally, a practical guideline is that $\rho$ should be chosen small enough that the consistency error $O(\rho)$ lies below the desired accuracy at the working polynomial degree $p$, and large enough to ensure bounded BDDC iterations when the ROM is active. Typical accuracy requirements in structural analysis are of order $10^{-2}$ or better, so $\rho = 5 \times 10^{-4}$ is comfortably below this threshold while still controlling the iteration count, as suggested by the results in Figure \figref{fig:test_2}{a}. A sharp analytical criterion, relating $\rho$ explicitly to $p$, the Jacobian of the geometric map, the trimmed area fraction, or the loading, remains an open question.}

\subsubsection{Computational performance} \label{sssec: RES performance}

\added[id=R2]{This section assesses the computational performance of the proposed solver through three complementary studies: a comparison against representative direct and iterative solvers, the benefit of the proposed acceleration layers against an unaccelerated BDDC baseline, and a scalability study up to $17{,}000$ subdomains. The geometry, materials, loads, and discretization are the same as in Section \ref{sssec: RES multiple}.}

\paragraph{Comparison with other solvers}
\begin{figure}[ht]
\centering
\includegraphics[width=0.35\textwidth]{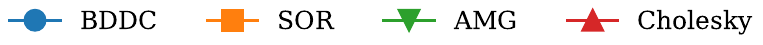}\\
\begin{subfigure}[b]{0.48\textwidth}\centering\includegraphics[width=\textwidth]{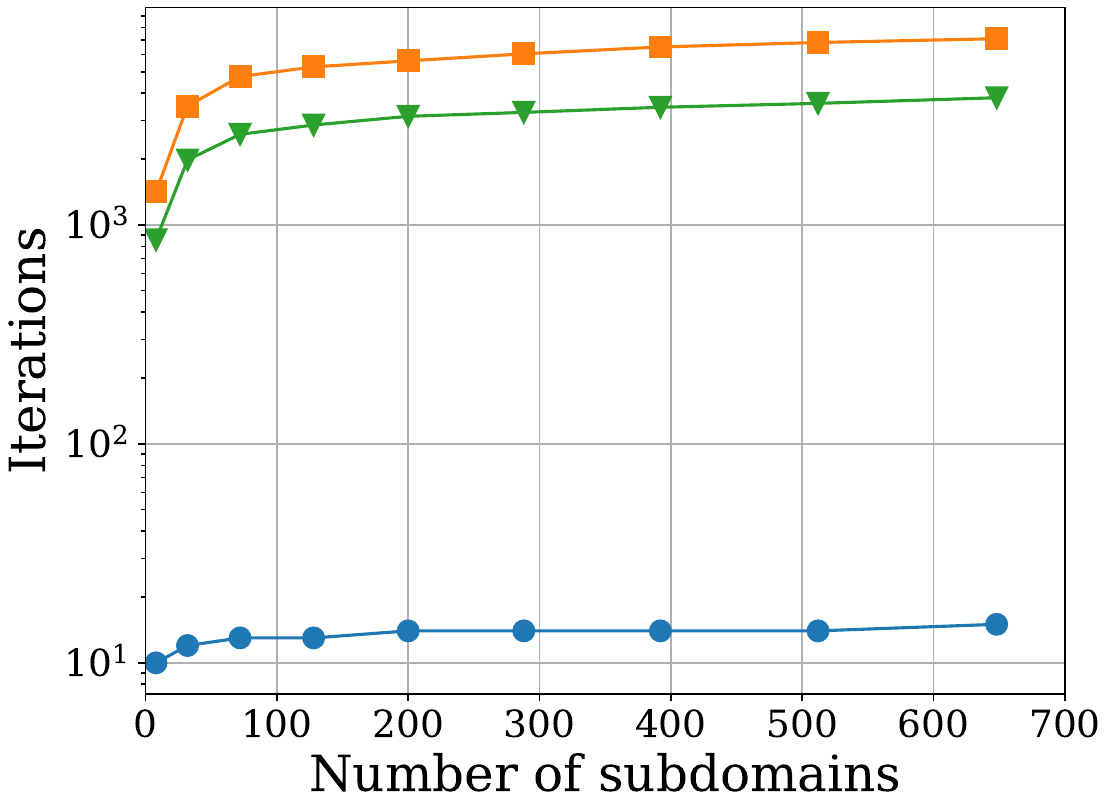}\caption{}\end{subfigure}\hfill
\begin{subfigure}[b]{0.48\textwidth}\centering\includegraphics[width=\textwidth]{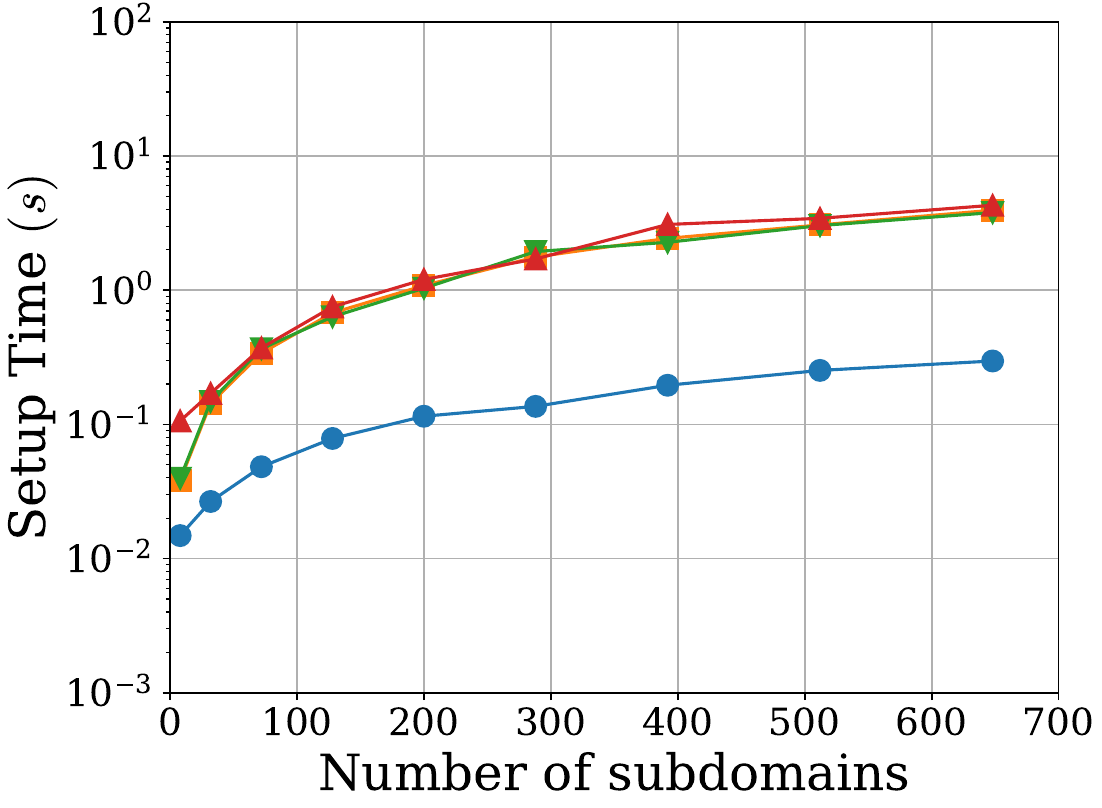}\caption{}\end{subfigure}
\begin{subfigure}[b]{0.48\textwidth}\centering\includegraphics[width=\textwidth]{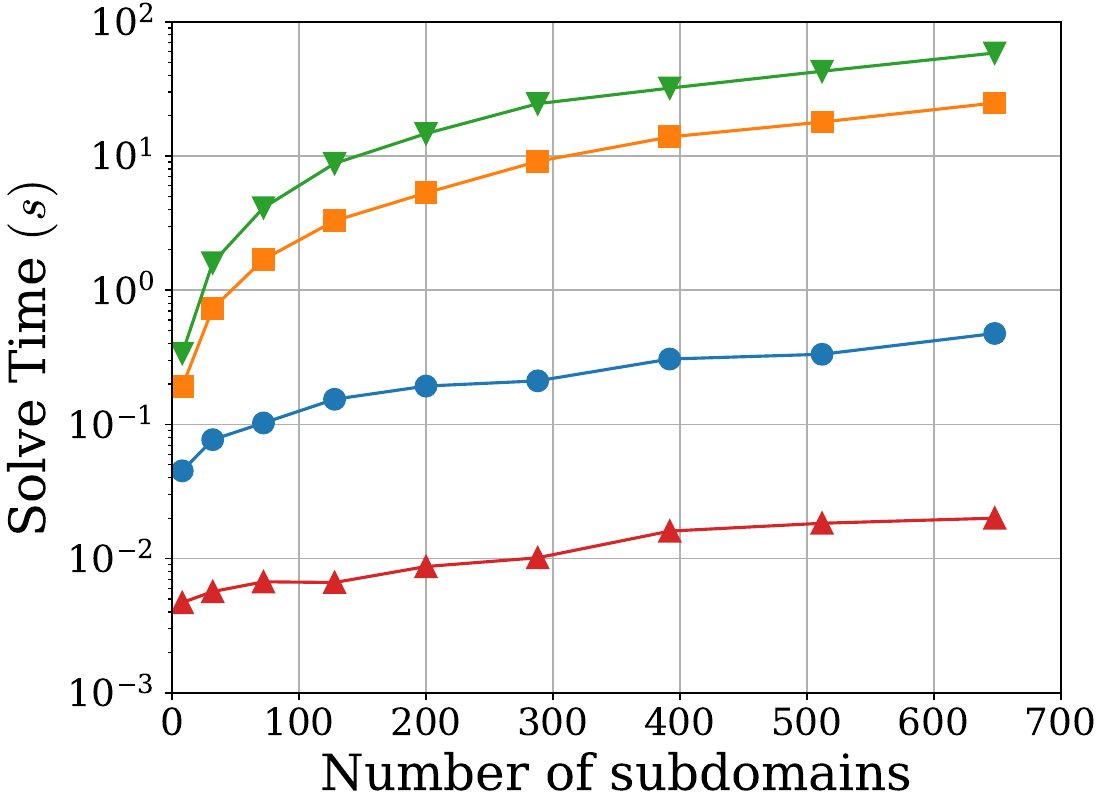}\caption{}\end{subfigure}\hfill
\begin{subfigure}[b]{0.48\textwidth}\centering\includegraphics[width=\textwidth]{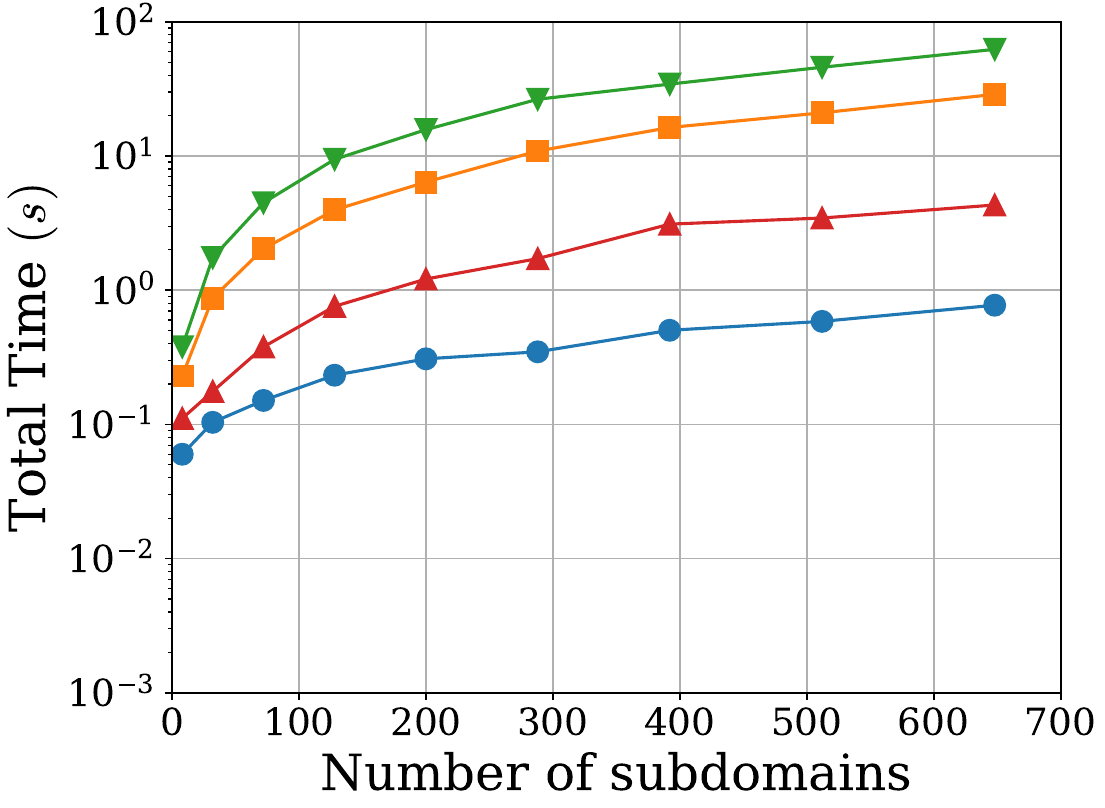}\caption{}\end{subfigure}
\caption{Performances of \replaced[id=R2]{four}{three} different solvers are compared: BDDC, Cholesky decomposition, \replaced[id=R2]{SOR, and AMG}{and SOR}. (a) Number of iterations required by the iterative solvers. (b), (c) Time required for the setup and solve phases, respectively. (d) Total computational time required by the solvers.}\label{fig:test_1}
\end{figure}

We first compare a basic version of the BDDC method with \replaced[id=R2]{three}{two} widely used solvers for linear systems: the Cholesky factorization, representing a direct method, and the conjugate gradient method combined with \added[id=R2]{either} the Successive Over-Relaxation (SOR) preconditioner \cite{young2014iterative} \added[id=R2]{or an algebraic multigrid (AMG) preconditioner \cite{ruge1987algebraic}, in its smoothed-aggregation variant (GAMG) as implemented in PETSc}, representing \replaced[id=R2]{iterative approaches}{an iterative approach}. All the solvers in this comparison rely on their PETSc implementations, accessed through petsc4py \cite{petsc-user-ref,petsc4py}, so that the results are not biased by implementation differences. In this comparison, the BDDC method is used without any of the acceleration techniques introduced in Section \ref{sec: SPEED}. Furthermore, the stabilization parameter $\rho$ is also set to zero. The performance is evaluated for an increasing number of cells. The trends for the \replaced[id=R2]{four}{three} solvers are reported in Figure \ref{fig:test_1}. In particular, Figure \figref{fig:test_1}{a} shows the number of iterations required by the iterative methods, highlighting that BDDC is more efficient than \replaced[id=R2]{both SOR and AMG}{SOR} by more than two orders of magnitude. It should be noted, however, that the iterations in these \deleted[id=R2]{two} approaches are performed on different systems. \replaced[id=R2]{The SOR and AMG methods operate}{The SOR method operates} on the full system given in Equation \eqref{eq: BDDC full system}, whereas BDDC is applied only to the condensed system defined in Equation \eqref{eq: BDDC system}.

Regarding computational time, Figures \figref{fig:test_1}{(b-d)} illustrate the setup time, solve time, and total time for the \replaced[id=R2]{four}{three} solvers. More precisely, the setup time does not include the assembly of the local stiffness matrices, as this cost is common to all methods. Although it represents the main portion of the computational effort, as shown in the remainder of this section, local assembly time can be reduced to a level comparable with setup and solve time by employing the acceleration techniques discussed in Section \ref{sec: SPEED}. Therefore, the setup time includes:
\begin{itemize}
\item BDDC: construction of the local operators (e.g., Schur complements, restriction/extension matrices) and preparation of the coarse problem solver (assembly of the coarse correction operator $S_C$ and its Cholesky factorization);
\item Cholesky: assembly of the global matrix and parallel computation of its Cholesky factorization;
\item \added[id=R2]{AMG: assembly of the global matrix and construction of the multigrid hierarchy (coarse levels, transfer operators, and smoothers);}
\item SOR: assembly of the global matrix and standard preprocessing required by the SOR method.
\end{itemize}
The most striking feature of these results is that BDDC outperforms \replaced[id=R2]{all the alternative}{both alternative} methods in terms of total computational time, being approximately one order of magnitude faster than the Cholesky decomposition and \replaced[id=R2]{about two orders of magnitude faster than SOR and AMG}{two orders of magnitude faster than SOR}. More specifically, for the Cholesky decomposition, the majority of the computational cost is associated with the setup phase, while the solution phase is comparatively fast, as expected for a direct method. In contrast, for SOR \added[id=R2]{and AMG}, most of the computational effort is concentrated in the solution phase due to the large number of iterations required. \added[id=R2]{In this regard, it is worth noting that, although AMG requires fewer iterations than SOR, each of its iterations is more expensive, as it involves the application of a full multigrid cycle, resulting in an overall longer solve time. Moreover, the relatively poor performance of AMG is not unexpected, since standard algebraic multigrid methods are known to lose efficiency for high-order discretizations, and the ill-conditioning introduced by trimming further degrades the effectiveness of the coarsening and smoothing procedures.} For BDDC, the setup and solution times are of comparable magnitude, leading to a more balanced computational cost distribution.

\paragraph{Efficiency of the accelerations}
\begin{figure}[ht]
\centering
\begin{subfigure}[b]{0.48\textwidth}\centering\includegraphics[width=\textwidth]{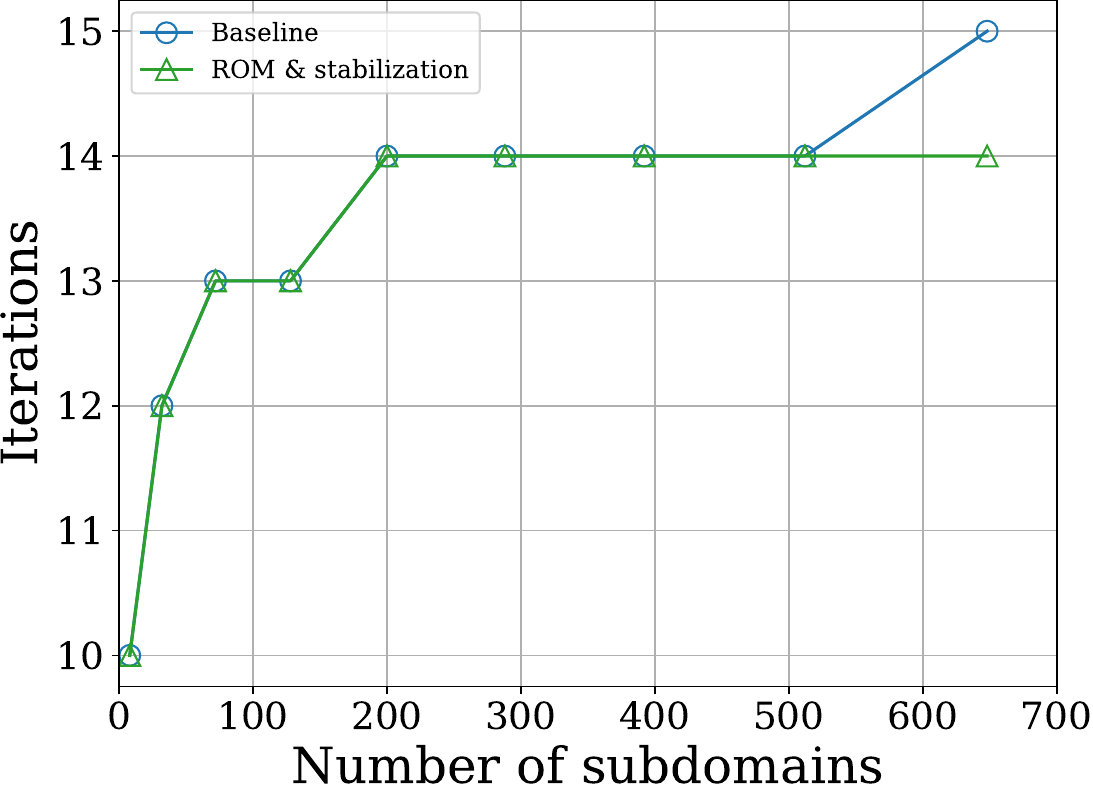}\caption{}\end{subfigure}\hfill
\begin{subfigure}[b]{0.48\textwidth}\centering\includegraphics[width=\textwidth]{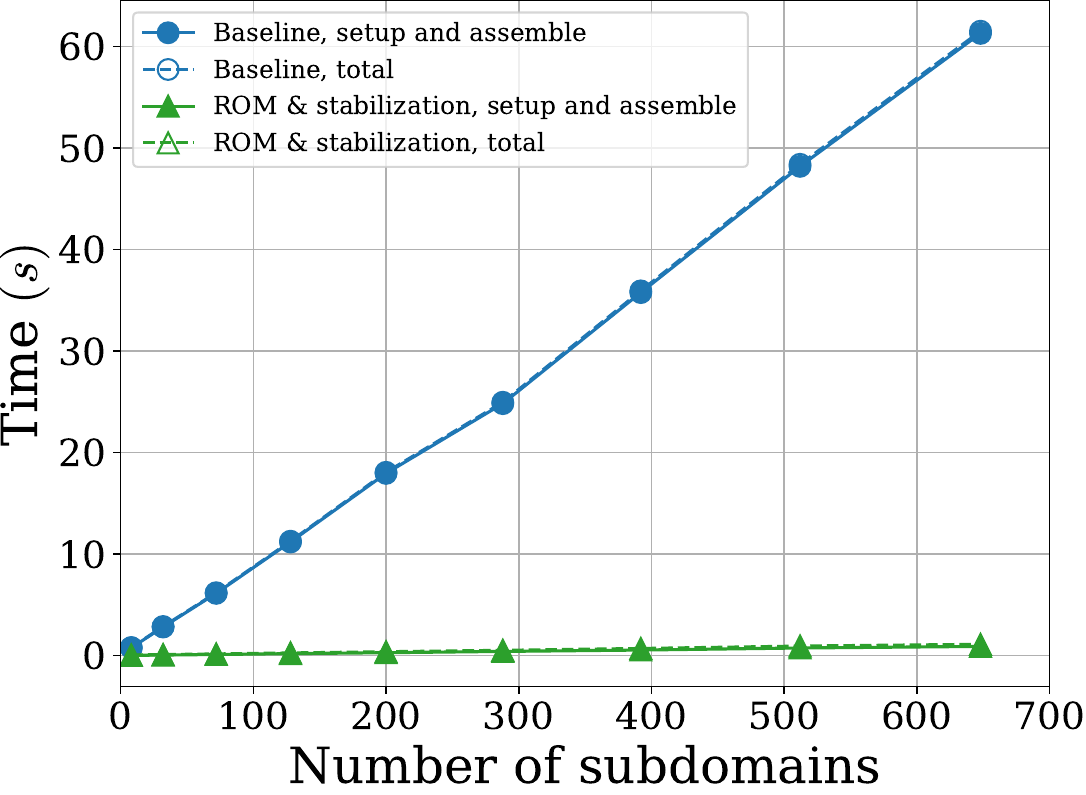}\caption{}\end{subfigure}
\caption{Iterations (a) and computational time (b) associated with the solution of the problem using a baseline BDDC approach\deleted[id=R1]{(based on full integration and without stabilization, fast assembly, or ROM)} and using the BDDC approach proposed here.}
\label{fig:test_4}
\end{figure}

We now turn to the efficiency of the proposed acceleration techniques with respect to the full integration approach. The performance is evaluated in terms of the number of iterations, setup time, and total computational time, as functions of the number of subdomains.

Figure \figref{fig:test_4}{a} shows the number of iterations for both the baseline case and the proposed method. Here, the term baseline refers to the use of full integration, i.e., no fast assembly technique and no ROM are employed. Additionally, stabilization is not included in the baseline configuration. It can be observed that there is no significant difference in the number of iterations between the baseline and the proposed method, except for the last data point, where the baseline requires one additional iteration.

In contrast, a substantial difference emerges when considering the setup and total computational times, as shown in Figure \figref{fig:test_4}{b}. In the baseline case, the majority of the computational cost is associated with the setup phase (that here includes also the assemble time), which is expected since full numerical integration must be performed. These results confirm the superiority of the proposed techniques with respect to a standard approach in terms of computational time.

\paragraph{Scalability of the method}\label{sssec: RES scalability}
\begin{figure}[ht]
\centering
\begin{subfigure}[b]{0.48\textwidth}\centering\includegraphics[width=\textwidth]{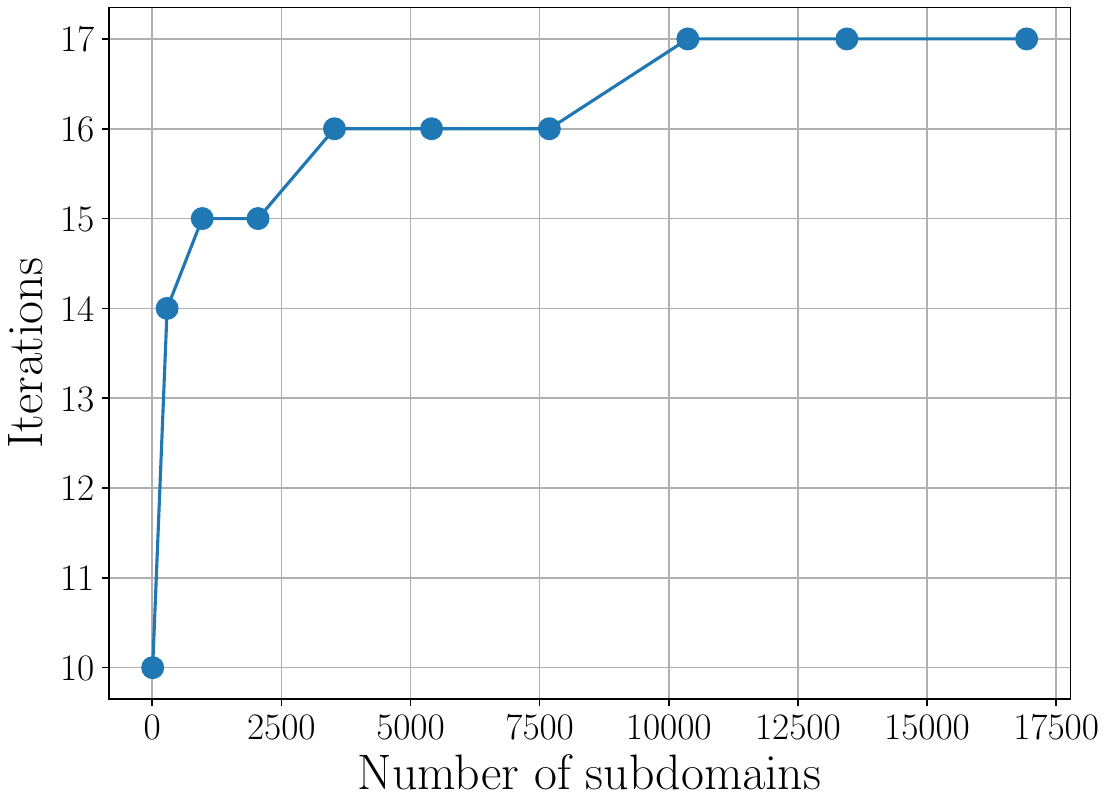}\caption{}\end{subfigure}\hfill
\begin{subfigure}[b]{0.48\textwidth}\centering\includegraphics[width=\textwidth]{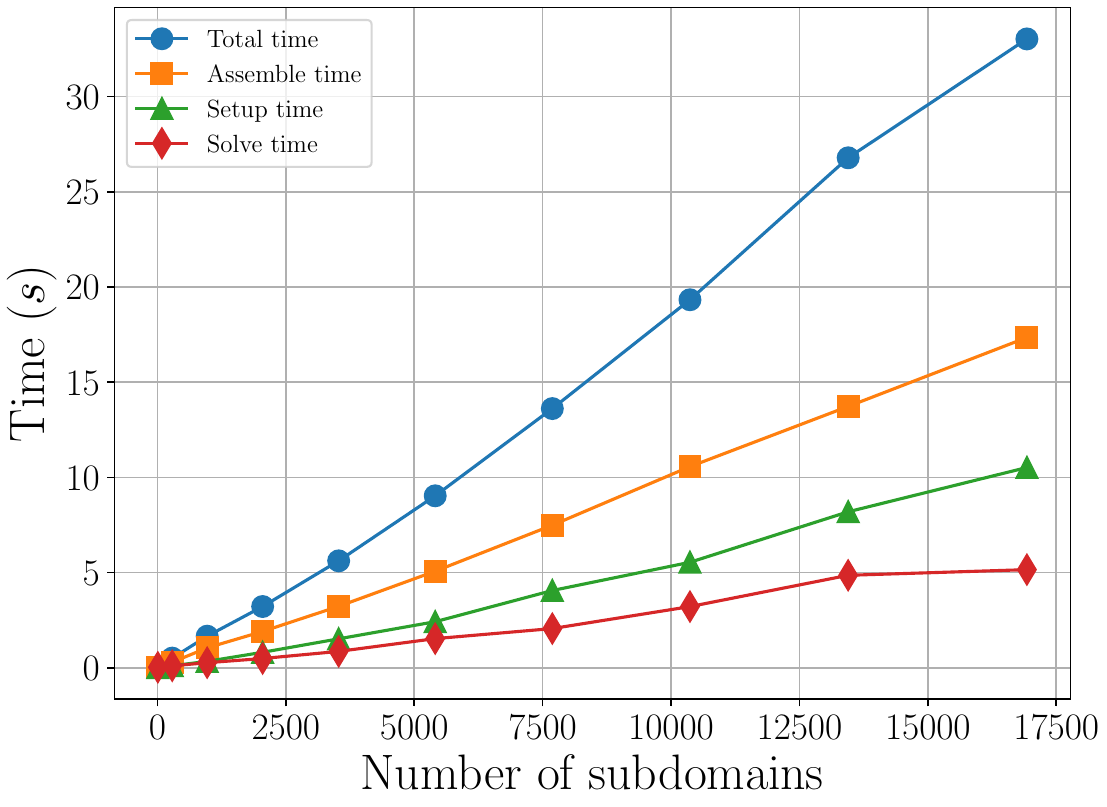}\caption{}\end{subfigure}
\caption{Iterations (a) and computational time (b) as a function of number of cells. In (b), the computational time is decomposed into assembly time, setup time, solve time, and total time (approximately equal to the sum of the former components). The test is performed up to a large number of cells, specifically 17,000.}
\label{fig:test_5}
\end{figure}

To examine scalability, the method is evaluated for increasingly large numbers of cells in the geometry, reaching up to 17,000 subdomains. From Figure \figref{fig:test_5}{a}, it can be seen that the number of iterations slightly increases with the number of cells, however reaching a constant value of 17. 

Figure \figref{fig:test_5}{b} reports the different components of the computational time, namely: the assembly time that includes the construction of the local matrices, the setup time associated with the computation of all quantities required by the BDDC method, the solve time corresponding to the iterative phase, and the total time, which is approximately the sum of the previous contributions. From this figure, it can be observed that the solve phase is the fastest component, being approximately twice as fast as the setup phase, which in turn is about twice as fast as the assembly phase. The overall computational complexity of the different components with respect to the number of cells appears to exhibit a linear trend. This is expected for problems that are not dominated by the size of the coarse correction. Remarkably, the most demanding analysis including 17,000 cells was still performed in approximately 30 seconds in a standard laptop \added[id=R2]{(online phase only; see Section~\ref{ssec: ROM cost} for the offline training cost).}

\added[id=R2]{For completeness, the problem size and statistics for the largest test are reported here. The total number of non-redundant degrees of freedom is $n_I + n_U = 2{,}169{,}728$, of which $n_U = 510{,}784$ belong to the skeleton and $n_I = 1{,}658{,}944$ are internal. The coarse problem has $n_C = 136{,}344$ degrees of freedom and is solved with a direct Cholesky factorization; its diagonally scaled condition number is $\kappa(S_C) = 6.0 \times 10^6$, which induces a relative error of at most $\kappa(S_C)\,\varepsilon \approx 10^{-9}$ in the coarse solve ($\varepsilon$ being the double-precision machine epsilon), well below the solver tolerance of $10^{-8}$, and therefore does not affect the robustness of the preconditioner. The peak resident memory is 30 GB (aggregate over 8 processes). It should be noted that this scalability test is performed on a structured quarter-annulus geometry, which is more regular than the application examples of Section~\ref{ssec: RES numerical test}; nevertheless, the latter exhibit a comparable iteration count. More challenging configurations, such as highly graded geometries or severe trimming, may increase the BDDC iteration count; their systematic assessment is left for future work.}

\paragraph{\added[id=R2]{Parallel efficiency and cost distribution}}
\added[id=R2]{So far the scalability has been assessed with respect to the number of subdomains at a fixed number of cores; we now turn to the parallel behaviour. To this end, a weak-scaling study, which is the natural measure of scalability for a domain-decomposition method, was carried out on a shared-memory server equipped with four Intel Xeon Gold 6148 sockets (20 cores per socket, 80 physical cores, four NUMA nodes, 1~TB of RAM). Each MPI rank is pinned to a distinct physical core with node-local memory, using a single computational thread per rank, and the reported times are the fastest of three repetitions}

\added[id=R2]{Table~\ref{tab: RES weak} reports the weak scaling at a fixed load of approximately $2{,}048$ subdomains per core, up to $131{,}072$ subdomains and $1.68\times10^{7}$ degrees of freedom on $64$ cores, over which the coarse problem grows proportionally from $16{,}704$ to $1{,}051{,}136$ degrees of freedom. Despite this $64\times$ growth, the BDDC iteration count stays between $15$ and $18$, the peak memory per rank remains $4.4$--$5.4$~GB, and the direct coarse solve stays below $1$~s, so the coarse problem does not become the bottleneck; this confirms the scalability of the coarse space. The parallel efficiency nonetheless decreases from $88\%$ on $2$ cores to $19\%$ on $64$ cores, and the distribution of costs reveals the origin. The coarse solve ($\le 1$~s) and the inter-process communication ($\le 17$~s) remain small throughout; the growth is concentrated in the solve phase and is dominated by the vector operations of the interface conjugate-gradient iteration, which reach about $160$~s of the $196$~s solve on $64$ cores. In the present implementation these operations are carried out redundantly on every process over the full interface vector, so that their cost scales with the global problem size; a distributed Krylov solver for the interface problem would remove them. A secondary contribution to the efficiency loss is the saturation of the per-socket memory bandwidth in the setup and assembly phases once more than about four ranks share a socket.}

\added[id=R2]{We emphasise that the parallel implementation is a research prototype and that its optimisation lies beyond the scope of this work, which concerns the overall solution strategy rather than parallel-implementation tuning. The principal route to improved efficiency at large core counts is a distributed Krylov solver for the interface problem, which would eliminate the redundant per-process vector operations identified above; at substantially larger scale a parallel or multilevel coarse solve would also become relevant. We finally observe that, in three dimensions, the per-cell workload grows substantially (the local problems scale as $(p+1)^{3}$ rather than $(p+1)^{2}$), so that each subdomain carries considerably more computation; this higher granularity is expected to improve the parallel efficiency of the dominant per-subdomain phases.}

\begin{table}[ht]\centering
\caption{\added[id=R2]{Weak scaling at a fixed load of $\approx 2{,}048$ subdomains per core on the four-socket Xeon server: number of subdomains, total (non-redundant) and coarse degrees of freedom, BDDC iteration count, the setup/assembly/solve decomposition of the wall time with the coarse-solve and communication contributions (both contained in the solve phase) listed separately, and the weak-scaling efficiency $E_w(p)=T(1)/T(p)$, expressed as a percentage.}}
\label{tab: RES weak}
{\footnotesize\setlength{\tabcolsep}{3pt}
\begin{tabular}{rrrrrrrrrrrr}\hline
Cores & Subdom. & DOFs & Coarse DOFs & Iters & Setup [s] & Asm.\ [s] & Solve [s] & Coarse [s] & Comm.\ [s] & Total [s] & Eff. \\ \hline
$1$  & $2{,}048$   & $263{,}168$      & $16{,}704$      & $15$ & $37.1$  & $17.1$  & $10.6$  & $0.09$ & $0.01$  & $64.8$  & $100\%$ \\
$2$  & $4{,}232$   & $543{,}168$      & $34{,}316$      & $16$ & $40.7$  & $19.4$  & $13.4$  & $0.14$ & $0.82$  & $73.5$  & $88\%$  \\
$4$  & $8{,}192$   & $1{,}050{,}624$  & $66{,}176$      & $16$ & $42.8$  & $19.4$  & $15.0$  & $0.17$ & $1.20$  & $77.2$  & $84\%$  \\
$8$  & $16{,}200$  & $2{,}076{,}480$  & $130{,}500$     & $17$ & $44.8$  & $20.5$  & $21.3$  & $0.23$ & $1.53$  & $86.6$  & $75\%$  \\
$16$ & $32{,}768$  & $4{,}198{,}400$  & $263{,}424$     & $17$ & $48.1$  & $24.2$  & $38.8$  & $0.30$ & $3.12$  & $111.1$ & $58\%$  \\
$32$ & $66{,}248$  & $8{,}485{,}568$  & $531{,}804$     & $17$ & $56.2$  & $30.4$  & $86.9$  & $0.52$ & $5.95$  & $173.4$ & $37\%$  \\
$64$ & $131{,}072$ & $16{,}785{,}408$ & $1{,}051{,}136$ & $18$ & $94.5$  & $46.7$  & $195.6$ & $0.99$ & $17.03$ & $336.8$ & $19\%$  \\ \hline
\end{tabular}}
\end{table}

\subsection{Application examples} \label{ssec: RES numerical test}
We now demonstrate the effectiveness of the proposed method on application-oriented examples. In particular, we consider two scenarios: a sandwich wing profile with a lattice core, and a wrench featuring a lattice structure.

\subsubsection{Sandwich wing} \label{sssec: RES wing}
\begin{figure}[ht]
\centering
\begin{subfigure}[b]{0.48\textwidth}\centering\includegraphics[width=\textwidth]{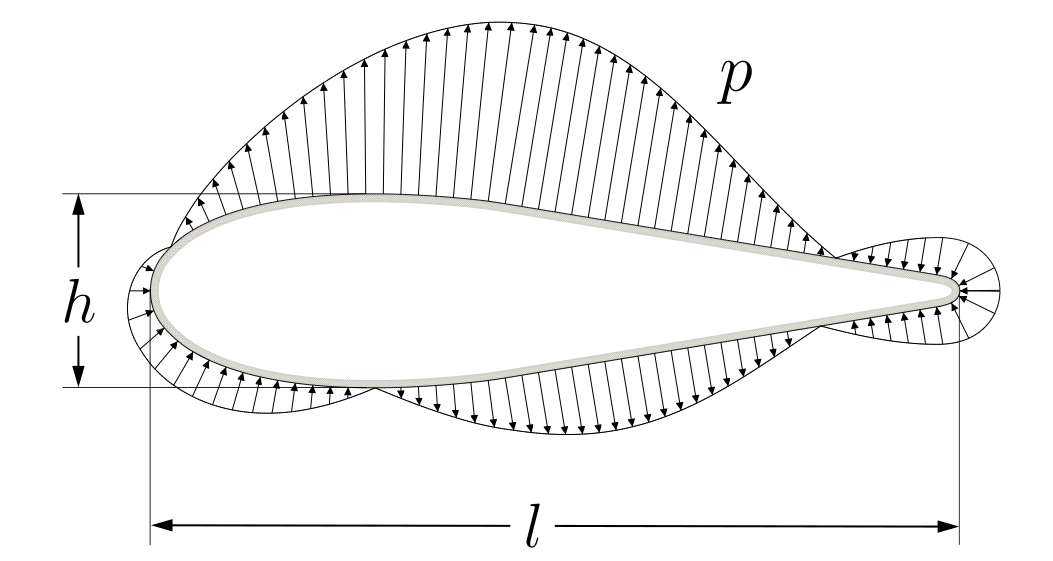}\caption{}\end{subfigure}\hfill
\begin{subfigure}[b]{0.48\textwidth}\centering\includegraphics[width=\textwidth]{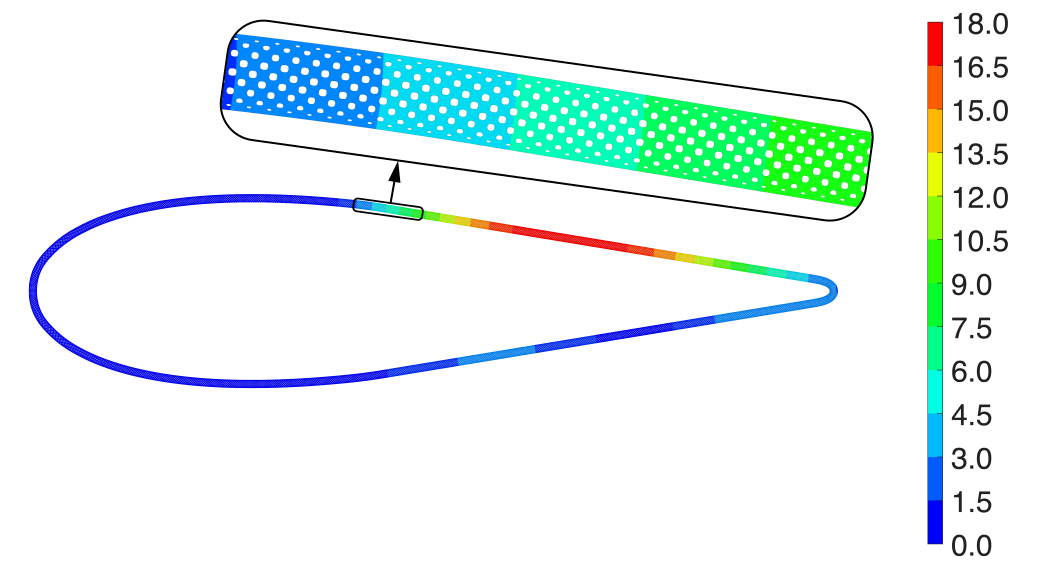}\caption{}\end{subfigure}
\caption{(a) Geometry of the sandwich wing example in Section \ref{sssec: RES wing}. (b) Contour of the magnitude of the displacement field with a zoom on the wing profile. Displacements are expressed here in mm.}
\label{fig: RES wing}
\end{figure}
The first application considers the analysis of a 2D wing profile featuring a sandwich structure, which is composed of a continuous external layer and a porous inner core. The geometry of the wing is depicted in Figure \figref{fig: RES wing}{a}, and it is characterized by a chord length of $l=10$ m, a maximum height of $h=220$ mm, and a profile thickness of $10$ mm.

The load applied to the structure, also illustrated in Figure \figref{fig: RES wing}{a}, simulates a typical aerodynamic relative pressure distribution over a wing profile with a maximum value of $p_{\mathrm{max}}=4.5$ N/mm. The structure is clamped at the nodes where the internal part of the wing section would be in contact with the main spar, representing a realistic constraint for a wing structure.

The domain is obtained by combining a total of 4,000 cells. The internal lattice structure is generated using a Schwarz Diamond level-set, whose threshold parameters are kept constant along the chord direction but are varied through the thickness within the parameter space $\mathbb{P}_2^\mathrm{SD}=[0.1,1.0]^4$. This specific selection of parameters ensures the creation of a continuous, solid external boundary, which is a critical requirement for aerodynamic performance. On the other hand, it produces a porous internal core, which is a characteristic of sandwich structures designed for lightweight applications. \deleted[id=R2]{For this level-set, the ROM employs 40 basis functions and 2 clusters per threshold parameter. The MDEIM approximation is based on Lagrangian interpolants of degree 6. The resulting fast assembly vector is characterized by a polynomial degree $q = 2$.}

The stabilization parameter for this analysis was chosen as $\rho=1\cdot10^{-4}$. It is worth mentioning that despite the geometric complexity of the problem, the BDDC solver demonstrated excellent performance, reaching convergence in just 18 iterations. Figure \figref{fig: RES wing}{b} presents a colormap of the magnitude of the displacement field superimposed, illustrating the structural response under the applied aerodynamic load. The internal structure of the sandwich wing is further illustrated in the zoomed view within the same figure.

\subsubsection{Lattice wrench} \label{sssec: RES wrench}
\begin{figure}[ht]
\centering
\begin{subfigure}[b]{0.48\textwidth}\centering\includegraphics[width=\textwidth]{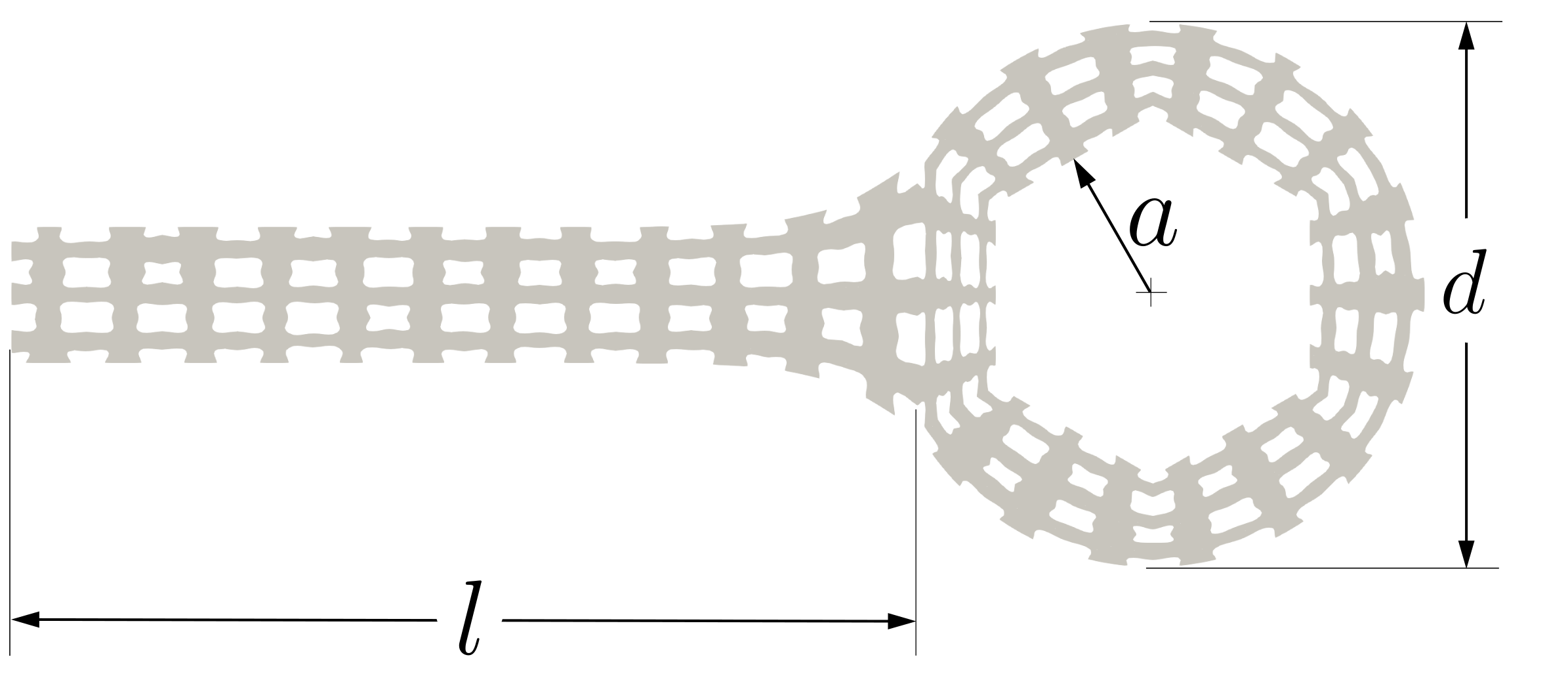}\caption{}\end{subfigure}\hfill
\begin{subfigure}[b]{0.48\textwidth}\centering\includegraphics[width=\textwidth]{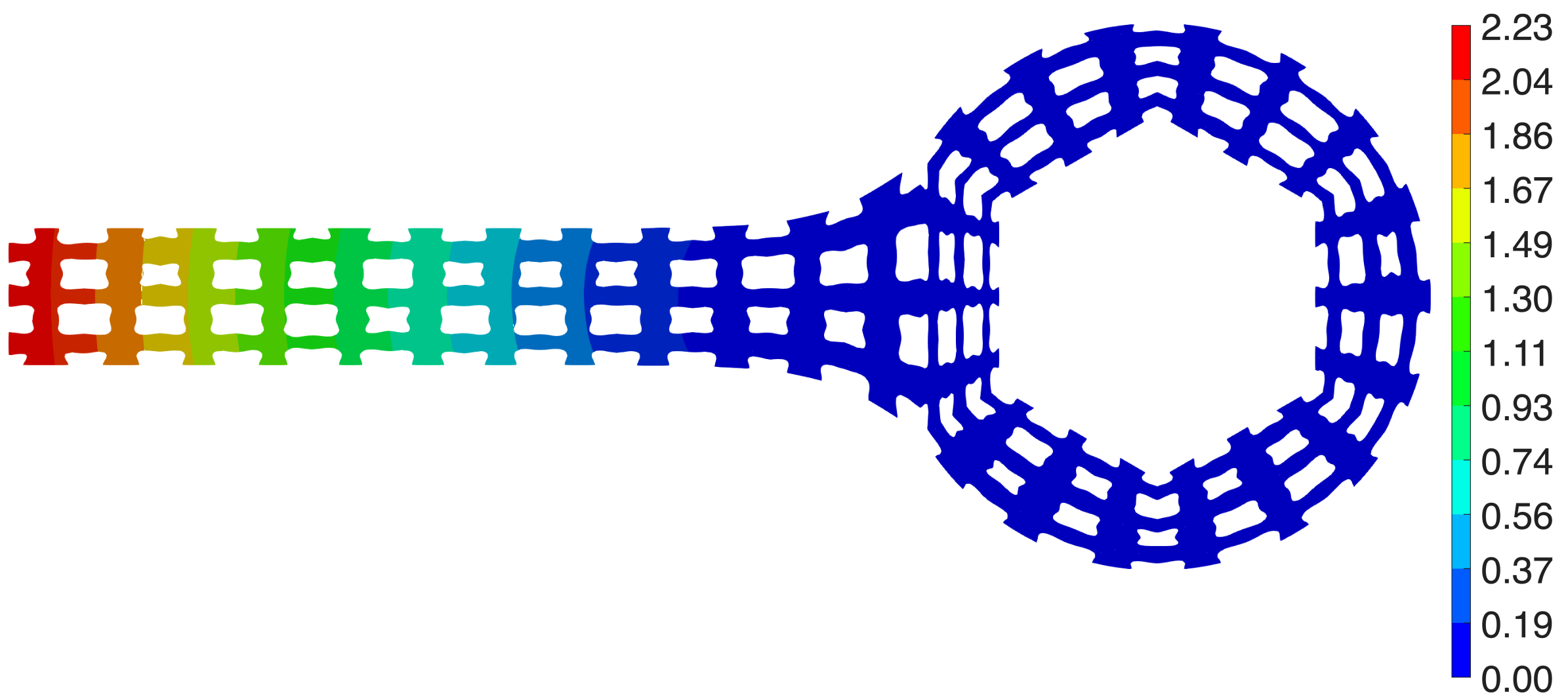}\caption{}\end{subfigure}\\
\begin{subfigure}[b]{0.48\textwidth}\centering\includegraphics[width=\textwidth]{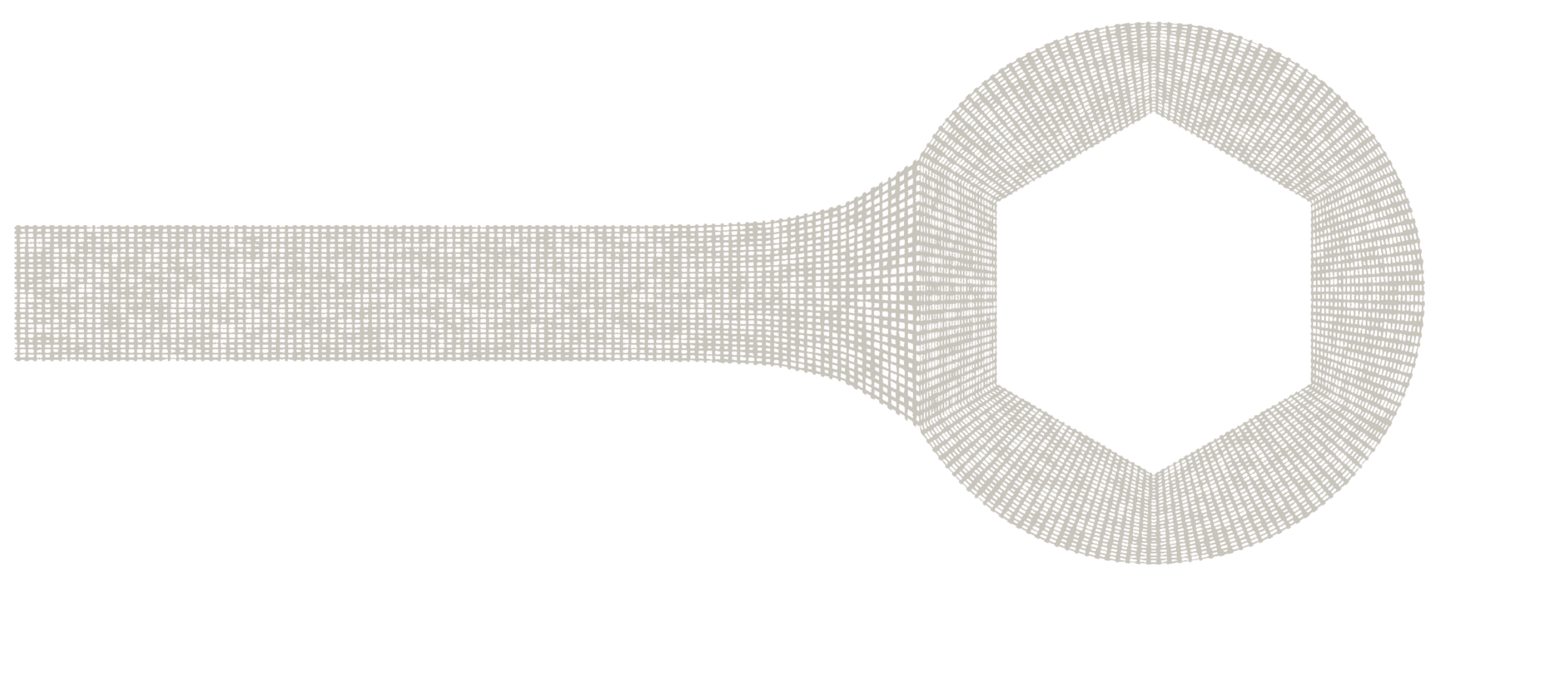}\caption{}\end{subfigure}\hfill
\begin{subfigure}[b]{0.48\textwidth}\centering\includegraphics[width=\textwidth]{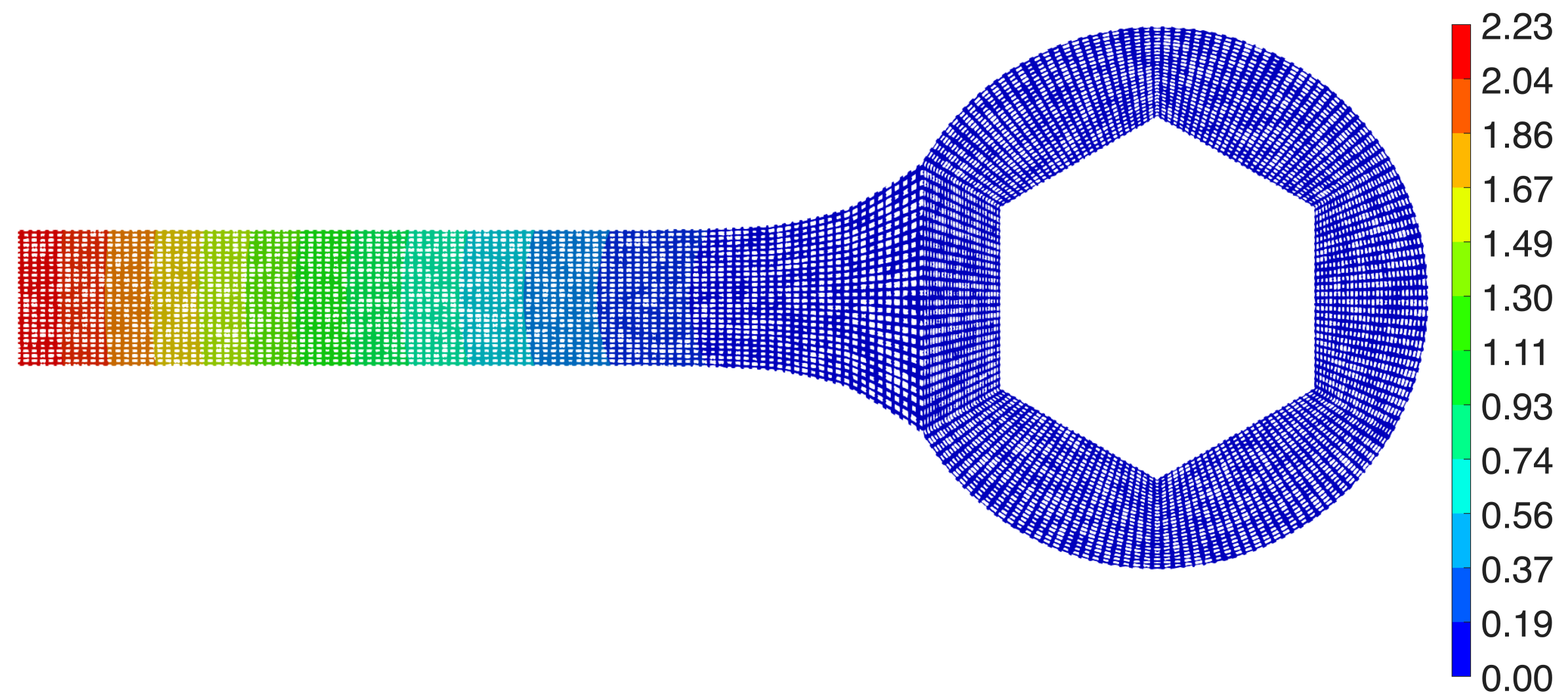}\caption{}\end{subfigure}

\caption{Geometry of lattice structure representing a wrench as described in Section \ref{sssec: RES wrench} \added[id=R2]{for the coarse (a) and the dense (c) lattice structures}. Contour of the magnitude of the displacement field  expressed in mm \added[id=R2]{for the coarse (b) and the dense (d) lattice structures}.}
\label{fig: RES wrench}
\end{figure}
In this second example, we consider a geometry that does not derive from a standard tensor-product-based disposition of cells. Specifically, we consider the wrench depicted in Figure \ref{fig: RES wrench}(a), partitioned into 90 subdomains. The global geometric features of the wrench are the handle length of $l = 50$ mm, the apothem of the hexagonal hole of $a=8.7$ mm, and the external diameter of $d = 60$ mm. The lattice micro-structure within each subdomain is defined by a Schoen IWP level-set, with threshold parameters chosen within the interval \replaced[id=R2]{$[-2.5, 2.5]^4$}{$[-2.5, 2.5]$}. \added[id=R2]{The threshold parameters are drawn randomly, independently for each cell, from a uniform distribution over this interval; the random seed is fixed in the corresponding scripts (see Table~\ref{tab: RES reproducibility}), so that the lattice geometries and the reported results are exactly reproducible.}

The mechanical behavior of the structure is modeled with a Young's modulus $E = 5$ N/mm and a Poisson's ratio $\nu = 0.25$. To simulate the action of the wrench, a homogeneous Dirichlet boundary condition is imposed on the interior hexagonal boundary of the wrench head, representing the fixed grip on a bolt. A downward traction force $\mathbf{f} = -1$ N/mm $\bm{e }_2$ is applied as a Neumann boundary condition on the vertical left side of the handle, while homogeneous Neumann boundary conditions are enforced on the remaining portions of the boundary.

\deleted[id=R2]{The fast assembly tensor employs a polynomial degree of $q=2$. Its surrogate model is constructed using 40 basis functions and 2 clusters per direction (16 in total). Within each cluster, the MDEIM interpolation uses a polynomial degree of 6. }The stabilization parameter is chosen as $\rho=5 \cdot 10^{-4}$. Despite the non-trivial layout of the subdomains and the variability of the random lattice parameters, the proposed BDDC solver proves to be highly robust, achieving convergence in only \deleted[id=R2]{16}\added[id=R2]{15} iterations. The resulting displacement field $u_h$ is illustrated in Figure \ref{fig: RES wrench}(b), exhibiting the expected bending behavior of the loaded wrench.

\added[id=R2]{To assess the accuracy of the computed solution, a reference solution is built by discretizing every cell of the wrench with a grid of $k \times k$ rectangular elements of degree $p=3$ within a CutFEM framework \cite{burman2015,burman2025cut}, with $k$ up to $6$. The level-set of each element describes the corresponding portion of the cell level-set, so that the lattice geometry is exactly the same for every $k$. These reference solutions employ full integration, no ROM acceleration, and a direct solver. Two quantities of interest are monitored: the compliance and the mean displacement of the loaded boundary. Since the traction is interpolated at the degrees of freedom of the loaded boundary, the resultant of the discrete load varies slightly with the discretization. For this reason, the compliance is normalized by the square of the resultant of the discrete load, and the mean displacement by the resultant itself.}

\begin{figure}[ht]
\centering
\begin{subfigure}[b]{0.70\textwidth}\centering\includegraphics[width=\textwidth]{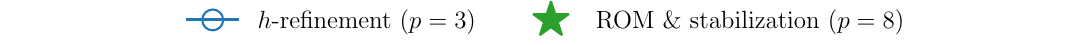}\end{subfigure}\\
\begin{subfigure}[b]{0.44\textwidth}\centering\includegraphics[width=\textwidth]{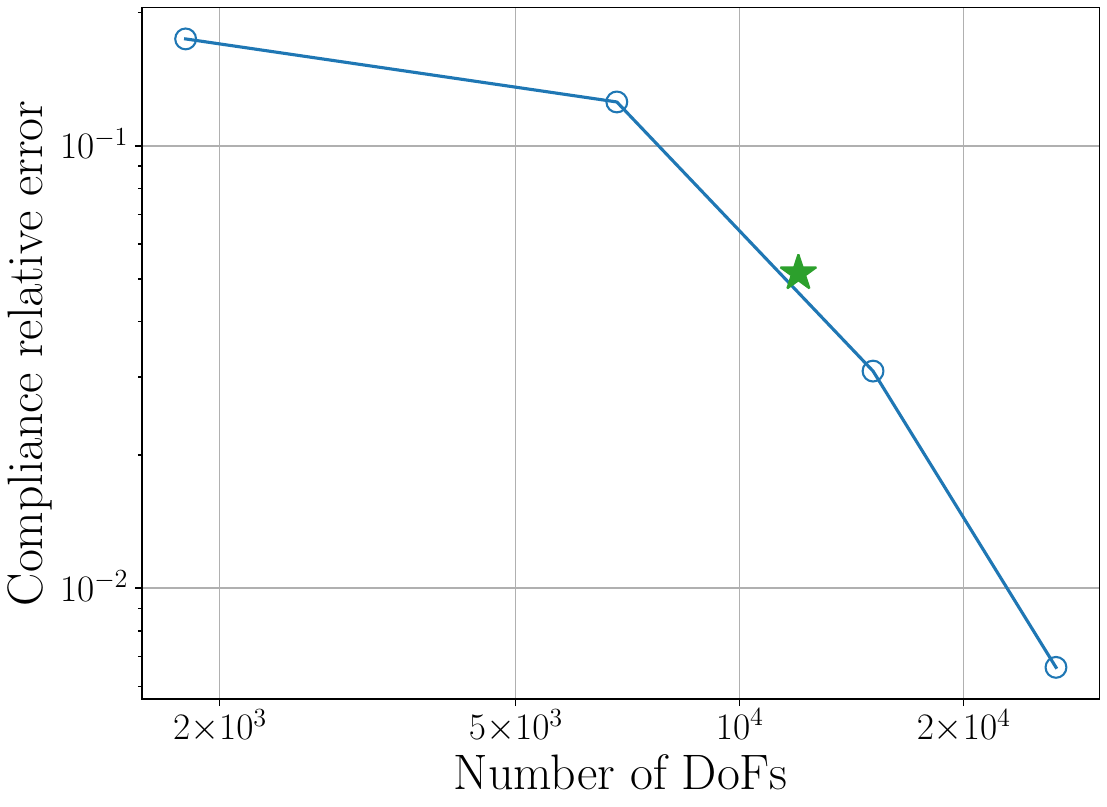}\caption{}\end{subfigure}\hfill
\begin{subfigure}[b]{0.44\textwidth}\centering\includegraphics[width=\textwidth]{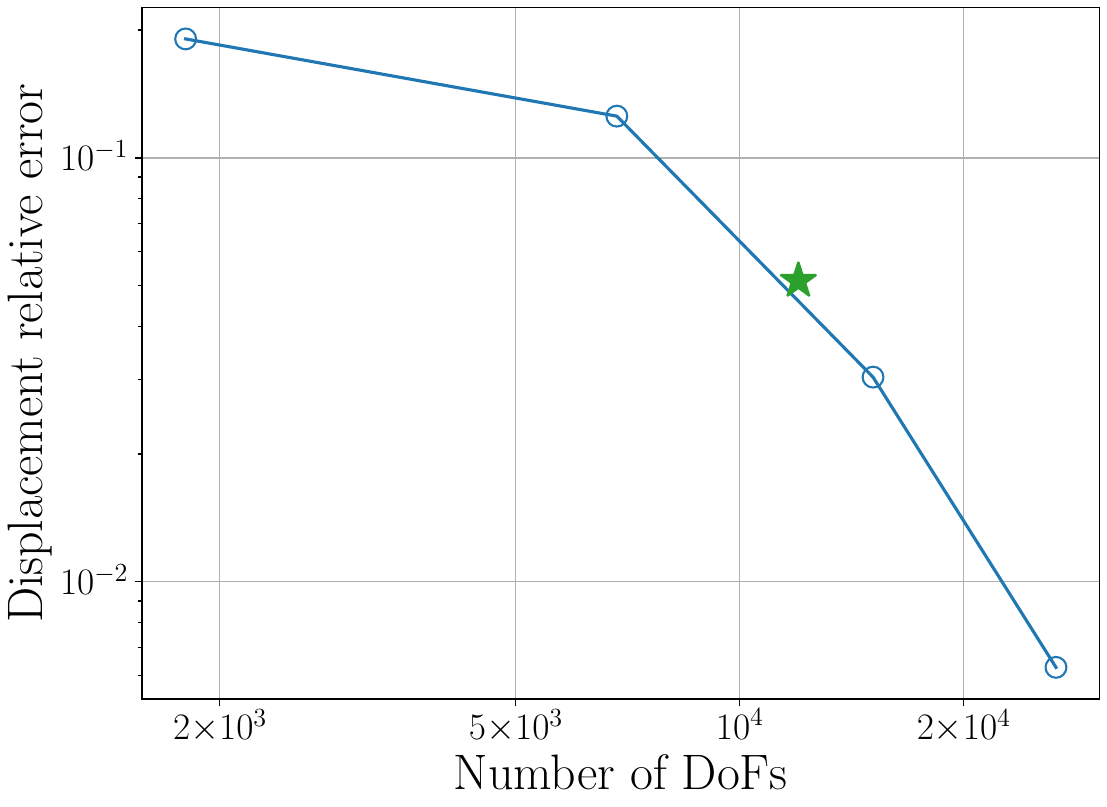}\caption{}\end{subfigure}
\caption{\added[id=R2]{Relative error of the normalized compliance (a) and of the normalized mean displacement of the loaded boundary (b) for the wrench example. In each plot, the blue curve corresponds to the $h$-refined solutions, obtained with a grid of $k \times k$ rectangular CutFEM elements of degree $p=3$ per cell, with $k=1,2,3,4$. The star corresponds to the proposed spectral method, with a single element of degree $p=8$ per cell, combined with the ROM acceleration, the $\alpha$-stabilization, and the BDDC solver. The errors are measured against the finest $h$-refined solution ($k=6$, corresponding to 59,400 DoFs).}}
\label{fig: RES wrench validation}
\end{figure}

\added[id=R2]{Figure~\ref{fig: RES wrench validation} reports the relative error of both quantities with respect to the finest reference solution, which counts 59,400 DoFs. The reference solutions converge monotonically as $k$ increases, with errors decreasing from about $17\%$ for $k=1$ down to $0.7\%$ for $k=4$. The result of the proposed method, which combines the ROM acceleration, the $\alpha$-stabilization, and the BDDC solver with degree $p=8$, lies on the same error curve, with an error close to $5\%$ for 12,000 DoFs. The impact of the individual ingredients of the method on both quantities is small compared to the discretization error. Replacing the ROM by full integration changes both quantities by $1.4\cdot10^{-4}$; removing the $\alpha$-stabilization changes them by $5.5\cdot10^{-3}$; and replacing the BDDC solver by a direct factorization changes them by $2.5\cdot10^{-10}$.}

\added[id=R2]{To further confirm the scalability observed in Section~\ref{sssec: RES scalability}, the similar wrench geometry shown in Figure~\ref{fig: RES wrench}(c) is considered, this time obtained as a lattice structure consisting of 7,290 subdomains, while all other properties remain unchanged. The BDDC solver remains robust and converges in a comparable number of iterations (18). The resulting displacement field is shown in Figure~\ref{fig: RES wrench}(d).}

\section{Conclusions} \label{sec: CONC}

In this work, we presented a novel domain decomposition method for fast simulation of large \added[id=R2]{two-dimensional} lattice structures geometrically described through implicit functions, without relying on homogenization or multiscale approaches.
Individual lattice cells, described through level-set functions that may change from cell to cell, are mapped using arbitrary order mappings (in a finite element spirit), allowing the creation of complex graded designs.

In order to solve such designs efficiently and with low memory requirements, without the need for homogenization or multiscale techniques, which rely on assumptions such as scale separation and periodicity, we use domain decomposition methods.
Specifically, the proposed framework integrates a Balanced Domain Decomposition by Constraints (BDDC) solver with high-order unfitted discretizations, which naturally handle level-set geometric descriptions.
This strategy exploits the geometric similarities between lattice cells by building a reduced order model (ROM) surrogate that assembles the cell stiffness matrices in a fraction of the time and without the need for expensive quadrature rules.
The different hypotheses and approximations introduced in the solver have been assessed through numerical experiments to validate the proposed framework.

This methodology is complemented with a $\alpha$-stabilization term, similar to the one applied in finite cell methods, required by the unfitted nature of the discretization.
When using ROM techniques, such stabilization is required to obtain a scalable solver with respect to the number of subdomains.
However, it breaks the consistency of the problem, introducing an error in the solution which remains small for moderate stabilization parameters.
As reported in the numerical experiments, the number of iterations remains asymptotically bounded as the ratio $H/h$ (subdomain versus mesh size) is kept constant, in agreement with the scalability properties of BDDC methods.
We demonstrated its performance on a\deleted[id=R2]{complex} 2D problem with 17,000 varying cell geometries, solving it in approximately 30 seconds on an off-the-shelf laptop. \added[id=R2]{A dedicated weak-scaling study further confirms this behaviour, the BDDC iteration count remaining bounded (between $15$ and $18$) as the problem is scaled up to $131{,}072$ subdomains and $1.68\times10^{7}$ degrees of freedom on $64$ cores.}

Such a solver could constitute a solid foundation for use in conjunction with optimization algorithms: it allows forward problems to be solved quickly, and the ROM surrogate for stiffness matrix assembly makes it straightforward to compute the gradient of the solution with respect to the design parameters.

We will also investigate alternative ROM surrogates for 3D problems.
The MDEIM coefficient interpolation algorithm proposed in this work relies on one parameter coefficient per interpolation mesh vertex, which in hexahedral meshes amounts to at least 8 coefficients per cell.
\added[id=R2]{In view of the additional complexity associated with a three-dimensional extension of the present framework, two possible research directions can be considered}: the use of neural network based surrogates that do not suffer the curse of dimensionality; and the discretization of the porosity field on tetrahedral meshes, which have fewer vertices per cell and thus lead to lower-dimensional parameter spaces.
Beyond the construction of the surrogates themselves, we plan to explore their use as approximate subdomain inverses for building efficient preconditioners. In this setting, the ROM-generated subdomain stiffness matrices would be applied directly inside a Krylov solver without the need to invert them, which would also reduce the need for the stabilization introduced in this work.
Finally, we plan to extend this framework to \added[id=R2]{support more complex boundary conditions, to} the case of hyperelastic materials under large deformations, and to the simulation of plates.

\section*{Acknowledgments}
The authors acknowledge the financial support of the Swiss National Science Foundation through
the project FLAS$_h$ (200021\_214987). The authors would like to acknowledge the support of Raul Rubio for his valuable scientific contributions, particularly in the area of domain decomposition.

\paragraph*{Declaration of generative AI and AI-assisted technologies in the writing process:}
During the preparation of this work, the authors used ChatGPT, Gemini, and Claude in order to improve language and readability.
After using this tool/service, the authors reviewed and edited the content as needed and take full responsibility for the content of the publication.

\appendix
\section{Algorithms}
In this appendix, some of the algorithms employed in the method are presented to provide additional clarity and to assist the reader in navigating its various aspects.

\begin{algorithm}[H]
\caption{Application of operator $S_{BDDC}^{-1}$ to $v$ (parallelizable over subdomains $i$).}\label{alg: SBDDC multiplication}
\begin{algorithmic}[1]
\State \textbf{Input:} vector $v \in \mathbb{R}^{n_U}$
\State \textbf{Output:} vector $w \in \mathbb{R}^{n_U}$
\State $w \gets 0 \in \mathbb{R}^{n_U}$
\ForAll{$i\in[1,\dots,n_c]$}\Comment{Executed in parallel}
    \State Extract subvector $v^{(i)}$       from $v$ 
    \State Solve         $\begin{bmatrix}S^{(i)}&C^{(i)^\top}\\C^{(i)}&0\end{bmatrix}\begin{bmatrix}x^{(i)}\\\lambda^{(i)}\end{bmatrix}=\begin{bmatrix}D^{(i)} v^{(i)}\\0\end{bmatrix}$ \Comment{Using Algorithm \ref{alg: Local solves a}}
    \State Prolong $D^{(i)} x^{(i)}$ to $x$
    \State Compute       $w = w + x$
\EndFor
\State $z \gets 0 \in \mathbb{R}^{n_c}$
\ForAll{$i\in[1,\dots,n_c]$}\Comment{Executed in parallel}
    \State Extract subvector $v^{(i)}$       from $v$ 
    \State Prolong ${\Psi^{(i)\top}} D^{(i)\top}v^{(i)}$ to $y$
    \State Compute       $z = z + y$
\EndFor
\State Solve $S_C g = z$\Comment{Executed with a parallel sparse Cholesky solver}
\ForAll{$i\in[1,\dots,n_c]$}\Comment{Executed in parallel}
    \State Extract subvector $g^{(i)}$       from $g$ 
    \State Prolong $D^{(i)}\Psi^{(i)} g^{(i)}$ to $x$
    \State Compute       $w = w + x$
\EndFor
\State \Return $w$
\end{algorithmic}
\end{algorithm}

\begin{algorithm}[H]
\caption{Local solve of the problem $\begin{bmatrix}S^{(i)}&C^{(i)^\top}\\C^{(i)}&0\end{bmatrix}\begin{bmatrix}x^{(i)}\\\lambda^{(i)}\end{bmatrix}=\begin{bmatrix}f^{(i)}\\g^{(i)}\end{bmatrix}$}
\label{alg: Local solves a}
\begin{algorithmic}[1]
\State \textbf{Input:} vector $f^{(i)} \in \mathbb{R}^{n_U^{(i)}}$, vector $g^{(i)} \in \mathbb{R}^{n_C^{(i)}}$
\State \textbf{Output:} vector $x^{(i)} \in \mathbb{R}^{n_U^{(i)}}$, vector $\lambda^{(i)} \in \mathbb{R}^{n_C^{(i)}}$.
\State Remove the rows and columns associated with cross-point coarse DoFs in $C^{(i)}$, $x^{(i)}$, $\lambda^{(i)}$, $f^{(i)}$, $g^{(i)}$, and $S^{(i)}$, yielding the system
\begin{equation*}
\begin{bmatrix}\hat{S}^{(i)}&\hat{C}^{(i)^\top}\\\hat{C}^{(i)}&0\end{bmatrix}\begin{bmatrix}\hat{x}^{(i)}\\\hat{\lambda}^{(i)}\end{bmatrix}=\begin{bmatrix}\hat{f}^{(i)}\\\hat{g}^{(i)}\end{bmatrix}.
\end{equation*}
\State Compute $\hat{S}^{(i)-1}\hat{C}^{(i)^\top}$ and $\hat{T}^{(i)}=\hat{C}^{(i)}\hat{S}^{(i)-1}\hat{C}^{(i)^\top}$
\Comment{Through solve of $\hat{S}^{(i)}X=\hat{C}^{(i)^\top}$}
\State Compute $\hat{r}^{(i)}=\hat{C}^{(i)}\hat{S}^{(i)-1}\hat{f}^{(i)}-\hat{g}^{(i)}$
\Comment{Through solve of $\hat{S}^{(i)}x=\hat{f}^{(i)}$}
\State Solve $\hat{T}^{(i)}\hat{\lambda}^{(i)}=\hat{r}^{(i)}$
\State Compute $\hat{h}^{(i)}=\hat{f}^{(i)}-\hat{C}^{(i)^\top}\hat{\lambda}^{(i)}$
\State Compute $\hat{x}^{(i)}=\hat{S}^{(i)-1}\hat{h}^{(i)}$ 
\Comment{Through solve of $\hat{S}^{(i)}x=\hat{h}^{(i)}$}
\State Assemble $x^{(i)}$  and $\lambda^{(i)}$  by reversing step 3
\State \Return $x^{(i)}$, $\lambda^{(i)}$
\end{algorithmic}
\end{algorithm}

\begin{algorithm}[H]
\caption{Select magic points}\label{alg: MDEIM magic points}
\begin{algorithmic}[1]
\State \textbf{Input:} matrix $U_r \in \mathbb{R}^{n_i\times n_r}$
\State \textbf{Output:} set of indices  $\{m_1,\dots,m_{n_r}\} \in \mathbb{N}$
\For{$j = 1$ to $n_r$}
    \If{$j = 1$}
        \State $m_1 \gets \operatorname{arg\,max}_i \left| U_{r\,1}^{[i]} \right|$
    \Else
        \State Compute projection of $U_{r\,j}$ onto $\mathrm{span}\{U_{r\,1},\dots,U_{r\,j-1}\}: U_{p\,j} = \sum_{k=1}^{j-1} U_{r\,k} U_{r\,k}^\top \, U_{r\,j}$
        \State $m_j \gets \operatorname{arg\,max}_i \left|U_{r\,j}^{[i]}-U_{p\,j}^{[i]}\right|$
    \EndIf
\EndFor
\State \textbf{return} $m_1, m_2, \dots, m_{n_r}$
\end{algorithmic}
\end{algorithm}

\begin{algorithm}[H]
\caption{Assembly of the local stiffness matrix}
\label{alg:local_stiffness}
\begin{algorithmic}[1]
\State \textbf{Input:} tensor function $ \hat{C}_{i_1 j_1 i_2 j_2}(\boldsymbol{\xi})$, threshold vector $\boldsymbol{\mu}_0$, stabilization parameter $\rho$
\State \textbf{Output:} stiffness tensor $K_{i_1i_2k_1k_2}$
\State Compute $\mathrm{C}_{i_1 j_1 i_2 j_2 k_3}$ by interpolating $ \hat{C}_{i_1 j_1 i_2 j_2}(\boldsymbol{\xi})$
\State Obtain $I( \boldsymbol{\mu}_0)$ \Comment{Using Algorithm \ref{alg: mdeim}}
\State Obtain $\mathcal{I}^{k_1k_2k_3}_{j_1j_2} \gets \text{Reshape}(I)$ 
\State Obtain $\tilde{\mathcal{I}}^{k_1k_2k_3}_{j_1j_2} \gets \int_{\hat{\Pi}} \fpd{\mathscr{B}_{k_1}}{\xi_{j_1}}
     \fpd{\mathscr{B}_{k_2}}{\xi_{j_2}}\mathscr{L}_{k_3}
     \,\dd\hat{\Omega} - \mathcal{I}^{k_1k_2k_3}_{j_1j_2}$
     \Comment{$\int_{\hat{\Pi}} \fpd{\mathscr{B}_{k_1}}{\xi_{j_1}}
     \fpd{\mathscr{B}_{k_2}}{\xi_{j_2}}\mathscr{L}_{k_3}
     \,\dd\hat{\Omega}$ is precomputed} 
\State Obtain $K_{i_1i_2k_1k_2} \gets \left(\mathcal{I}^{k_1k_2k_3}_{j_1j_2} + \rho \cdot \tilde{\mathcal{I}}^{k_1k_2k_3}_{j_1j_2}\right) \mathrm{C}_{i_1 j_1 i_2 j_2 k_3}$ 
\State \Return $K_{i_1i_2k_1k_2}$
\end{algorithmic}
\end{algorithm}

\begin{algorithm}[H]
\caption{MDEIM}
\label{alg: mdeim}
\begin{algorithmic}[1]
\State \textbf{Input:} list(tensor $U_r \in \mathbb{R}^{n_i\times n_r}$), list(reduced-order functions $I_r(\bm{\mu}): \mathbb{R}^{4} \to \mathbb{R}^{r}$), threshold vector $\boldsymbol{\mu}_0$
\State \textbf{Output:} vector $I \in \mathbb{R}^{n_i}$
\State $i_k \gets \text{FindCluster}(\boldsymbol{\mu}_0)$
\State $I \gets U_r^{i_k} I_r^{i_k}(\boldsymbol{\mu}_0)$
\State \Return $I$
\end{algorithmic}
\end{algorithm}

\bibliography{Bibliography}
\end{document}